\documentclass{article}

\usepackage{amsfonts,amsmath,tikz,xcolor,mathtools,extarrows, amsthm}
\usepackage[margin=1in]{geometry}
\usepackage[colorlinks=true,linkcolor=red!70!black,citecolor=green!70!black,urlcolor=magenta!70!black,backref]{hyperref}

\colorlet{darkblue}{blue!70!black}
\colorlet{darkred}{red!70!black}
\colorlet{darkgreen}{green!70!black}
\colorlet{darkmagenta}{magenta!70!black}

\newenvironment{tric}
    {\begin{tikzpicture}[scale= 0.867,semithick,draw=darkblue,double distance=1.1,
        baseline={([yshift=-.8ex]current bounding box.center)}] }
    {\end{tikzpicture}}

\newenvironment{trics}
    {\begin{tikzpicture}[scale= 0.5,semithick,draw=darkblue,double distance=1.1,
        baseline={([yshift=-.8ex]current bounding box.center)}] }
    {\end{tikzpicture}}

\newtheorem{thm}{Theorem}[section]
\newtheorem{lemma}[thm]{Lemma}
\newtheorem{proposition}[thm]{Proposition}
\newtheorem{corollary}[thm]{Corollary}
\newtheorem{cor}[thm]{Corollary}
\newtheorem{conjecture}[thm]{Conjecture}
\newtheorem{conj}[thm]{Conjecture}

\newtheorem{notation}[thm]{Notation}
\newtheorem{example}[thm]{Example}

\theoremstyle{definition}
\newtheorem{definition}[thm]{Definition}
\newtheorem{defn}[thm]{Definition}

\theoremstyle{remark}
\newtheorem{remark}[thm]{Remark}
\newtheorem{rmk}[thm]{Remark}

\DeclareMathOperator{\im}{im}
\DeclareMathOperator{\id}{id}
\DeclareMathOperator{\ot}{\otimes}

\DeclareMathOperator{\wt}{\rm{wt}}
\DeclareMathOperator{\Uq}{U_q(\mathfrak{g}_2)}
\DeclareMathOperator{\unw}{\underline{w}}
\DeclareMathOperator{\unu}{\underline{u}}
\DeclareMathOperator{\unv}{\underline{v}}

\DeclareMathOperator{\unx}{\underline{x}}
\DeclareMathOperator{\uny}{\underline{y}}
\DeclareMathOperator{\unz}{\underline{z}}

\DeclareMathOperator{\Fund}{\bf{Fund}}
\DeclareMathOperator{\Rep}{\bf{Rep}}
\DeclareMathOperator{\Hom}{\rm{Hom}}
\DeclareMathOperator{\End}{\rm{End}}
\DeclareMathOperator{\Kar}{\rm{Kar}}

\DeclareMathOperator{\DD}{\textbf{Web}_q(\mathfrak{g}_2)}
\DeclareMathOperator{\smd}{\overset{\oplus}{\subset}}
\DeclareMathOperator{\Xa}
{\mathsf{H}_{\varpi_1\varpi_2}^{\varpi_2\varpi_1}}
\DeclareMathOperator{\Xb}
{\mathsf{H}_{\varpi_2\varpi_1}^{\varpi_1\varpi_2}}
\DeclareMathOperator{\triv}{\mathbf{1}}
\DeclareMathOperator{\Kab}{\text{Ker}_{(a, b)}}

\author{Bodish, Elijah\\
  \texttt{ebodish@mit.edu}
  \and
  Wu, Haihan\\
  \texttt{hihwu@ucdavis.edu}}
\date{}
\title{Triple clasp formulas for $G_2$}
    
\begin{document}

\maketitle
\begin{abstract}
We use Kuperberg's diagrammatic description of homomorphisms between fundamental representations of $G_2$ to give explicit recursive formulas for the idempotent projecting to the highest weight irreducible summand in each tensor product of fundamental representations. 
\end{abstract}
\tableofcontents

\section{Introduction}
\label{intro}

\subsection{History of Clasp Formulas}
\label{history}

\def\TLCa{
\begin{tric}
\draw (0.75,0.2)--(0.75,1.5) node[midway,left,scale=0.7,black]{$n$}
      (0.75,-0.2)--(0.75,-1.5) node[midway,left,scale=0.7,black]{$n$};
\draw[darkred,thick] (0,-0.2)rectangle(1.5,0.2)  
                     node[pos=0.5,scale=0.7,black]{$n$};
\end{tric}
}

\def\TLCb{
\begin{tric}
\draw (0.6,0.2)--(0.6,1.5) node[midway,left,scale=0.7,black]{$n-1$}
      (0.6,-0.2)--(0.6,-1.5) node[midway,left,scale=0.7,black]{$n-1$};
\draw (1.6,-1.5)--(1.6,1.5);
\draw[darkred,thick] (0,-0.2)rectangle(1.2,0.2)  
                     node[pos=0.5,scale=0.7,black]{$n-1$};
\end{tric}
}

\def\TLCc{
\begin{tric}
\draw (0.4,0.9)--(0.4,1.5) node[midway,left,scale=0.7,black]{$n-1$}
      (0.4,-0.9)--(0.4,-1.5) node[midway,left,scale=0.7,black]{$n-1$}
      (0.4,0.5)--(0.4,-0.5) node[midway,left,scale=0.7,black]{$n-2$};
\draw (0.9,0.5)..controls(0.9,0.1)and(1.6,-0.3)..(1.6,1.5)
      (0.9,-0.5)..controls(0.9,-0.1)and(1.6,0.3)..(1.6,-1.5);
\draw[darkred,thick] (0,0.5)rectangle(1.2,0.9)  
                     node[pos=0.5,scale=0.7,black]{$n-1$}
                     (0,-0.5)rectangle(1.2,-0.9)  
                     node[pos=0.5,scale=0.7,black]{$n-1$};
\end{tric}
}

The discovery of the Jones polynomial in the early 1980's \cite{JonesPolynomial} triggered mathematical developments in areas including knot theory and quantum algebra. One way to define the Jones polynomial is by using the braiding in the Temperley-Lieb Category \cite[Chapter 2]{Kauffman-Lins}. Half a decade earlier Rummer-Teller-Weyl found a description of morphisms between tensor products of the vector representation of $SL_2(\mathbb{C})$ in terms of cup and cap diagrams \cite[Equation 3]{Weyl1932}. The $q$-analogue of their result is that the Temperley-Lieb Category is monoidally equivalent to the full monoidal subcategory of $\Rep(U_q(\mathfrak{sl}_2))$ generated by the $q$-analogue of the vector representation.

Let $V\in \Rep(U_q(\mathfrak{sl}_2))$ denote the $q$-analogue of the vector representation of $SL_2(\mathbb{C})$. For each $n\in \mathbb{Z}_{\ge 0}$, there is an irreducible representation $V(n)$, which is a direct summand of $V^{\otimes n}$ and which is not a direct summand of $V^{\otimes m}$ for $m< n$. Note that $V(1)\cong V$. So for each $n$, there is an idempotent in the Temperley-Lieb category which can be viewed as the idempotent in $\End_{U_q(\mathfrak{sl}_2)}(V^{\otimes n})$ with image $V(n)$. The condition that $V(n)$ is not a summand of $V^{\otimes m}$ for $m< n$ implies that composing a projector with any cap diagram will result in zero. 

These idempotents are usually called Jones-Wenzl projectors, as they were first considered by Jones \cite[Section 4.2]{Jon2}, and the following explicit inductive formula was first given by Wenzl \cite{Wenzl}. 
\begin{equation}\label{JonesWenzl}
\TLCa \ =\TLCb \ \ \ +\ \ \frac{[n-1]}{[n]}\ \TLCc
\end{equation}
Here we use the notation $[m]$ to denote the quantum integer $[m]_q:=\frac{q^m-q^{-m}}{q-q^{-1}}$, for each $m\in \mathbb{Z}$. Our convention is that a red box with label $n$ is the morphism in the Temperley-Lieb category which corresponds to the idempotent with image $V(n)$. 

The Jones-Wenzl projectors and the recursive formula in Equation \eqref{JonesWenzl} describing them have been proven useful in link homology \cite{MR2901969}, Soergel bimodules \cite{MR2873427}, and the theory of subfactors and planar algebras \cite{MR2679382}. The present work is concerned with generalizing Equation \eqref{JonesWenzl} from $\mathfrak{sl}_2$ to the Lie algebra $\mathfrak{g}_2$. However, many things we say in the introduction make sense for all semisimple Lie algebras. 

Fix a finite dimensional semisimple Lie algebra $\mathfrak{g}$. There is an associated quantum enveloping algebra $U_q(\mathfrak{g})$, which is a $\mathbb{C}(q)$ algebra defined by generators and relations which ``quantize" the Serre presentation of the usual enveloping algebra \cite[Chapter 4]{JantzenQgps}. The finite dimensional irreducible type-$\textbf{1}$ representations\footnote{This means that for all simple roots $\alpha$, the element $K_{\alpha}$ acts on the $\mu$ weight space of any representation by $+ q^{(\alpha, \mu)}$ \cite[Section 5.2]{JantzenQgps}. We will only consider type-$\textbf{1}$ representations in this paper.} are in bijection with the finite dimensional irreducible representations of $\mathfrak{g}$, i.e. for each dominant integral weight $\lambda$ there is a finite dimensional irreducible module of $U_q(\mathfrak{g})$, which we denote by $V(\lambda)$. We will abuse notation and write $\Rep(U_q(\mathfrak{g}))$ to refer to the category of finite dimensional type-$\textbf{1}$ representations of $U_q(\mathfrak{g})$. The algebra $U_q(\mathfrak{g})$ is a Hopf algebra \cite[Section 4.8]{JantzenQgps}, and it turns out that $\Rep(U_q(\mathfrak{g}))$ is closed under taking tensor products. Furthermore, since we are working over $\mathbb{C}(q)$, where $q$ is an indeterminant or a generic element of $\mathbb{C}$, the category $\Rep(U_q(\mathfrak{g}))$ is a semisimple tensor category, and the Grothendieck ring of $\Rep(U_q(\mathfrak{g}))$ is isomorphic to the Grothendieck ring of the category of finite dimensional representations of $\mathfrak{g}$. 

The recursion in Equation \eqref{JonesWenzl} is expressed in terms of the Temperley-Lieb category, which describes the full monoidal subcategory of $\Rep(U_q(\mathfrak{sl}_2))$ generated by $V$. One way to generalize this subcategory to arbitrary quantum groups is proposed in the following definition.

\begin{defn}
The category $\Fund(U_q(\mathfrak{g}))$ is the full monoidal subcategory of $\Rep(U_q(\mathfrak{g}))$ generated by the irreducible representations with highest weight a fundamental weight. 
\end{defn}
 
We denote the set of fundamental weights of $\mathfrak{g}$ by $\{\varpi_i\}$. Let $\lambda$ be a dominant integral weight. Then we can write $\lambda = \sum n_i \varpi_i$, where $n_i\in \mathbb{Z}_{\ge 0}$. There is a partial order on all weights, where $\mu\le \lambda$ when $\lambda- \mu$ is a $\mathbb{Z}_{\ge 0}$-linear combination of positive roots. With respect to this partial order, the irreducible representation $V(\lambda)$ has highest weight $\lambda$. Also, $V(\lambda)$ is a direct summand of the tensor product $\bigotimes_i {V(\varpi_i)^{\otimes n_i} }$. Thus, there are projection and inclusion maps $\bigotimes_i V(\varpi_i)^{\otimes n_i}\rightarrow V(\lambda) \rightarrow \bigotimes_i V(\varpi_i)^{\otimes n_i}$ such that the composition is an idempotent $C_{\lambda}$ in $\Fund(U_q(\mathfrak{g}))$. We are interested in finding explicit descriptions of these idempotents, generalizing Formula \eqref{JonesWenzl}.  

Unless $\lambda$ is a fundamental weight or zero, $V(\lambda)$ will not be an object in $\Fund(U_q(\mathfrak{g}))$. However, $C_{\lambda}$ is a morphism in $\Fund(U_q(\mathfrak{g}))$ and we think of it as a replacement for $V(\lambda)$. Analogous to how $V(\lambda)$ is characterized as the irreducible representation with highest weight $\lambda$, the morphism $C_{\lambda}$ is characterized as the non-zero idempotent endomorphism of $\bigotimes_i V(\varpi_i)^{\otimes n_i}$ such that if $f: \bigotimes_i V(\varpi_i)^{\otimes n_i}\rightarrow V(\varpi_{i_1})\otimes \dots \otimes V(\varpi_{i_r})$ is a morphism in $\Fund(U_q(\mathfrak{g}))$, and $\sum_{k=1}^r\varpi_{i_k} < \lambda$, then $f\circ C_{\lambda} = 0$.

The Temperley-Lieb category gives a generators and relations description of $\Fund(U_q(\mathfrak{sl}_2))$. Kuperberg's paper \cite{Kupe} was the first attempt to generalize this to other $\mathfrak{g}$. In this work, Kuperberg gives generators and relations descriptions of the monoidal categories $\Fund(U_q(\mathfrak{g}))$ when $\mathfrak{g}$ is a rank two simple Lie algebra, i.e. $\mathfrak{g} = \mathfrak{sl}_3$ , $\mathfrak{sp}_4$, or $\mathfrak{g}_2$. This work was later extended to types $A$ \cite[Theorem 3.3.1]{CKM} and $C$ \cite[Theorem 1.4]{bodish2021type}. Following these authors we refer to a diagrammatic category which is monoidally equivalent to $\Fund(U_q(\mathfrak{g}))$ as $\textbf{Web}_q(\mathfrak{g})$. It remains an open problem to define web categories for all simple $\mathfrak{g}$. 

Kuperberg also introduced the terminology \emph{clasp} to refer to an idempotent projecting to the highest weight irreducible summand of a tensor product of fundamental representations, viewed as a morphism in $\textbf{Web}_q(\mathfrak{g})$. If the highest weight of this irreducible summand is $\lambda$, then we will call this idempotent a $\lambda$-clasp. To generalize the Jones-Wenzl recursion, a first step is to find recursive formulas of clasps in the rank two cases.

In the $\mathfrak{sl}_3$ case, a recursive formula was given by Ohtsuki and Yamada \cite[Definition 2.4]{OhtsukiYamada}, where they called a clasp a ``magic element". Later, Dongseok Kim found other recursive formulas for the $\mathfrak{sl}_3$ case as well \cite[Theorem 3.3]{Kim07}. 

In \cite[Conjecture 3.16]{elias2015light}, Elias made his type $A$ clasp conjecture, which implies a recursive description of each $\mathfrak{sl}_n$ clasp using the language of $\mathfrak{sl}_n$ webs. Also, \cite[Theorem 2.57]{elias2015light} provides a basis for all homomorphism spaces between fundamental representations for $\mathfrak{sl}_n$. These bases have a particularly nice form which reduces the validity of the type $A$ clasp conjecture to an explicit (but unwieldy) calculation. In \cite[Section 3.4 and Section 3.5]{elias2015light} Elias proved the $\mathfrak{sl}_n$ clasp conjecture by hand for $n\le 4$.

Since $\mathfrak{sp}_4$ is rank two, there are two simple roots: one short $\alpha_1$ and one long $\alpha_2$. We write $\varpi_1$ and $\varpi_2$ for the corresponding fundamental weights. In the $\mathfrak{sp}_4$ case, Kim gave recursive clasp formulas for the $a\varpi_1$-clasp \cite[Corollary 4.3]{Kim07} and the $b\varpi_2$-clasp \cite[Corollary 4.5]{Kim07}. However, an inductive formula for the $\mathfrak{sp}_4$ $a\varpi_1 + b\varpi_2$-clasp remained unknown until recently, when the first named author derived formulas generalizing Elias's type $A$ clasp conjecture to type $C_2$ \cite[Theorem 1.5]{Sp4Clasp}.

In the $\mathfrak{g}_2$ case, little was known before the present work. Attempts at getting the $\mathfrak{g}_2$-clasp formulas have been made, including a few base-case calculations by Sakamoto and Yonezawa \cite[Section 5]{G2basecasecalc}. In this paper, we give triple clasp expansions for the $\mathfrak{g}_2$ $\lambda$-clasps for all dominant integral weights $\lambda$. Our results are summarized in the following theorem.

\begin{thm}\label{MainResult}

\def\neckcollarl{
\begin{trics}
\draw [darkred,thick,scale=1.7] 
      (0,0.2)rectangle(1.8,-0.2) node[pos=0.5,black,scale=0.7] {$\lambda+\mu$}
      (0,0.7)--(-0.3,1.3)--(2.1,1.3)--(1.8,0.7)--cycle 
      (0.9,1)node[black,scale=0.7] {$\mathbb{D}(^iELL_{\lambda, \varpi}^{\lambda+ \mu})$}
       (0,-0.7)--(-0.3,-1.3)--(2.1,-1.3)--(1.8,-0.7)--cycle 
      (0.9,-1)node[black,scale=0.7] {$^jELL_{\lambda, \varpi}^{\lambda+ \mu}$}
      (-0.3,-1.8)rectangle(0.9,-2.2) node[pos=0.5,black,scale=0.7] {$\lambda$}
     (-0.3,1.8)rectangle(0.9,2.2) node[pos=0.5,black,scale=0.7] {$\lambda$}
      ;             
\draw [scale=1.7] (0.9,0.2)--(0.9,0.7) node[midway,left,scale=0.7] {$\lambda+\mu$}
       (0.9,-0.2)--(0.9,-0.7) node[midway,left,scale=0.7] {$\lambda+\mu$}
      (0.4,-1.3)--(0.4,-1.8) node[midway,left,scale=0.7] {$\lambda$}
       (0.4,1.3)--(0.4,1.8) node[midway,left,scale=0.7] {$\lambda$}
       (0.4,-2.2)--(0.4,-2.7) node[midway,left,scale=0.7] {$\lambda$}
       (0.4,2.2)--(0.4,2.7) node[midway,left,scale=0.7] {$\lambda$}
       (1.5,1.3)--(1.5,2.7) node[midway,left,scale=0.7] {$\varpi$}
       (1.5,-1.3)--(1.5,-2.7) node[midway,left,scale=0.7] {$\varpi$};
\end{trics}
}

\def\neckcollarlo{
\begin{trics}
\draw [darkred,thick,scale=1.7] 
      (0,0.2)rectangle(1.8,-0.2) node[pos=0.5,black,scale=0.7] {$\lambda+\varpi$}; 
      ;             
\draw [scale=1.7](0.9,0.2)--(0.9,2.7) node[midway,left,scale=0.7] {$\lambda+\varpi$}
       (0.9,-0.2)--(0.9,-2.7) node[midway,left,scale=0.7] {$\lambda+\varpi$}; 
\end{trics}
}

\def\neckcollarlol{
\begin{trics}
\draw [darkred,thick,scale=1.7] 
      (-0.3,-0.2)rectangle(0.9,0.2) node[pos=0.5,black,scale=0.7] {$\lambda$};
      ;             
\draw  [scale=1.7](0.4,0.2)--(0.4,2.7) node[midway,left,scale=0.7] {$\lambda$}
       (0.4,-0.2)--(0.4,-2.7) node[midway,left,scale=0.7] {$\lambda$}
       (1.5,2.7)--(1.5,-2.7) node[midway,left,scale=0.7] {$\varpi$};
\end{trics}
}

Let $\lambda$ be a dominant integral weight for $\mathfrak{g}_2$. Then the $\lambda+\varpi$ clasp is given by the following recursive formula.
 
 $$  \neckcollarlo =   \neckcollarlol \ \ \ - \sum_{\mu\in V(\varpi)\backslash\lbrace{\varpi\}}} \left( \left( {^{p\ell}t}_{\lambda, \varpi}^{\mu} \right)^{-1} \right) _{ij}\cdot  \neckcollarl$$

Here a red box labelled by $\chi$ denotes the $\chi$-clasp, i.e. the morphism in $\textbf{Web}_q(\mathfrak{g}_2)$ which corresponds to the idempotent with image $V(\chi)$. The diagrams $\mathbb{D}(^iELL_{\lambda, \varpi}^{\lambda+ \mu})$ and  $^jELL_{\lambda, \varpi}^{\lambda+ \mu}$ are given explicitly in Formula \eqref{FirstTripleClaspsExpansionGraph} and Formula \eqref{SecondTripleClaspsExpansionGraph} in Section \ref{results}. The coefficients  ${^{p\ell}t}_{\lambda, \varpi}^{\mu}$ are given explicitly by Equation \eqref{FirstExplicit} to Equation \eqref{DeterminantR00} in Section \ref{results} and Equation \eqref{explicitinappendix} to Equation \eqref{LastExplicit} in Appendix \ref{Relations}.
    
 \end{thm}

\subsection{Connection to the Clasp Conjecture}
\label{sec-claspconj}

Let $\mathbb{F}$ be a field. Consider objects $X$ and $S$ in an additive $\mathbb{F}$-linear Karoubian category with duality $\mathbb{D}$, i.e. a contravariant endofunctor with $\mathbb{D}^2\cong \id$, such that $\End(S)= \mathbb{F}\cdot \id_{S}$, $\mathbb{D}(X) = X$, and $\mathbb{D}(S) = S$. Given $\pi: X\rightarrow S$, we obtain a map $\iota = \mathbb{D}(\pi): S\rightarrow X$ and 
\[
\pi\circ \iota = \kappa\id_{S},
\]
for some $\kappa \in \mathbb{F}$. If $\kappa\ne 0$, then $S$ is isomorphic to the image of the idempotent 
\[
e = \dfrac{1}{\kappa} \iota \circ \pi.
\]
In \cite[Definition 3.8]{elias2015light}, the coefficient $\kappa$, computed in the $\mathfrak{sl}_n$ web category, is called a local intersection form. We carry out analogous calculations in the $\mathfrak{g}_2$ web category and find an analogue of Elias's clasp conjecture holds for $G_2$. 

In fact, we expect that something in general will hold. Let $\mathfrak{g}$ be a simple Lie algebra and let $U_q(\mathfrak{g})$ be the associated quantum group. Let $W$ denote the Weyl group associated to $\mathfrak{g}$. For $V\in \Rep(U_q(\mathfrak{g}))$ we will write $\wt V$ to denote the set of all weights $\mu$ such that the $\mu$ weight space of $V$ is non-zero. Fix a fundamental weight $\varpi$. For each $\mu \in \wt V(\varpi)$, such that $\mu$ is in the same $W$ orbit as $\varpi$, there should be a clasped \emph{elementary light ladder map}\footnote{In examples, this is most easily defined using web categories. But with some care should make sense in general, even without having a generators and relations presentation of $\Fund(\mathfrak{g})$. The main feature should be that the map is the composition of projectors and neutral maps with some fixed map from $V(\lambda_{min})\ot V(\varpi)\rightarrow V(\lambda_{min}+ \mu)$, where $\lambda_{min}$ is the smallest dominant weight so that $V(\lambda)\otimes V(\varpi)$ contains a copy of $V(\lambda+ \mu)$.}
\[
LL_{\lambda, \varpi}^{\lambda+ \mu}: V(\lambda)\otimes V(\varpi) \rightarrow V(\lambda+ \mu).
\]

\begin{defn}
For an extremal weight $\mu$ in a fundamental representation (i.e. a weight in the $W$ orbit of $\varpi_i$ for some $i$) we write $d_{\mu}$ to denote the minimal length element $w\in W$ so that $w(\mu)$ is dominant. We also define $\Phi_{\mu}$ to be the set of positive roots which are sent to negative roots by $d_{\mu}$. 
\end{defn}

\begin{conjecture}
If we denote by $\mathbb{D}$ the duality\footnote{Again, this is most easily defined in terms of webs, in which case it is just flipping the diagram upside down.} on $\Fund(\mathfrak{g})$, and write 
\[
LL_{\lambda, \varpi}^{\lambda+ \mu}\circ \mathbb{D}(LL_{\lambda, \varpi}^{\lambda+ \mu}) = \kappa_{\lambda, \varpi}^{\mu}\id_{V(\lambda+ \mu)},
\]
then
\begin{equation}
\kappa_{\lambda, \varpi}^{\mu} = \pm\prod_{\alpha\in {\Phi_{\mu}}}\frac{[(\alpha^{\vee}, \lambda + \rho)]_{q^{\ell(\alpha)}}}{[(\alpha^{\vee}, \lambda +\mu + \rho)]_{q^{\ell(\alpha)}}}.\label{claspconjround}
\end{equation}
Here $\rho$ is the sum of the fundamental weights and $l(\alpha) =(\alpha, \alpha)/2$.
\end{conjecture}

\begin{rmk}
The conjecture is known to be true in types $A_n$, for $n\le 4$ \cite[Section 3]{elias2015light} and type $C_2$ \cite[Corollary 3.7]{Sp4Clasp}. 
\end{rmk}

The following Proposition is an elementary consequence of our main theorem.

\begin{proposition}
The conjecture is true for $\mathfrak{g}_2$. 
\end{proposition}

\begin{proof}
See Corollary \ref{cor-claspconj}.
\end{proof}

\begin{rmk}
We also expect there to be a more general form of the conjecture which describes what happens for $\mu\in V(\varpi)$ which are not in the extremal Weyl orbit. The work in this paper and \cite{Sp4Clasp} could give enough data to guess the answer when $V(\varpi)_{\mu}$ is one dimensional, but we have not yet carried this out. We also hope the general form of the conjecture will give rise to a product formula which computes the elementary divisors of the matrix of local intersection forms when $\dim V(\varpi)_{\mu}> 1$. 

\end{rmk}


\subsection{Roundness conjecture for the determinant of trihedron coefficients}
Notice that in the clasp conjecture, the right hand side of Formula \eqref{claspconjround} is written as a quotient of products of quantum integers. An element of $\mathbb{C}(q)$ is called q-round when it can be written this way. One way to generalize the clasp conjecture, without giving explicit formulas, is to look for q-roundness in other structure constants in web categories. 

\begin{definition}
Suppose $V_{\lambda_1}, V_{\lambda_2}, V_{\lambda_3}  \in \Rep(U_q(\mathfrak{g}))$, we define the trihedron coefficients $M_\theta(V_{\lambda_1},V_{\lambda_2},V_{\lambda_3})$ as the following bilinear pairing:
\begin{align*}
    M_\theta(V_{\lambda_1},V_{\lambda_2},V_{\lambda_3}): \Hom(\mathbb{C}(q),V_{\lambda_1}\ot V_{\lambda_2}\ot V_{\lambda_3}) \ot \Hom(V_{\lambda_1}\ot V_{\lambda_2}\ot V_{\lambda_3},\mathbb{C}(q)) &\longrightarrow \mathbb{C}(q)\\
    f\ot g &\mapsto g\circ f(1)
\end{align*}
\end{definition}

\begin{remark}
Lusztig defined a canonical basis $\mathbb{B}$ in $U_q(\mathfrak{g})^+$ \cite{MR1035415}, which descends to a canonical basis in any irreducible representation. In \cite{MR1180036}, he further explains how the canonical basis in an irreducible representation can be used to obtain a canonical basis for a tensor product of irreducible representations. Moreover, in \cite{LuszQuantum}[27.2.5], it is explained how the canonical basis in a tensor product naturally gives rise to a basis for the space of coinvariants. In \cite{MR1180036}, a dual canonical basis is defined for each irreducible representation and their tensor products. The dual canonical bases naturally give rise to a basis for the space of invariant vectors in a tensor product of representations. 

Write $\lambda^*:= -w_0(\lambda)$, where $w_0$ is the longest element in the Weyl group. Since $V(\lambda)^*\cong V(\lambda^*)$, we can use the isomorphism
\[
\Hom(V_{\lambda_1}\otimes V_{\lambda_2}\otimes V_{\lambda_3}, \mathbb{C}(q))\cong \Hom(\mathbb{C}(q), V_{\lambda_3^*}\otimes V_{\lambda_2^*}\otimes V_{\lambda_1^*})
\]
to transport the dual canonical basis for the invariant space to a ``dual canonical basis" for the space of coinvariants.
We can then write the  matrix of the pairing $M_{\theta}$ with respect to these dual canonical bases. Let us write $M_{\theta}^{\mathbb{B}^*}$ to denote this matrix.

See \cite{MR2695927} and \cite{MR1680395} for more discussion of how the web basis and dual canonical basis are related.
\end{remark}

We have the following conjecture given by Greg Kuperberg \cite{KupConjecture}: 
\begin{conj}\label{Kuperconj}
The determinant of the matrix $M_{\theta}^{\mathbb{B}^*}$ is $q$-round.
\end{conj}

\begin{remark}
Since the web basis and the dual canonical basis are the same for $U_q(\mathfrak{sl}_2)$ \cite{MR2695927}, it follows from \cite{MR1272656} that Conjecture \ref{Kuperconj} is true in this case.
\end{remark}

\begin{remark}
It is not the case that the web basis and the dual canonical basis are the same in general. One might ask if the change of basis matrix between the web basis and the dual canonical basis is unitriangular. This is believed to be true in the rank two cases, especially for $\mathfrak{sl}_3$ \cite{MR1680395}.
\end{remark}
\begin{remark}
One might expect that the determinant of the pairing matrix $M_{\theta}$ will not change if it is computed using bases which are a unitriangular change of basis from the dual canonical bases. If both the above desired claims hold, then one could compute the determinant of trihedron coefficients by using the graphical calculus, and verify the conjecture. 
\end{remark}

\begin{remark}
When $\lambda_3=\varpi$, the clasped light ladder diagrams from Definition \ref{claspedbasis} give rise to a basis for  $\Hom(\mathbb{C}(q),V_{\lambda_1}\ot V_{\lambda_2}\ot V_{\lambda_3})$. In this case, the trihedron coefficients $M_\theta(V_{\lambda_1},V_{\lambda_2},V_{\lambda_3})$
is equal to the local intersection form matrix $\Bigg( {^{pq}\kappa_{\lambda_2, \varpi}^{\lambda_1}} \Bigg)$ multiplied by the quantum dimension $\dim(V_{\lambda_1})$. The quantum dimension $\dim(V_{\lambda_1})$ must be q-round, and in the $\mathfrak{g}_2$ case, we know that the determinant of  $\Bigg( {^{pq}\kappa_{\lambda_2, \varpi}^{\lambda_1}} \Bigg)$ is q-round from Section \ref{results}. So the determinant of $M_\theta(V_{\lambda_1},V_{\lambda_2},V_{\lambda_3})$ is also q-round. 

If we could show that the clasped web basis was a unitriangular change of basis away from the dual canonical basis, then using the present work it should be possible to argue that Conjecture \ref{Kuperconj} holds in the $\mathfrak{g}_2$ case when $\lambda_3= \varpi$. Similarly, the $\mathfrak{sl}_3$ case and the $\mathfrak{sp}_4$ case would follow from results in \cite{elias2015light,Sp4Clasp} when $\lambda_3=\varpi$.

In the case when $\lambda_3$ is a fundamental weight and $\lambda_2 - \lambda_1$ is in the $W$ orbit of $\lambda_3$, the invariant space in which we view the clasped web basis element is one dimensional. Thus, if the clasped web basis is a unitriangular change of basis matrix away from the dual canonical basis element, then the clasped web basis is in fact equal to the dual canonical basis element. It would be interesting to determine whether the clasped web basis always agreed with the dual canonical basis in this case.
\end{remark}

\begin{remark}
When $\lambda_i \neq \varpi$, calculations for some of the trihedron coefficients $M_\theta(V_{\lambda_1},V_{\lambda_2},V_{\lambda_3})$ are also done in the $\mathfrak{sl}_3$ case, which provides further evidence to the conjecture \cite{MR2221529}.
\end{remark}

\subsection{Reshetikhin-Turaev invariants via skein-theory}

Let $\mathfrak{g}$ be a simple Lie algebra over $\mathbb{C}$. In order to define an analogue of the Jones polynomial for $\mathfrak{g}$, Reshetikhin and Turaev defined a link invariant using the category $\textbf{Rep}(U_q(\mathfrak{g}))$ \cite{RTlink}. Their construction gives a knot invariant for every type-$\textbf{1}$ representation of $U_q(\mathfrak{g})$. More generally, one can label each component of a link with an object in $\textbf{Rep}(U_q(\mathfrak{g}))$ and their construction gives the colored Reshetikhin-Turaev link invariant.

Kuperberg's original motivation for studying $\DD$ was to compute the Reshitikhin-Turaev link invariant associated to $\mathfrak{g}_2$. Originally, he gave a diagrammatic method to compute the $\mathfrak{g}_2$ link invariant when each component is colored by the first fundamental representation \cite{Kupe-first-G2}. Soon after, using $\DD$, he gave diagrammatic tools for computing the $\mathfrak{g}_2$ link invariant colored by both fundamental representations \cite[Section 4]{Kupe}. In this paper, we give explicit formulas for idempotents projecting to each irreducible $U_q(\mathfrak{g}_2)$ module. Combined with Kuperberg's earlier work this gives a diagrammatic approach to computing the Reshitikhin-Turaev invariant of a link with components colored by any irreducible. 

Reshitkhin-Turaev's paper about their link invariant was intended as a prequel to their work which gave an associated $3$-manifold invariant \cite{RTinv}. The first step one takes to make sense of their $3$-manifold invariant is to leave behind representation theory of $U_q(\mathfrak{g})$ for generic $q$ and work instead with $q$ specialized to a root of unity. 

Let $U_{\mathbb{Z}[q,q^{-1}]}(\mathfrak{g}_2)$ be the $\mathbb{Z}[q,q^{-1}]$-subalgebra of $\Uq$ generated by $\frac{E_{\alpha}^{(a)}}{[a]_{q^{\ell(\alpha)}}!}$, $\frac{F_{\alpha}^{(a)}}{[a]_{q^{\ell(\alpha)}}!}$, and $K_{\alpha}^{\pm 1}$, for all simple roots $\alpha$, and all $a\in \mathbb{Z}_{\ge 0}$. When $\xi$ is a root of unity in $\mathbb{C}$, we can study the relation between $\mathbb{C}\ot_{q= \xi}\DD$, and the category of tilting modules of $\mathbb{C}\ot_{q= \xi}U_{\mathbb{Z}[q,q^{-1}]}(\mathfrak{g}_2)$. It is possible to adapt the approach from \cite{Sp4tilt}, which itself is based on \cite{elias2015light}, to prove that the Karoubi envelope of $\mathbb{C}\ot_{q= \xi}\DD$ is equivalent to the category of tilting modules as long as $[2]_{\xi}, [3]_{\xi}\ne 0$. The same result is work in progress of Victor Ostrik and Noah Snyder, but they propose a slightly different approach.


When $\xi$ is a root of unity of order greater than $5$, the generators of the negligible ideal in the category of $U_{\xi}(\mathfrak{g}_2)$ tilting modules are (identity morphisms of) certain irreducible tilting modules with quantum dimension zero. Irreducible tilting modules are also Weyl modules, so these generating objects correspond to clasps in $\DD$. Moreover, the objects which survive in the negligible quotient are the irreducible Weyl modules with non-zero quantum dimension. Once the equivalence between the web category and the category of tilting modules is established, one can give a generators and relations presentation of the associated modular tensor category\footnote{In the case that $\xi$ is a root of unity that actually gives rise to a modular category as the negligible quotient \cite{rowellqmod}.} using the $\mathfrak{g}_2$ triple clasp formulas for the negligible clasps. Combined with our description of the clasps corresponding to the irreducible Weyl modules with non-zero dimension, this gives an explicit way to compute the quantum $\mathfrak{g}_2$ $3$-manifold invariant \cite{RTinv}.

On the other hand, topologists have tried to understand the quantum $3$-manifold invariants with graphical categories. A construction of the quantum $\mathfrak{sl}_2$ $3$-manifold invariant using the Temperley-Lieb category was given by Lickorish \cite{Lickorish}. This work was generalized to the $\mathfrak{sl}_3$ case by Ohtsuki and Yamada \cite{OhtsukiYamada} with $\mathfrak{sl}_3$ webs. A self-contained proof of invariance under Kirby moves \cite{KirbyMove} using the graphical category was given in both cases. One can now give similar constructions and proofs in the $\mathfrak{g}_2$ case by using our clasp formulas.

\subsection{Structure of the Paper}
\label{structure}

In Section \ref{g2spider}, we recall the definition of $\mathfrak{g}_2$ webs, and also the relation between $\mathfrak{g}_2$ webs and representations of the quantum group associated to $\mathfrak{g}_2$ \cite{Kupe}. Then we give the definition of $\mathfrak{g}_2$ clasps, and introduce the elementary light ladders as building blocks of the $\mathfrak{g}_2$ triple clasp expansions.

In Section \ref{FormulaProof}, we give the triple clasp formulas for $\mathfrak{g}_2$ explicitly, and then prove that a linear combination of webs given inductively by the $\mathfrak{g}_2$ triple clasp expansions satisfies the definition of a $\mathfrak{g}_2$ clasp. 

In Appendix \ref{Recursions}, we give a complete list of the recursive formulas which determine the coefficients in the $\mathfrak{g}_2$ triple clasp formulas. All the graphical calculations that lead to the recursive formulas can be found in Appendix \ref{GraphicalCalc}. The SAGE code verifying that the triple clasp coefficients satisfy the recursive formulas is included with the source file of the arXiv submission of this paper.


\subsection{Acknowledgments}
\label{ack}

Thank you to our advisors Greg Kuperberg and Ben Elias. Elijah Bodish was supported by the University of Oregon's Lokey Graduate Science Fellowship. Haihan Wu was partially supported by NSF grant CCF-1716990.

\section{Webs and Clasps for $\mathfrak{g}_2$}
\label{g2spider}

\subsection{Definition of $\mathfrak{g}_2$ Webs}
\label{spiderdefn}

\def\Skeina
{\begin{tric}
\draw (0.75,0) circle (0.75);
\end{tric}
}

\def\Skeinb
{\begin{tric}
\draw[double,thin] (0.75,0) circle (0.75);
\end{tric}
}

\def\Skeinc
{\begin{tric}
\draw (0,0.7)--(0,1.5) (0,-0.7)--(0,-1.5)
      (0,0.7)..controls(-0.5,0.7)and(-0.5,-0.7)..(0,-0.7)  (0,0.7)..controls(0.5,0.7)and(0.5,-0.7)..(0,-0.7);
\end{tric}
}

\def\Skeinca
{\begin{tric}
\draw (0,1.5)--(0,-1.5);
\end{tric}
}

\def\Skeind
{\begin{tric}
\draw [scale=0.65] (-1,0)--(1,0)--(0,1.732)--cycle
                  (-1,0)--(-2,-0.577) (1,0)--(2,-0.577) (0,1.732)--(0,2.887);
\end{tric}
}

\def\Skeinda
{\begin{tric}
\draw[scale=0.65]
    (0,0.577)--(-2,-0.577) (0,0.577)--(2,-0.577) (0,0.577)--(0,2.887);
\end{tric}
}

\def\Skeindb
{\begin{tric}
\draw [scale=0.65] (0,0.7)..controls(0.5,0.7)and(0.5,1.732)..(0,1.732)  (0,0.7)..controls(-0.5,0.7)and(-0.5,1.732)..(0,1.732)  
(-1,-0.5)--(0,0)--(1,-0.5) (0,1.732)--(0,2.887);
\draw[double,thin,scale=0.65]  (0,0)--(0,0.7) ;
\end{tric}
}

\def\Skeindc
{\begin{tric}
\draw [scale=0.65] (0,0.7)..controls(0.5,0.7)and(0.5,1.732)..(0,1.732)  (0,0.7)..controls(-0.5,0.7)and(-0.5,1.732)..(0,1.732)  
(-1,-0.5)--(0,0)--(1,-0.5) (0,0)--(0,0.7) (0,1.732)--(0,2.887)  ;
\end{tric}
}

\def\Skeindd
{\begin{tric}
\draw [scale=0.65] (0,0.7)..controls(0.5,0.7)and(0.5,1.732)..(0,1.732)  (0,0.7)..controls(-0.5,0.7)and(-0.5,1.732)..(0,1.732)  
(-1,-0.5)..controls(0,0)..(1,-0.5)  (0,1.732)--(0,2.887)  ;
\end{tric}
}

\def\Skeine
{\begin{tric}
\draw [scale=0.5](-1,-1)--(-1,1)--(1,1)--(1,-1)--cycle;
\draw [scale=0.5](-1,-1)--(-2,-2) (-1,1)--(-2,2) (1,1)--(2,2) (1,-1)--(2,-2);
\end{tric}
}

\def\Skeinea
{\begin{tric}
\draw[scale=0.5] (-2,-1.732)--(-1,0)--(-2,1.732) (2,1.732)--(1,0)--(2,-1.732) (-1,0)--(1,0);
\end{tric}
}

\def\Skeineb
{\begin{tric}
\draw [scale=0.5]  (-1.732,2)--(0,1)--(1.732,2) (1.732,-2)--(0,-1)--(-1.732,-2) (0,-1)--(0,1);
\end{tric}
}

\def\Skeinec
{\begin{tric}
\draw [scale=0.45] (-2,-2)..controls(-1,-1)and(-1,1)..(-2,2) (2,2)..controls(1,1)and(1,-1)..(2,-2);
\end{tric}
}

\def\Skeined
{\begin{tric}
\draw  [scale=0.45] (-2,2)..controls(-1,1)and(1,1)..(2,2) (2,-2)..controls(1,-1)and(-1,-1)..(-2,-2) ;
\end{tric}
}

\def\Skeinee
{\begin{tric}
\draw [scale=0.5](0,-0.2)--(-0.7,1)--(0.7,1)--cycle;
\draw [scale=0.5] (-0.7,1)--(-1.5,2) (0.7,1)--(1.5,2) 
                   (-1.5,-2)--(0,-1)--(1.5,-2);
\draw [double,thin,scale=0.5] (0,-1)--(0,-0.2);
\end{tric}
}

\def\Skeinef
{\begin{tric}
\draw [scale=0.5](0,-0.2)--(-0.7,1)--(0.7,1)--cycle;
\draw [scale=0.5] (-0.7,1)--(-1.5,2) (0.7,1)--(1.5,2) 
                   (-1.5,-2)--(0,-1)--(1.5,-2)  (0,-1)--(0,-0.2);
\end{tric}
}

\def\Skeineg
{\begin{tric}
\draw [scale=0.5](-0.7,1)--(0.7,1) (-0.7,1)..controls(-0.2,-0.2)and(0.2,-0.2) ..(0.7,1);
\draw [scale=0.5] (-0.7,1)--(-1.5,2) (0.7,1)--(1.5,2) 
                   (-1.5,-2)..controls(0,-1)..(1.5,-2)  ;
\end{tric}
}

\def\Skeinf
{\begin{tric}
\draw [scale=0.6] (18:1)--(90:1)--(162:1)--(234:1)--(306:1)--cycle
      (18:1)--(18:2) (90:1)--(90:2) (162:1)--(162:2) (234:1)--(234:2) (306:1)--(306:2);
\end{tric}
}

\def\Skeinfa
{\begin{tric}
\draw [scale=0.6] (18:1)--(90:1)--(162:1)
      (18:1)--(18:2) (90:1)--(90:2) (162:1)--(162:2) (162:1)..controls(234:1)..(234:2) (18:1)..controls(306:1)..(306:2);
\end{tric}
}

\def\Skeinfb
{\begin{tric}
\draw [scale=0.6] (90:1)--(162:1)--(234:1)
      (90:1)..controls(18:1)..(18:2) (90:1)--(90:2) (162:1)--(162:2) (234:1)--(234:2) (234:1)..controls(306:1)..(306:2);
\end{tric}
}

\def\Skeinfc
{\begin{tric}
\draw [scale=0.6] (162:1)--(234:1)--(306:1)
      (306:1)..controls(18:1)..(18:2) (162:1)..controls(90:1)..(90:2) (162:1)--(162:2) (234:1)--(234:2) (306:1)--(306:2);
\end{tric}
}

\def\Skeinfd
{\begin{tric}
\draw [scale=0.6] (234:1)--(306:1)--(18:1)
      (18:1)--(18:2) (18:1)..controls(90:1)..(90:2) (234:1)..controls(162:1)..(162:2) (234:1)--(234:2) (306:1)--(306:2);
\end{tric}
}

\def\Skeinfe
{\begin{tric}
\draw [scale=0.6] (306:1)--(18:1)--(90:1)
      (18:1)--(18:2) (90:1)--(90:2) (90:1)..controls(162:1)..(162:2) (306:1)..controls(234:1)..(234:2) (306:1)--(306:2);
\end{tric}
}

\def\Skeinff
{\begin{tric}
\draw  [scale=0.6] (90:1)..controls(18:1)..(18:2) (90:1)--(90:2)   
      (90:1)..controls(162:1)..(162:2) (234:2)..controls(234:1)and(306:1)..(306:2);
\end{tric}
}

\def\Skeinfg
{\begin{tric}
\draw    [scale=0.6] (162:1)..controls(90:1)..(90:2)
        (162:1)--(162:2) (162:1)..controls(234:1)..(234:2) (18:2)..controls(18:1)and(306:1)..(306:2);
\end{tric}
}

\def\Skeinfh
{\begin{tric}
\draw    [scale=0.6]  (18:2)..controls(18:1)and(90:1)..(90:2) 
         (234:1)..controls(162:1)..(162:2) (234:1)--(234:2) (234:1)..controls(306:1)..(306:2);
\end{tric}
}

\def\Skeinfi
{\begin{tric}
\draw  [scale=0.6] (306:1)..controls(18:1)..(18:2)  
      (90:2)..controls(90:1)and(162:1)..(162:2)
      (306:1)..controls(234:1)..(234:2) (306:1)--(306:2);
\end{tric}
}

\def\Skeinfj
{\begin{tric}
\draw      [scale=0.6]  (18:1)--(18:2) (18:1)..controls(90:1)..(90:2)
           (234:2)..controls(234:1)and(162:1)..(162:2)  (18:1)..controls(306:1)..(306:2);
\end{tric}
}

\def\Skeinfk
{\begin{tric}
\draw [scale=0.6] (18:1)--(90:1)--(162:1)
      (18:1)--(18:2) (90:1)--(90:2) (162:1)--(162:2) (306:2)--(270:1)--(234:2) 
      (18:1)--(0,-0.3)--(162:1);
\draw [scale=0.6,double,thin] (0,-0.3)--(0,-1);
\end{tric}
}

\def\Skeinfl
{\begin{tric}
\draw [scale=0.6] (18:1)--(90:1)--(162:1)
      (18:1)--(18:2) (90:1)--(90:2) (162:1)--(162:2) (306:2)--(270:1)--(234:2) 
      (18:1)--(0,-0.3)--(162:1)  (0,-0.3)--(0,-1);
\end{tric}
}

\def\Skeinfm
{\begin{tric}
\draw [scale=0.6] (18:1)--(90:1)--(162:1)
      (18:1)--(18:2) (90:1)--(90:2) (162:1)--(162:2) (306:2)..controls(270:1)..(234:2) 
      (18:1)..controls(0,-0.3)..(162:1) ;
\end{tric}
}

\def\Skeinfn
{\begin{tric}
\draw [scale=0.6] (18:1)..controls(90:1)..(90:2)
      (0,-0.3)--(18:1)--(18:2)  (306:2)--(270:1)--(234:2) 
      (0,-0.3)..controls(162:1)..(162:2);
\draw [scale=0.6,double,thin] (0,-0.3)--(0,-1);
\end{tric}
}

\def\Skeinfo
{\begin{tric}
\draw [scale=0.6] (90:2)..controls(90:1)..(162:1)
      (0,-0.3)..controls(18:1)..(18:2)  (162:1)--(162:2) (306:2)--(270:1)--(234:2) 
      (0,-0.3)--(162:1);
\draw [scale=0.6,double,thin] (0,-0.3)--(0,-1);
\end{tric}
}

\def\Skeing
{\begin{tric}
\draw    (0,-1.5)--(0,0) (0,0)..controls(0.7,0.5)and(0.7,1.5)..(0,1.5)..controls(-0.7,1.5)and(-0.7,0.5)..(0,0);
\end{tric}
}

\def\Skeinh
{\begin{tric}
\draw[scale=0.7,double,thin] (0,0.75)--(0,-0.75);
\draw [scale=0.7] (-0.5,1.616)--(0,0.75)--(0.5,1.616) (-0.5,-1.616)--(0,-0.75)--(0.5,-1.616);
\end{tric}
}

\def\Verta
{\begin{tric}
\draw [scale=0.8] (0,0)--(90:1) (0,0)--(210:1) (0,0)--(330:1);
\draw (90:1)node[black,anchor=south]{$\varpi_1$}
      (210:1)node[black,anchor=north]{$\varpi_1$}
      (330:1)node[black,anchor=north]{$\varpi_1$}; 
\end{tric}
}

\def\Vertb
{\begin{tric}
\draw[scale=0.8,double,thin] (0,0)--(90:1) ;
\draw[scale=0.8] (0,0)--(240:1) (0,0)--(300:1);
\draw (90:1)node[black,anchor=south]{$\varpi_2$}
      (240:1)node[black,anchor=north]{$\varpi_1$}
      (300:1)node[black,anchor=north]{$\varpi_1$};
\end{tric}
}

\def\Skeini{
\begin{tric}
\draw [double,thin](0,0.7)--(0,1.5) (0,-0.7)--(0,-1.5);
\draw (0,0.7)..controls(-0.5,0.7)and(-0.5,-0.7)..(0,-0.7)  
      (0,0.7)..controls(0.5,0.7)and(0.5,-0.7)..(0,-0.7);
\end{tric}
}

\def\Skeinia{
\begin{tric}
\draw[double,thin] (0,1.5)--(0,-1.5);
\end{tric}
}

\def\Skeinj{
\begin{tric}
\draw[double,thin](0,-1.5)--(0,0);
\draw (0,0)..controls(0.7,0.5)and(0.7,1.5)..(0,1.5)..controls(-0.7,1.5)and(-0.7,0.5)..(0,0);
\end{tric}
}

\def\Skeinja{
\begin{tric}
\draw[double,thin](0,-0.5)--(0,0) (0,-1.4)--(0,-2);
\draw (0,0)..controls(0.4,0.5)and(0.4,1)..(0,1)..controls(-0.4,1)and(-0.4,0.5)..(0,0);
\draw(0,-1.4)..controls(0.4,-1.4)and(0.4,-0.5)..(0,-0.5)
     (0,-1.4)..controls(-0.4,-1.4)and(-0.4,-0.5)..(0,-0.5);
\end{tric}
}

\def\Skeinjb{
\begin{tric}
\draw[double,thin] (0,-1.4)--(0,-2);
\draw (0,0)..controls(0.4,0.5)and(0.4,1)..(0,1)..controls(-0.4,1)and(-0.4,0.5)..(0,0)  (0,-0.5)--(0,0);
\draw(0,-1.4)..controls(0.4,-1.4)and(0.4,-0.5)..(0,-0.5)
     (0,-1.4)..controls(-0.4,-1.4)and(-0.4,-0.5)..(0,-0.5);
\end{tric}
}

\def\Skeinjc{
\begin{tric}
\draw[double,thin] (0,-1.4)--(0,-2);
\draw (0,0.5)ellipse(0.3 and 0.5);
\draw(0,-1.4)..controls(0.4,-1.4)and(0.4,-0.5)..(0,-0.5)
     (0,-1.4)..controls(-0.4,-1.4)and(-0.4,-0.5)..(0,-0.5);
\end{tric}
}

\def\Skeinjd{
\begin{tric}
\draw[double,thin](0,-1.5)--(0,-0.5);
\draw (0,-0.5)..controls(0.4,0.5)and(0.4,1.5)..(0,1.5)..controls(-0.4,1.5)and(-0.4,0.5)..(0,-0.5);
\end{tric}
}

\def\Skeinje{
\begin{tric}
\draw[double,thin] (0,-1.5)--(0,-0.5);
\draw (0.4,0.5)..controls(0.4,1.5)..(0,1.5)..controls(-0.4,1.5)..(-0.4,0.5)
      (0.4,0.5)--(0,-0.5)--(-0.4,0.5) (0.4,0.5)--(-0.4,0.5);
\end{tric}
}

\def\Skeink{
\begin{tric}
\draw [double,thin](0,-0.7)--(0,-1.5);
\draw (0,0.7)--(0,1.5) 
      (0,0.7)..controls(-0.5,0.7)and(-0.5,-0.7)..(0,-0.7)  
      (0,0.7)..controls(0.5,0.7)and(0.5,-0.7)..(0,-0.7);
\end{tric}
}

\def\Skeinl
{\begin{tric}
\draw [scale=0.65] (-1,0)--(1,0)--(0,1.732)--cycle
                   (0,1.732)--(0,2.887);
\draw[double,thin,scale=0.65] (-1,0)--(-2,-0.577) (1,0)--(2,-0.577);
\end{tric}
}

\def\Skeinm
{\begin{tric}
\draw [scale=0.65] (-1,0)--(1,0)--(0,1.732)--cycle
                  (-1,0)--(-2,-0.577) (1,0)--(2,-0.577) ;
\draw[double,thin,scale=0.65] (0,1.732)--(0,2.887);
\end{tric}
}

\def\Skeinmaa
{\begin{tric}
\draw[scale=0.65]
    (0,0.577)--(-2,-0.577) (0,0.577)--(2,-0.577) ;
\draw[double,thin,scale=0.65] (0,0.577)--(0,2.887);
\end{tric}
}

\def\Skeinma
{\begin{tric}
\draw[scale=0.65]
    (0,0.577)--(-1,-0.577) (0,0.577)--(1,-0.577) ;
\draw[double,thin,scale=0.65] (0,0.577)--(0,2.887);
\end{tric}
}

\def\Skeinmb
{\begin{tric}
\draw [scale=0.65] (0,0.7)..controls(0.5,0.7)and(0.5,1.732)..(0,1.732)  (0,0.7)..controls(-0.5,0.7)and(-0.5,1.732)..(0,1.732)  
(-1,-0.5)--(0,0)--(1,-0.5);
\draw[double,thin,scale=0.65] (0,1.732)--(0,2.887) (0,0)--(0,0.7) ;
\end{tric}
}

\def\Skeinmc
{\begin{tric}
\draw [scale=0.65] (0,0.7)..controls(0.5,0.7)and(0.5,1.732)..(0,1.732)  (0,0.7)..controls(-0.5,0.7)and(-0.5,1.732)..(0,1.732)  
(-1,-0.5)--(0,0)--(1,-0.5) (0,0)--(0,0.7);
\draw[double,thin,scale=0.65] (0,1.732)--(0,2.887)  ;
\end{tric}
}

\def\Skeinmd
{\begin{tric}
\draw [scale=0.65] (0,0.7)..controls(0.5,0.7)and(0.5,1.732)..(0,1.732)  (0,0.7)..controls(-0.5,0.7)and(-0.5,1.732)..(0,1.732)  
(-1,-0.5)..controls(0,0)..(1,-0.5) ;
\draw[double,thin,scale=0.65] (0,1.732)--(0,2.887)  ;
\end{tric}
}

\def\Skeinn
{\begin{tric}
\draw[scale=0.65] (1,0)--(0,1)--(-1,0)--(0,-1)--cycle 
      (1,0)--(2,0) (0,1)--(0,2) (-1,0)--(-2,0) ;
\draw[double,thin,scale=0.65] (0,-1)--(0,-2);
\end{tric}
}

\def\Skeinna
{\begin{tric}
\draw [scale=0.65](1,0)..controls(0,1)..(0,2)    (0,-1)..controls(-1,0)..(-2,0) 
      (0,-1)--(1,0)--(2,0)   ;
\draw[double,thin,scale=0.65] (0,-1)--(0,-2);
\end{tric}
}

\def\Skeinnb
{\begin{tric}
\draw [scale=0.65](-1,0)..controls(0,1)..(0,2)    (0,-1)..controls(1,0)..(2,0) 
      (0,-1)--(-1,0)--(-2,0) ;
\draw[double,thin,scale=0.65] (0,-1)--(0,-2);
\end{tric}
}

\def\Skeinnc
{\begin{tric}
\draw[scale=0.63] (1,0)--(0,-1)--(-0.2,0.2)--cycle  
      (1,0)--(2,0)  (0,2)--(-0.7,0.7)--(-2,0) ;
\draw[double,thin,scale=0.63] (0,-1)--(0,-2)  (-0.7,0.7)--(-0.2,0.2);
\end{tric}
}

\def\Skeinnd
{\begin{tric}
\draw[scale=0.63] (1,0)--(0,-1)  (1,0)..controls(-0.2,0.2)..(0,-1)  
      (1,0)--(2,0)  (0,2)..controls(-0.7,0.7)..(-2,0) ;
\draw[double,thin,scale=0.63] (0,-1)--(0,-2) ;
\end{tric}
}

\def\Skeinne
{\begin{tric}
\draw[scale=0.63] (1,0)--(0,-1)--(-0.2,0.2)--cycle  
      (1,0)--(2,0)  (0,2)--(-0.7,0.7)--(-2,0)  (-0.7,0.7)--(-0.2,0.2) ;
\draw[double,thin,scale=0.63] (0,-1)--(0,-2)  ;
\end{tric}
}

\begin{defn}\label{GraphicalSpider} 
The category $\DD$ is the strict pivotal $\mathbb{C}(q)$-linear category, whose objects are generated by self-dual objects $ \varpi_1$ and $\varpi_2$, and whose morphisms are generated by the following two trivalent vertices:
 \begin{align*}
     \Verta \in {\Hom}_{\DD}(\varpi_1\otimes\varpi_1,\varpi_1 )\ \ \ \ \ \ \text{and} \ \ \ \ \ \ \ \ \ \ \ \  \Vertb \in {\Hom}_{\DD}(\varpi_1\otimes\varpi_1,\varpi_2),
 \end{align*}
modulo the tensor-ideal generated by the following relations.

\begin{align*}
&(S1)\ \ \ \Skeina=\frac{[2][7][12]}{[4][6]} \ \ \ \ \ \ \ \ \ \ \ \ \ \ \ \ \ \ \ \ \ \ \ \ \ \ \ \ \ \ \ \ \ \ \ \ \ \ \ \
(S2)\ \ \ \Skeinb=\frac{[7][8][15]}{[3][4][5]}\\
&(S3)\Skeing=0 \ \ \ \ \ \ \ \ \ \ \ \ \ \ \ \ 
(S4)\ \ \ \Skeink=0\ \ \ \ \ \ \ \ \ \ \ \ \ \ \ \ (S5)\Skeinl=0
\\
&(S6)\ \ \ \Skeinc= -\frac{[3][8]}{[2][4]}\ \ \ \Skeinca\ \ \ \ \ \ \ \ \ \ \ \ \ \ \ \ \ \ \ \ \ \ \ \ \ \ \ \ \ \ \ \ \ \ \ \ \ \ \ \ \ \
(S7)\ \ \ \Skeini= \ -[2] \ \ \ \Skeinia\\
&(S8)\ \ \Skeinea\ \   = \  \frac{1}{[2]} \Skeinec \ \ + \ \ \ \Skeinh 
   \ \ \  - \  \frac{1}{[3]}\Skeineb \ - \  \frac{[4][6]}{[2]^2[12]}\Skeined
\end{align*}
\end{defn}

The tensor product of objects in $\DD$ is concatenation of words. Tensor product of morphisms is horizontal concatenation. Composition of morphisms is vertical stacking.

\def\compoexamplea
{\begin{tric}
\draw [scale=0.7] (0,0)--(0,2);
\draw[double,thin, scale=0.7] (1.5,1)--(1.5,2);
\draw [scale=0.7](1,0)..controls(1,0.5)..(1.5,1) (2,0)..controls(2,0.5)..(1.5,1);
\end{tric}
}

\def\compoexampleb
{\begin{tric}
\draw[scale=0.7](2,0)--(2,-2);
\draw[scale=0.7](0.5,-1)--(0.5,-2);
\draw[scale=0.7](1,0)..controls(1,-0.5)..(0.5,-1) (0,0)..controls(0,-0.5)..(0.5,-1);
\end{tric}
}

\def\compoexamplec
{\begin{tric}
\draw [scale=0.7] (-2.5,0)--(-2.5,2);
\draw[double,thin, scale=0.7] (-1,1)--(-1,2);
\draw [scale=0.7](-1.5,0)..controls(-1.5,0.5)..(-1,1) (-0.5,0)..controls(-0.5,0.5)..(-1,1);
\draw[scale=0.7](2,0)--(2,2);
\draw[scale=0.7](0.5,1)--(0.5,0);
\draw[scale=0.7](1,2)..controls(1,1.5)..(0.5,1) (0,2)..controls(0,1.5)..(0.5,1);
\end{tric}
}

\def\compoexampled
{\begin{tric}
\draw [scale=0.7] (0,2)--(0,4);
\draw[double,thin, scale=0.7] (1.5,3)--(1.5,4);
\draw [scale=0.7](1,2)..controls(1,2.5)..(1.5,3) (2,2)..controls(2,2.5)..(1.5,3);
\draw[scale=0.7](2,0)--(2,2);
\draw[scale=0.7](0.5,1)--(0.5,0);
\draw[scale=0.7](1,2)..controls(1,1.5)..(0.5,1) (0,2)..controls(0,1.5)..(0.5,1);
\end{tric}
}

\begin{example}
Let \ \ 
$f=\ \compoexamplea \in {\Hom}_{\DD}(\varpi_1 ^ {\otimes 3},\varpi_1 \otimes \varpi_2)
\ ,  \ \  g= \ \compoexampleb \in {\Hom}_{\DD}(\varpi_1 ^{\otimes 2}, \varpi_1 ^ {\otimes 3} ).
$
Then  
\[
f\otimes g =\ \compoexamplec \in {\Hom}_{\DD}(\varpi_1 ^ {\otimes 5},\varpi_1 \otimes \varpi_2 \otimes \varpi_1 ^{\otimes 3})
\]
and
\[
f\circ g =\ \compoexampled \in {\Hom}_{\DD}(\varpi_1 ^ {\otimes 2},\varpi_1 \otimes \varpi_2).
\]

\end{example}

\begin{rmk}
Definition \ref{GraphicalSpider} is equivalent to the definition of $\mathfrak{g}_2$ Spider given in \cite[Section 4]{Kupe}, with some rescaling of both types of trivalent vertices. See Lemma \ref{SkeinRelations}.
\end{rmk}

\begin{notation}
When $\unw = \varpi_1^{\otimes a}\otimes \varpi_2 ^{\otimes b}$ we call $\unw$ segregated. For example, $\varpi_1\varpi_1\varpi_2$ is segregated, while $\varpi_1\varpi_2\varpi_1$ is not. 

\def\labedgea{
\begin{tric}
\draw (0,-0.5)--(0,0.5) node[left,midway,scale=0.7] {$(a,b)$};
\end{tric}
}

\def\labedgeb{
\begin{tric}
\draw (0,-0.5)--(0,0.5) (0.3,-0.5)--(0.3,0.5) (1.1,-0.5)--(1.1,0.5);
\filldraw[black] (0.5,0) circle (1pt) (0.7,0) circle (1pt) (0.9,0) circle (1pt) (1.9,0) circle (1pt) (2.1,0) circle (1pt) (2.3,0) circle (1pt);
\draw[double,thin] (1.4,-0.5)--(1.4,0.5) (1.7,-0.5)--(1.7,0.5) (2.5,-0.5)--(2.5,0.5);
\end{tric}
}
We will write
\[
\underline{s}_{(a,b)}:=\varpi_1^{\otimes a}\otimes \varpi_2 ^{\otimes b},
\]
and we will denote $\id_{\underline{s}_{(a,b)}}=\labedgeb \in End(\underline{s}_{(a,b)})$ as a labelled edge $\labedgea$. 
\end{notation}

\subsection{Kuperberg's Results on Equinumeration}
\label{equinemeration}

We recall the results of \cite{Kupe} which describe the relation between $\mathfrak{g}_2$ webs and representations of the quantum group associated to $\mathfrak{g}_2$. We will only work over the field $\mathbb{C}(q)$ where $q$ is either an indeterminant or a generic element of $\mathbb{C}^{\times}$. 

\begin{notation}
Let $\Phi$ be the root system of type $\mathfrak{g}_2$ with Weyl group $W$ and simple roots $\alpha_1$ and $\alpha_2$, where $\alpha_1$ is the short root. It follows that the positive roots are
\[
\Phi_+= \lbrace\alpha_1, 3\alpha_1 + \alpha_2, 2\alpha_1 + \alpha_2, 3\alpha_1 + 2\alpha_2, \alpha_1 + \alpha_2, \alpha_2\rbrace.
\]
Equip $\mathbb{Z}\Phi$ with the $W$ invariant symmetric form determined by 
\[
(\alpha_1, \alpha_1) = 2, \quad (\alpha_1, \alpha_2) = -3 = (\alpha_2, \alpha_1), \quad \text{and} \quad (\alpha_2, \alpha_2) = 6. 
\]

We write $X$ for the integral weight lattice and $X_+$ for the dominant integral weights. The fundamental weights are $\varpi_1= 2\alpha_1+ \alpha_2$ and $\varpi_2 = 3\alpha_1 + 2\alpha_2$. We may use the notation $(a, b)$ for $a\varpi_1 + b\varpi_2$, in particular $X_+= \lbrace (a, b) \ | \ a, b\ge 0\rbrace$. 
\end{notation}
 
 \begin{definition}
Let $\lambda, \mu\in X_+$. We define $\mu\le \lambda$ if $\lambda- \mu$ is a non-negative linear combination of positive roots. We also write $\mu< \lambda$ if $\mu\le \lambda$ and $\mu \ne \lambda$.  
\end{definition}

\begin{defn}\cite[Section 4.3]{JantzenQgps}
The algebra $\Uq$ is the $\mathbb{C}(q)$ algebra generated by elements $F_1$, $F_2$, $K_1^{\pm 1}$, $K_2^{\pm 1}$, $E_1$, and $E_2$ subject to the following relations:
\begin{itemize}
    \item $K_1K_2= K_2K_1$,
    \item $K_2E_2= q^{6}E_2K_2$, $K_2E_1= q^{-3}E_1K_2$, \ \ \ \ \ $K_1E_1= q^2E_1K_1$, $K_1E_2= q^{-3}K_1E_1$,
    \item $K_2F_2= q^{-6}F_2K_2$, $K_2F_1 = q^3F_1K_2$, \ \ \ \ \ $K_1F_1= q^{-2}F_1K_1$, $K_1F_2= q^3F_2K_1$,
    \item $E_1F_1 = F_1E_1 + \dfrac{K_1- K_1^{-1}}{q- q^{-1}}$, \ \ \ \ \ $E_2F_2= F_2E_2 + \dfrac{K_2- K_2^{-1}}{q^3- q^{-3}}$,
    \item $E_2^2E_1 + E_1E_2^2 = \dfrac{[6]}{[3]}E_2E_1E_2$, \ \ \ \ \  $F_2^2F_1 + F_1F_2^2 = \dfrac{[6]}{[3]}F_2F_1F_2$,
    \item $E_1^4E_2+ \dfrac{[4][3]}{[2]}E_1^2E_2E_1^2 + E_2E_1^4 = [4]E_1^3E_2E_1 + [4]E_1E_2E_1^3$,
    \item $F_1^4F_2+ \dfrac{[4][3]}{[2]}F_1^2F_2F_1^2 + F_2F_1^4 = [4]F_1^3F_2F_1 + [4]F_1F_2F_1^3$.
\end{itemize}
\end{defn}

The irreducible, finite dimensional, type-$\textbf{1}$ representations of $\Uq$ are in bijection with the finite dimensional irreducible representations of $\mathfrak{g}_2(\mathbb{C})$. For each $\lambda\in X_+$ we write $V(\lambda)$ for the $\Uq$ module which corresponds to the $\mathfrak{g}_2$ representation with highest weight $\lambda$. 

The algebra $\Uq$ is a Hopf algebra, so its representation category is a monoidal category. We are only interested in type-$\textbf{1}$ $\Uq$ modules, that is modules such that $K_{\alpha_1}$ and $K_{\alpha_2}$ act diagonalizably with eigenvalues in $+q^m$ for $m\in \mathbb{Z}$. It is not hard to see that the condition of being type-$\textbf{1}$ is closed under taking tensor product. 

\begin{notation}
We write $\Rep(\Uq)$ for the monoidal category of finite dimensional type-$\textbf{1}$ $\Uq$ modules. 
\end{notation}

The category $\Rep(U_q(\mathfrak{g}_2))$ is completely reducible \cite[Theorem 5.17]{JantzenQgps}. Moreover, we can determine how a module in $\Rep(U_q(\mathfrak{g}_2))$ decomposes by looking at its weight space decomposition.

The modules $V(\lambda)$ are type-$\textbf{1}$. Also, we have
\[
V(\lambda)\otimes V(\mu) \cong \bigoplus_{\nu \in X_+} V(\nu)^{\oplus m_{\nu}^{\lambda, \mu}},
\]
where the integers $m_{\nu}^{\lambda, \mu}$ are the same as those describing the tensor product decomposition of the analogous $\mathfrak{g}_2(\mathbb{C})$ modules. So the tensor product of type-$\textbf{1}$ modules are also type-$\textbf{1}$.

\begin{defn}[{\cite[Section 5.1]{JantzenQgps}}]
A module $W\in \Rep(\Uq)$ decomposes as a direct sum 
\[
W = \oplus_{\mu\in X}W_{\mu},
\]
where 
\[
W_{\mu} = \lbrace w\in W \ | \ K_1w = q^{(\alpha_1, \mu)}w, K_2w= q^{(\alpha_2, \mu)}w\rbrace. 
\]
We will call this direct sum decomposition the \textbf{weight space decomposition} of $W$, say that $W_{\mu}$ is the \textbf{$\mu$ weight space} of $W$, and call $w\in W_{\mu}$ a \textbf{weight vector} of weight $\mu$. We say that
\[
\wt W:= \lbrace \mu \ | \ W_{\mu}\ne 0\rbrace
\]
is the set of \textbf{weights} of $W$.
\end{defn}

\begin{remark}
It is well known that $\dim V(\lambda)_{\mu}$ is determined by the dimension of the $\mu$ weight space of the corresponding $\mathfrak{g}_2(\mathbb{C})$ module. In particular, $\dim V(\lambda)_{\lambda} = 1$ for all $\lambda\in X_+$. 
\end{remark}

\begin{notation}
Let $W$ be a module in $\Rep(U_q(\mathfrak{g}_2))$. For each $\lambda\in X_+$ there are non-negative integers $m_{\lambda}(W)$ such that
\[
W\cong \bigoplus_{\lambda\in X_+}V(\lambda)^{m_{\lambda}(W)}.
\]
We write $[W:V(\lambda)]:=m_{\lambda}(W)$ in this case. 
\end{notation}

\begin{defn}
The category of \textbf{fundamental representations}, $\Fund(\Uq)$ is the full monoidal subcategory of $\Rep(\Uq)$ generated by the objects $V(\varpi_1)$ and $V(\varpi_2)$. 
\end{defn}

\begin{remark}
The objects in the category $\Fund(U_q(\mathfrak{g}_2)$ are all isomorphic to iterated tensor products of fundamental representations. This includes the empty tensor product, which we take to be the trivial module, denoted by $\triv$. The category is $\mathbb{C}(q)$-linear additive, but is not closed under taking direct summands.
\end{remark}

\begin{definition}
Let $\unw$ be an object in $\DD$. Then $\unw = w_1w_2\ldots w_n$ for $w_i \in \lbrace \varpi_1, \varpi_2 \rbrace$. We define 
\[
V(\unw) := V(w_1)\otimes V(w_2)\otimes \ldots \otimes V(w_n).
\]
\end{definition}

\begin{remark}
Note that
\[
V(\unw)\cong \bigoplus_{\mu\in X_+} {V(\mu)^{\oplus m_{\mu}^{\unw}}}_{.}
\]
The integers $m_{\mu}^{\unw} :=m_{\mu}(V(\unw))= [V(\unw):V(\mu)]$ are the same as those describing the tensor product decomposition of the analogous $\mathfrak{g}_2(\mathbb{C})$ modules. 
\end{remark}

\begin{notation}
Given an object $\unw = w_1w_2\ldots w_n$, we write
\[
\wt \unw = \sum_{i=1}^n \wt w_i. 
\]
Note that $\wt \unw\in X_+$ for all $\unw$. 
\end{notation}

\begin{thm}[{\cite[Theorem 5.1]{Kupe}}]\label{T:Kupes-functor-Phi}
There is an essentially surjective monoidal functor 
\[
\Phi: \DD\rightarrow \Fund(\Uq)
\]
such that $\Phi(\varpi) = V(\varpi)$ for $\varpi\in \lbrace \varpi_1, \varpi_2\rbrace$. 
\end{thm}

\begin{thm}[{\cite[Theorem 6.10]{Kupe}}]\label{equivalencethm}
Let $\unw$ and $\unu$ be objects in $\DD$. Then
\[
\dim\Hom_{\DD}(\unw, \unu) = \dim\Hom_{\Fund(\Uq)}(V(\unw), V(\unu)),
\]
and it follows that the functor $\Phi$ is an equivalence of monoidal categories. 
\end{thm}

Recall that given a category $\mathcal{C}$, the Karoubi envelope of $\mathcal{C}$, is the category with objects: pairs $(X, e)$, where $X$ is an object in $\mathcal{C}$ and $e\in \End_{\mathcal{C}}(X)$ is an idempotent, and morphisms: triples $(e', f, e): (X, e)\rightarrow (Y, e')$, where $f:X\rightarrow Y$ is a morphism in $\mathcal{C}$ so that $e'\circ f\circ e= f$. Given a $\mathbb{C}(q)$-linear category $\mathcal{C}$, the additive envelope of $\mathcal{C}$ is the category with objects formal direct sums of objects in $\mathcal{C}$ and morphisms matrices of morphisms in $\mathcal{C}$. 

\begin{definition}
Let $\mathcal{C}$ be a $\mathbb{C}(q)$-linear category. Define the \emph{Karoubi completion} of $\mathcal{C}$ to be the additive envelope of the Karoubi envelope of $\mathcal{C}$.
\end{definition}

\begin{cor}\label{karoubi}
The functor $\Phi$ induces an equivalence of monoidal categories
\[
\Kar(\DD)\rightarrow \Rep(\Uq)
\]
such that $(\unw, e)\mapsto \im \Phi(e)$ and $(e', f, e): (\unw, e)\rightarrow (\unw, e')\mapsto \Phi(e'\circ f\circ e)$.
\end{cor}

\begin{proof}
Since every object in $\Rep(\Uq)$ is a direct sum of direct summands of objects in $\Fund(\Uq)$, this follows from $\Phi$ being an equivalence.
\end{proof}

\begin{lemma}\label{lowertermslemma}
Let $D\in \End_{\DD}(\unw)$ such that $\Phi(D)$ acts as zero on $V(\wt \unw)\smd V(\unw)$. Then we can write $D$ as linear combination
\[
D= \sum_i A_i\circ B_i,
\]
where $B_i\in \Hom_{\DD}(\unw, \unu_i)$ and $A_i\in \Hom_{\DD}(\unu_i, \unw)$ for some $\unu_i$ with $\wt \unu_i< \wt \unw$.
\end{lemma}
\begin{proof}
Since 
\[
V(\unw) = V(\wt \unw) \bigoplus_{\mu< \wt \unw}V(\mu)^{\oplus m_{\mu}^{\unw}},
\]
\\
our hypothesis on $\Phi(D)$ implies that we can write $\Phi(D) = \sum_i \iota_i\circ \pi_i$, where for each $i$ there is some $\mu_i< \wt \unw$ such that $\pi_i$ is a projection $V(\unw)\rightarrow V(\mu_i)$ and $\iota_i$ is an inclusion $V(\mu_i)\rightarrow V(\unw)$. 

For each $\mu_i$ fix an object $\unu_i$ in $\DD$ with $\wt \unu_i = \mu_i$. There is a projection $\gamma_i: V(\unu_i)\rightarrow V(\mu_i)$ and inclusion $\gamma^i: V(\mu_i)\rightarrow V(\unu_i)$ so that $\gamma_i\circ \gamma^i= \id_{V(\mu_i)}$. Then we can write 
\[
\iota_i\circ \pi_i = \iota_i \circ \id_{V(\mu_i)} \circ \pi_i = \iota_i \circ \gamma_i\circ \gamma^i \circ \pi_i.
\]
Thus, $\iota_i\circ \gamma_i\in \Hom_{\Uq}(V(\unu_i), V(\unw))$ and $\gamma^i\circ \pi_i\in \Hom_{\Uq}(V(\unw), V(\unu_i))$. The desired  result now follows from $\Phi$ being an equivalence. 
\end{proof}

\begin{defn}
The \textbf{neutral coefficient} of a diagram $D\in \End(\unw)$ is the scalar by which $\Phi(D)$ acts on the one dimensional weight space $V(\unw)_{\wt\unw}$. We write $\Phi(D)|_{V(\unw)_{\wt\unw}}= N_{D}\cdot \id$. 
\end{defn}

\begin{lemma}\label{factorizationlemma}
Let $D\in\End_{\DD}(\unw)$. Then we can express $D$ as a linear combination of diagrams 
\[
D= N_D\cdot \id_{\unw} + \sum_i A_i\circ B_i,
\]
where $B_i\in \Hom_{\DD}(\unw, \unu_i)$ and $A_i\in \Hom_{\DD}(\unu_i, \unw)$ for some $\unu_i$ with $\wt \unu_i< \wt \unw$.
\end{lemma}

\begin{proof}
Consider $\Phi(D)- N_D\id_{V(\unw)}\in \End(V(\unw))$. This endomorphism has $V(\unw)_{\wt\unw}$ in its kernel and therefore also acts as zero on $V(\wt \unw)\smd V(\unw)$. The desired result now follows from Lemma \ref{lowertermslemma}. 
\end{proof}

\subsection{Definition of clasps}

\def\defaa{
\begin{tric}
\draw[darkred,thick] (0,-1.2) rectangle (1.5,-0.8)
                      node[pos=0.5,black,scale=0.7] {$(a,b)$}; 
\draw[darkred,thick] (0,1.2) rectangle (1.5,0.8) 
                      node[pos=0.5,black,scale=0.7] {$(a,b)$};
\draw (0.75,-0.8)--(0.75,0.8) node[left,midway,scale=0.7] {$(a,b)$};
\draw (0.75,-1.2)--(0.75,-2) node[left,midway,scale=0.7] {$(a,b)$};
\draw (0.75,1.2)--(0.75,2) node[left,midway,scale=0.7] {$(a,b)$};
\end{tric}
}

\def\thindefaa{
\begin{tric}
\draw[darkred,thick] (0,-1.2) rectangle (1.2,-0.8)
                      node[pos=0.5,black,scale=0.7] {$(a,b)$}; 
\draw[darkred,thick] (0,1.2) rectangle (1.2,0.8) 
                      node[pos=0.5,black,scale=0.7] {$(a,b)$};
\draw (0.75,-0.8)--(0.75,0.8) node[left,midway,scale=0.7] {$(a,b)$};
\draw (0.75,-1.2)--(0.75,-2.5) node[left,midway,scale=0.7] {$(a,b)$};
\draw (0.75,1.2)--(0.75,2.5) node[left,midway,scale=0.7] {$(a,b)$};
\end{tric}
}

\def\defab{
\begin{tric}
\draw[darkred,thick] (0,-0.2) rectangle (1.5,0.2)
                      node[pos=0.5,black,scale=0.7] {$(a,b)$}; 
\draw (0.75,-0.2)--(0.75,-2) node[left,midway,scale=0.7] {$(a,b)$};
\draw (0.75,0.2)--(0.75,2) node[left,midway,scale=0.7] {$(a,b)$};
\end{tric}
}

\def\defabshort{
\begin{tric}
\draw[darkred,thick] (0,-0.2) rectangle (1.5,0.2)
                      node[pos=0.5,black,scale=0.7] {$(a,b)$}; 
\draw (0.75,-0.2)--(0.75,-1.5) node[left,midway,scale=0.7] {$(a,b)$};
\draw (0.75,0.2)--(0.75,1.5) node[left,midway,scale=0.7] {$(a,b)$};
\end{tric}
}

\def\defac{
\begin{tric}
\draw[darkred,thick] (0,-1.2) rectangle (1.5,-0.8)
                      node[pos=0.5,black,scale=0.7] {$(a,b)$}; 
\draw[darkred,thick] (0.75,1) ellipse (0.75 and 0.2)
                      node[black,scale=0.7] {$W$};
\draw (0.75,-0.8)--(0.75,0.8) node[left,midway,scale=0.7] {$(a,b)$};
\draw (0.75,-1.2)--(0.75,-2) node[left,midway,scale=0.7] {$(a,b)$};
\draw (0.75,1.2)--(0.75,2) node[left,midway,scale=0.7] {$(m,n)$};
\end{tric}
}

\def\markedstring{
\begin{tric}
\draw (0.75,-0.8)--(0.75,0.8) node[left,midway,scale=0.7] {$(a,b)$};
\end{tric}
}

\def\lowweb{
\begin{tric}
\draw[darkred,thick] (0.75,1) ellipse (0.75 and 0.2)
                      node[black,scale=0.7] {$W$};
\draw (0.75,0.8)--(0.75,0.4) node[left,midway,scale=0.7] {$(a,b)$};
\draw (0.75,1.2)--(0.75,1.6) node[left,midway,scale=0.7] {$(m,n)$};
\end{tric}                      
}

\def\neck{
\begin{tric}
\draw [darkred,thick] (0.75,1) ellipse (1.1 and 0.5)
                      node[black] {${\  }^1 D_{(c,d)}^{i{(c,d)}}$}
                      (0.75,-0.5) ellipse (1.1 and 0.5)
                      node[black] {${\  }^2 D_{(c,d)}^{i{(c,d)}}$};
\draw (0.75,0)--(0.75,0.5) node[midway,left,scale=0.7] {$(c,d)$}
      (0.75,-1)--(0.75,-1.5) node[midway,left,scale=0.7] {$(a,b)$}
       (0.75,1.5)--(0.75,2) node[midway,left,scale=0.7] {$(a,b)$};
\end{tric}
}

\def\idab{
\begin{tric}
\draw (0.75,-1.5)--(0.75,1.5) node[left,midway,scale=0.7] {$(a,b)$};
\end{tric}
}

\def\neckonclasp{
\begin{tric}
\draw [darkred,thick] (0.75,1.2) ellipse (1.1 and 0.5)
                      node[black] {${\  }^1 D_{(c,d)}^{i{(c,d)}}$}
                      (0.75,-0.3) ellipse (1.1 and 0.5)
                      node[black] {${\  }^2 D_{(c,d)}^{i{(c,d)}}$};
\draw (0.75,0.2)--(0.75,0.7) node[midway,left,scale=0.7] {$(c,d)$}
      (0.75,-0.8)--(0.75,-1.3) node[midway,left,scale=0.7] {$(a,b)$}
       (0.75,1.7)--(0.75,2.2) node[midway,left,scale=0.7] {$(a,b)$}
        (0.75,-1.7)--(0.75,-2.2) node[midway,left,scale=0.7] {$(a,b)$};
\draw[darkred,thick](0,-1.3) rectangle (1.5,-1.7) 
                      node[pos=0.5,black,scale=0.7] {$(a,b)$};
\end{tric}
}

\def\idabonclasp{
\begin{tric}
\draw (0.75,-1)--(0.75,1.5) node[left,midway,scale=0.7] {$(a,b)$}
      (0.75,-1.4)--(0.75,-2.4) node[left,midway,scale=0.7] {$(a,b)$};
\draw[darkred,thick](0,-1) rectangle (1.5,-1.4) 
                      node[pos=0.5,black,scale=0.7] {$(a,b)$};
\end{tric}
}

\def\necktop{
\begin{tric}
\draw [darkred,thick] (0.75,1) ellipse (1.1 and 0.5)
                      node[black] {${\  }^1 D_{(c,d)}^{i{(c,d)}}$};
\draw (0.75,1.5)--(0.75,1.9) node[left,midway,scale=0.7] {$(a,b)$};
\draw (0.75,0.5)--(0.75,0.1) node[left,midway,scale=0.7] {$(c,d)$};
\end{tric}
}

\def\neckbutt{
\begin{tric}
\draw [darkred,thick] (0.75,1) ellipse (1.1 and 0.5)
                      node[black] {${\  }^2 D_{(c,d)}^{i{(c,d)}}$};
\draw (0.75,1.5)--(0.75,1.9) node[left,midway,scale=0.7] {$(c,d)$};
\draw (0.75,0.5)--(0.75,0.1) node[left,midway,scale=0.7] {$(a,b)$};
\end{tric}
}

\def\sameclaspa{
\begin{tric}
\draw[darkred,thick](0,-1) rectangle (1.5,-1.4) 
                      node[pos=0.5,black,scale=0.7] {Clasp\ 1};
\draw (0.75,-1)--(0.75,-0.6) node[left,midway,scale=0.7] {$(a,b)$};
\draw (0.75,-1.4)--(0.75,-1.8) node[left,midway,scale=0.7] {$(a,b)$};
\end{tric}
}

\def\sameclaspb{
\begin{tric}
\draw[darkred,thick](0,-1) rectangle (1.5,-1.4) 
                      node[pos=0.5,black,scale=0.7] {Clasp\ 2};
\draw (0.75,-1)--(0.75,-0.6) node[left,midway,scale=0.7] {$(a,b)$};
\draw (0.75,-1.4)--(0.75,-1.8) node[left,midway,scale=0.7] {$(a,b)$};
\end{tric}
}

\def\sameclaspc{
\begin{tric}
\draw[darkred,thick](0,0.3) rectangle (1.5,0.7) 
                      node[pos=0.5,black,scale=0.7] {Clasp\ 1}
                     (0,-0.3) rectangle (1.5,-0.7) 
                      node[pos=0.5,black,scale=0.7] {Clasp\ 2};
\draw (0.75,-0.3)--(0.75,0.3) node[left,midway,scale=0.7] {$(a,b)$}
      (0.75,-0.7)--(0.75,-1.2) node[left,midway,scale=0.7] {$(a,b)$}
      (0.75,0.7)--(0.75,1.2) node[left,midway,scale=0.7] {$(a,b)$};
\end{tric}
}

\def\sameclaspd{
\begin{tric}
\draw (0.75,-1)--(0.75,1.5) node[left,midway,scale=0.7] {$(a,b)$}
      (0.75,-1.4)--(0.75,-2.4) node[left,midway,scale=0.7] {$(a,b)$};
\draw[darkred,thick](0,-1) rectangle (1.5,-1.4) 
                      node[pos=0.5,black,scale=0.7] {Clasp 2};
\end{tric}
}

\def\sameclaspe{
\begin{tric}
\draw [darkred,thick] (0.75,1.2) ellipse (1.1 and 0.5)
                node[black] {${\  }^1 {\overline{D}}_{(c,d)}^{j{(c,d)}}$}
                      (0.75,-0.3) ellipse (1.1 and 0.5)
                node[black] {${\  }^2 {\overline{D}}_{(c,d)}^{j{(c,d)}}$};
\draw (0.75,0.2)--(0.75,0.7) node[midway,left,scale=0.7] {$(c,d)$}
      (0.75,-0.8)--(0.75,-1.3) node[midway,left,scale=0.7] {$(a,b)$}
       (0.75,1.7)--(0.75,2.2) node[midway,left,scale=0.7] {$(a,b)$}
        (0.75,-1.7)--(0.75,-2.2) node[midway,left,scale=0.7] {$(a,b)$};
\draw[darkred,thick](0,-1.3) rectangle (1.5,-1.7) 
                      node[pos=0.5,black,scale=0.7] {Clasp 2};
\end{tric}
}

\def\sameclaspf{
\begin{tric}
\draw (0.75,-1)--(0.75,1.5) node[left,midway,scale=0.7] {$(a,b)$}
      (0.75,1.9)--(0.75,2.9) node[left,midway,scale=0.7] {$(a,b)$};
\draw[darkred,thick](0,1.5) rectangle (1.5,1.9) 
                      node[pos=0.5,black,scale=0.7] {Clasp 1};
\end{tric}
}

\def\sameclaspg{
\begin{tric}
\draw [darkred,thick] (0.75,1.2) ellipse (1.1 and 0.5)
                      node[black] {${\  }^1 {\overline{\overline{D}}}_{(c,d)}^{k{(c,d)}}$}
                      (0.75,-0.3) ellipse (1.1 and 0.5)
                      node[black] {${\  }^2 {\overline{\overline{D}}}_{(c,d)}^{k{(c,d)}}$};
\draw (0.75,0.2)--(0.75,0.7) node[midway,left,scale=0.7] {$(c,d)$}
      (0.75,-0.8)--(0.75,-1.3) node[midway,left,scale=0.7] {$(a,b)$}
       (0.75,1.7)--(0.75,2.2) node[midway,left,scale=0.7] {$(a,b)$}
        (0.75,2.6)--(0.75,3.1) node[midway,left,scale=0.7] {$(a,b)$};
\draw[darkred,thick](0,2.2) rectangle (1.5,2.6) 
                      node[pos=0.5,black,scale=0.7] {Clasp 1};
\end{tric}
}

\def\idab{
\begin{tric}
\draw (0.75,-1.5)--(0.75,1.5) node[left,midway,scale=0.7] {$(a,b)$};
\end{tric}
}

\def\sameclasph{
\begin{tric}
\draw (0.75,-0.2)--(0.75,-1.5) node[left,midway,scale=0.7] {$(a,b)$}
      (0.75,0.2)--(0.75,1.5) node[left,midway,scale=0.7] {$(a,b)$};
\draw[darkred,thick](0,-0.2) rectangle (1.5,0.2) 
                      node[pos=0.5,black,scale=0.7] {Clasp 1};
\end{tric}
}

\def\sameclaspi{
\begin{tric}
\draw (0.75,-1.5)--(0.75,1.5) node[left,midway,scale=0.7] {$(a,b)$};
\end{tric}
}

\def\sameclaspj{
\begin{tric}
\draw [darkred,thick] (0.75,1.2) ellipse (1.1 and 0.5)
                      node[black] {${\  }^1 
                      {\overline{D}}_{(c,d)}^{j{(c,d)}}$}
                      (0.75,-0.3) ellipse (1.1 and 0.5)
                      node[black] {${\  }^2 {\overline{D}}_{(c,d)}^{j{(c,d)}}$};
\draw (0.75,0.2)--(0.75,0.7) node[midway,left,scale=0.7] {$(c,d)$}
      (0.75,-0.8)--(0.75,-1.3) node[midway,left,scale=0.7] {$(a,b)$}
       (0.75,1.7)--(0.75,2.2) node[midway,left,scale=0.7] {$(a,b)$};
\end{tric}
}

\def\sameclaspk{
\begin{tric}
\draw (0.75,-0.2)--(0.75,-1.5) node[left,midway,scale=0.7] {$(a,b)$}
      (0.75,0.2)--(0.75,1.5) node[left,midway,scale=0.7] {$(a,b)$};
\draw[darkred,thick](0,-0.2) rectangle (1.5,0.2) 
                      node[pos=0.5,black,scale=0.7] {Clasp 2};
\end{tric}
}

\def\sameclaspl{
\begin{tric}
\draw (0.75,-1.5)--(0.75,1.5) node[left,midway,scale=0.7] {$(a,b)$};
\end{tric}
}

\def\sameclaspm{
\begin{tric}
\draw [darkred,thick] (0.75,1.2) ellipse (1.1 and 0.5)
                      node[black] {${\  }^1 {\overline{\overline{D}}}_{(c,d)}^{k{(c,d)}}$}
                      (0.75,-0.3) ellipse (1.1 and 0.5)
                      node[black] {${\  }^2 {\overline{\overline{D}}}_{(c,d)}^{k{(c,d)}}$};
\draw (0.75,0.2)--(0.75,0.7) node[midway,left,scale=0.7] {$(c,d)$}
      (0.75,-0.8)--(0.75,-1.3) node[midway,left,scale=0.7] {$(a,b)$}
       (0.75,1.7)--(0.75,2.2) node[midway,left,scale=0.7] {$(a,b)$};
\end{tric}
}

\def\Hmovea{
\begin{tric}
\draw (0.5,-0.2)--(0.5,-1)  (0.5,0.2)--(0.5,1) ;
\draw[double,thin] (1,-0.2)--(1,-1)  (1,0.2)--(1,1) ;
\draw[darkred,thick](-0.5,-0.2) rectangle (2,0.2) 
                      node[pos=0.5,black,scale=0.7] {$(a,b)$};
\draw (0,-0.2)--(0,-1) node[midway,left,scale=0.7] {$(a_1,b_1)$}                      (1.5,-0.2)--(1.5,-1) node[midway,right,scale=0.7] {$(a_2,b_2)$}
       (0,0.2)--(0,1) node[midway,left,scale=0.7] {$(a_1,b_1)$}               (1.5,0.2)--(1.5,1) node[midway,right,scale=0.7] {$(a_2,b_2)$};
\end{tric}
}

\def\Hmoveb{
\begin{tric}
\draw [double,thin] (0.5,-0.2)--(0.5,-1)  (0.5,0.2)--(0.5,1) ;
\draw (1,-0.2)--(1,-1)  (1,0.2)--(1,1) ;
\draw[darkred,thick](-0.5,-0.2) rectangle (2,0.2) 
                      node[pos=0.5,black,scale=0.7] {$(a,b)$};
\draw (0,-0.2)--(0,-1) node[midway,left,scale=0.7] {$(a_1,b_1)$}                      (1.5,-0.2)--(1.5,-1) node[midway,right,scale=0.7] {$(a_2,b_2)$}
       (0,0.2)--(0,1) node[midway,left,scale=0.7] {$(a_1,b_1)$}               (1.5,0.2)--(1.5,1) node[midway,right,scale=0.7] {$(a_2,b_2)$};
\end{tric}
}

\def\Hmovec{
\begin{tric}
\draw (0.5,-0.2)--(0.5,-0.6)  (0.5,0.2)--(0.5,0.6)
       (1,-1)--(1,-0.6)      (1,1)--(1,0.6)   
       (0.5,-0.6)--(1,-0.6) (0.5,0.6)--(1,0.6);
\draw[double,thin] (1,-0.2)--(1,-0.6)  (1,0.2)--(1,0.6) 
                   (0.5,-0.6)--(0.5,-1)  (0.5,0.6)--(0.5,1);
\draw[darkred,thick](-0.5,-0.2) rectangle (2,0.2) 
                      node[pos=0.5,black,scale=0.7] {$(a,b)$};
\draw (0,-0.2)--(0,-1) node[midway,left,scale=0.7] {$(a_1,b_1)$}                      (1.5,-0.2)--(1.5,-1) node[midway,right,scale=0.7] {$(a_2,b_2)$}
       (0,0.2)--(0,1) node[midway,left,scale=0.7] {$(a_1,b_1)$}               (1.5,0.2)--(1.5,1) node[midway,right,scale=0.7] {$(a_2,b_2)$};
\end{tric}
}

\def\Hmoved{
\begin{tric}
\draw (0.5,-0.4)--(0.5,-0.9)  (0.5,0.2)--(0.5,0.6)
             (1,1)--(1,0.6)   
       (0.5,0.6)--(1,0.6)
       (0.5,1)--(1,1) (0.5,1)--(0.5,1.4);
       
\draw[double,thin] (1,-0.4)--(1,-0.9)  (1,0.2)--(1,0.6) 
                    (0.5,0.6)--(0.5,1)
                   (1,1)--(1,1.4);
\draw[darkred,thick]  (-0.5,-0.4) rectangle (2,0.2) 
                      node[pos=0.5,black] {$C_{\unw}$};
\filldraw[black] 
(-0.2,-0.7) circle (1pt) (0,-0.7) circle (1pt) (0.2,-0.7) circle (1pt)
(1.3,-0.7) circle (1pt) (1.5,-0.7) circle (1pt) (1.7,-0.7) circle (1pt)
(-0.2,0.8) circle (1pt) (0,0.8) circle (1pt) (0.2,0.8) circle (1pt)
(1.3,0.8) circle (1pt) (1.5,0.8) circle (1pt) (1.7,0.8) circle (1pt);

\end{tric}
}

\def\Hmovee{
\begin{tric}
\draw (0.5,-0.4)--(0.5,-0.9)  (0.5,0.2)--(0.5,0.6)
             (1,1)--(1,0.6)   
       (0.5,0.6)--(1,0.6)
       (0.5,1)--(1,1) (0.5,1)--(0.5,1.4)
       (0.5,0.6)--(0.5,1);
       
\draw[double,thin] (1,-0.4)--(1,-0.9)  (1,0.2)--(1,0.6) 
                   (1,1)--(1,1.4);
\draw[darkred,thick]  (-0.5,-0.4) rectangle (2,0.2) 
                      node[pos=0.5,black] {$C_{\unw}$};
\filldraw[black] 
(-0.2,-0.7) circle (1pt) (0,-0.7) circle (1pt) (0.2,-0.7) circle (1pt)
(1.3,-0.7) circle (1pt) (1.5,-0.7) circle (1pt) (1.7,-0.7) circle (1pt)
(-0.2,0.8) circle (1pt) (0,0.8) circle (1pt) (0.2,0.8) circle (1pt)
(1.3,0.8) circle (1pt) (1.5,0.8) circle (1pt) (1.7,0.8) circle (1pt);

\end{tric}
}

\def\Hmovef{
\begin{tric}
\draw (0.5,-0.4)--(0.5,-0.9)  
             (1,1)--(1,0.6)   
       (0.5,0.2)..controls(0.5,0.6)..(1,0.6)
         (1,1)..controls(0.5,1)..(0.5,1.4);
       
\draw[double,thin] (1,-0.4)--(1,-0.9)  (1,0.2)--(1,0.6) 
                   (1,1)--(1,1.4);
\draw[darkred,thick]  (-0.5,-0.4) rectangle (2,0.2) 
                      node[pos=0.5,black] {$C_{\unw}$};
\filldraw[black] 
(-0.2,-0.7) circle (1pt) (0,-0.7) circle (1pt) (0.2,-0.7) circle (1pt)
(1.3,-0.7) circle (1pt) (1.5,-0.7) circle (1pt) (1.7,-0.7) circle (1pt)
(-0.2,0.8) circle (1pt) (0,0.8) circle (1pt) (0.2,0.8) circle (1pt)
(1.3,0.8) circle (1pt) (1.5,0.8) circle (1pt) (1.7,0.8) circle (1pt);

\end{tric}
}

\def\Hmoveg{
\begin{tric}
\draw  (0.5,-0.4)--(0.5,-0.9)  
         (1,1)--(1,0.6)   

(0.5,0.2)--(0.5,0.8)--(0.5,1.4) (1,1)--(0.8,0.8)--(1,0.6) (0.8,0.8)--(0.5,0.8);
       
\draw[double,thin] (1,-0.4)--(1,-0.9)  (1,0.2)--(1,0.6) 
                   (1,1)--(1,1.4);
\draw[darkred,thick]  (-0.5,-0.4) rectangle (2,0.2) 
                      node[pos=0.5,black] {$C_{\unw}$};
\filldraw[black] 
(-0.2,-0.7) circle (1pt) (0,-0.7) circle (1pt) (0.2,-0.7) circle (1pt)
(1.3,-0.7) circle (1pt) (1.5,-0.7) circle (1pt) (1.7,-0.7) circle (1pt)
(-0.2,0.8) circle (1pt) (0,0.8) circle (1pt) (0.2,0.8) circle (1pt)
(1.3,0.8) circle (1pt) (1.5,0.8) circle (1pt) (1.7,0.8) circle (1pt);

\end{tric}
}

\def\Hmoveh{
\begin{tric}
\draw  (0.5,-0.4)--(0.5,-0.9)  (0.5,0.2)--(0.5,1.4)     (1,1)--(1,0.6)   
        (1,0.6)..controls(0.6,0.6)and(0.6,1)..(1,1);
       
\draw[double,thin] (1,-0.4)--(1,-0.9)  (1,0.2)--(1,0.6) 
                   (1,1)--(1,1.4);
\draw[darkred,thick]  (-0.5,-0.4) rectangle (2,0.2) 
                      node[pos=0.5,black] {$C_{\unw}$};
\filldraw[black] 
(-0.2,-0.7) circle (1pt) (0,-0.7) circle (1pt) (0.2,-0.7) circle (1pt)
(1.3,-0.7) circle (1pt) (1.5,-0.7) circle (1pt) (1.7,-0.7) circle (1pt)
(-0.2,0.8) circle (1pt) (0,0.8) circle (1pt) (0.2,0.8) circle (1pt)
(1.3,0.8) circle (1pt) (1.5,0.8) circle (1pt) (1.7,0.8) circle (1pt);

\end{tric}
}

\def\Hmovei{
\begin{tric}
\draw  (0.5,-0.4)--(0.5,-0.9)  (0.5,0.2)--(0.5,1.4) ;
       
\draw[double,thin] (1,-0.4)--(1,-0.9)  (1,0.2)--(1,1.4);
\draw[darkred,thick]  (-0.5,-0.4) rectangle (2,0.2) 
                      node[pos=0.5,black] {$C_{\unw}$};
\filldraw[black] 
(-0.2,-0.7) circle (1pt) (0,-0.7) circle (1pt) (0.2,-0.7) circle (1pt)
(1.3,-0.7) circle (1pt) (1.5,-0.7) circle (1pt) (1.7,-0.7) circle (1pt)
(-0.2,0.8) circle (1pt) (0,0.8) circle (1pt) (0.2,0.8) circle (1pt)
(1.3,0.8) circle (1pt) (1.5,0.8) circle (1pt) (1.7,0.8) circle (1pt);

\end{tric}
}

\def\Hmovej{
\begin{tric}
\draw (0.5,-0.2)--(0.5,-0.6)  (0.5,0.2)--(0.5,0.6)
       (1,-1)--(1,-0.6)      (1,1)--(1,0.6)   
       (0.5,-0.6)--(1,-0.6) (0.5,0.6)--(1,0.6);
\draw[double,thin] (1,-0.2)--(1,-0.6)  (1,0.2)--(1,0.6) 
                   (0.5,-0.6)--(0.5,-1)  (0.5,0.6)--(0.5,1);
\draw[darkred,thick](-0.5,-0.2) rectangle (2,0.2) 
                      node[pos=0.5,black,scale=0.7] {$(a,b)$};
\draw (0,-0.2)--(0,-1) node[midway,left,scale=0.7] {$(a_1,b_1)$}                      (1.5,-0.2)--(1.5,-1) node[midway,right,scale=0.7] {$(a_2,b_2)$}
       (0,0.2)--(0,1) node[midway,left,scale=0.7] {$(a_1,b_1)$}               (1.5,0.2)--(1.5,1) node[midway,right,scale=0.7] {$(a_2,b_2)$}
         (0.75,1.45)--(0.75,1.85) node[midway,left,scale=0.7] {$(m,n)$};
\filldraw[darkred,thick,fill=white] (0.75,1.2) ellipse (1.25 and 0.25)   
                      node[black,scale=0.7] {$W$} ;
\end{tric}
}

\def\Hmovek{
\begin{tric}
\draw  (0.5,0.2)--(0.5,0.6)     (1,1)--(1,0.6)   (0.5,0.6)--(1,0.6);
\draw[double,thin]   (1,0.2)--(1,0.6)  (0.5,0.6)--(0.5,1);
\draw  (0,0.2)--(0,1) node[midway,left,scale=0.7] {$(a_1,b_1)$}  
       (1.5,0.2)--(1.5,1) node[midway,right,scale=0.7] {$(a_2,b_2)$}
         (0.75,1.45)--(0.75,1.85) node[midway,left,scale=0.7] {$(m,n)$};
\filldraw[darkred,thick,fill=white] (0.75,1.2) ellipse (1.25 and 0.25)   
                      node[black,scale=0.7] {$W$} ;
\end{tric}
}

\def\Hmovel{
\begin{tric}
\draw [double,thin] (0.5,-0.2)--(0.5,-0.6)  (0.5,0.2)--(0.5,0.6) 
                     (1,-0.6)--(1,-1) (1,0.6)--(1,1);
\draw (1,-0.2)--(1,-0.6)  (1,0.2)--(1,0.6) 
       (0.5,0.6)--(0.5,1) (0.5,-0.6)--(0.5,-1)
       (0.5,0.6)--(1,0.6) (0.5,-0.6)--(1,-0.6) ;
\draw[darkred,thick](-0.5,-0.2) rectangle (2,0.2) 
                      node[pos=0.5,black,scale=0.7] {$(a,b)$};
\draw (0,-0.2)--(0,-1) node[midway,left,scale=0.7] {$(a_1,b_1)$}                      (1.5,-0.2)--(1.5,-1) node[midway,right,scale=0.7] {$(a_2,b_2)$}
       (0,0.2)--(0,1) node[midway,left,scale=0.7] {$(a_1,b_1)$}               (1.5,0.2)--(1.5,1) node[midway,right,scale=0.7] {$(a_2,b_2)$};
\end{tric}
}

\def\Ha{
\begin{tric}
\draw (0,-0.4)--(0,0)--(0.7,0)--(0.7,0.4);
\draw[double,thin] (0,0)--(0,0.4) (0.7,0)--(0.7,-0.4);
\end{tric}
}

\def\Hb{
\begin{tric}
\draw [double,thin](0,-0.4)--(0,0)  (0.7,0)--(0.7,0.4);
\draw (0,0.4)--(0,0)--(0.7,0)--(0.7,-0.4);
\end{tric}
}

\def\neckcollar{
\begin{tric}
\draw [darkred,thick] (0.75,1.2) ellipse (1.4 and 0.5)
                      node[black] {${\  }^1 F_{(a,b)(c,d)}^{j{(c,d)}}$}
                      (0.75,-1.2) ellipse (1.4 and 0.5)
                      node[black] {${\  }^2 F_{(a,b)(c,d)}^{j{(c,d)}}$}
                     (0,0.2)rectangle(1.5,-0.2)node[pos=0.5,black,scale=0.7] {$(c,d)$};
\draw (0.75,0.2)--(0.75,0.7) node[midway,left,scale=0.7] {$(c,d)$}
       (0.75,-0.2)--(0.75,-0.7) node[midway,left,scale=0.7] {$(c,d)$}
      (0.75,-1.7)--(0.75,-2.2) node[midway,left,scale=0.7] {$(a,b)$}
       (0.75,1.7)--(0.75,2.2) node[midway,left,scale=0.7] {$(a,b)$};
\end{tric}
}

\def\neckcollara{
\begin{tric}
\draw [darkred,thick] (0.75,0) ellipse (1.4 and 0.5)
                      node[black] {${\  }^1 F_{(a,b)(c,d)}^{j{(c,d)}}$};
\draw (0.75,0.5)--(0.75,1) node[midway,left,scale=0.7] {$(a,b)$}
      (0.75,-0.5)--(0.75,-1) node[midway,left,scale=0.7] {$(c,d)$};
\end{tric}
}

\def\neckcollarb{
\begin{tric}
\draw [darkred,thick] (0.75,0) ellipse (1.4 and 0.5)
                      node[black] {${\  }^2 F_{(a,b)(c,d)}^{j{(c,d)}}$};
\draw (0.75,0.5)--(0.75,1) node[midway,left,scale=0.7] {$(c,d)$}
      (0.75,-0.5)--(0.75,-1) node[midway,left,scale=0.7] {$(a,b)$};
\end{tric}
}

\def\neckcollarc{
\begin{tric}
\draw[darkred,thick] (0,-0.2) rectangle (1.5,0.2)
                      node[pos=0.5,black,scale=0.7] {$(1,0)$}; 
\draw (0.75,-0.2)--(0.75,-1) (0.75,0.2)--(0.75,1) ;
\end{tric}
}

\def\neckcollard{
\begin{tric}
\draw (0.75,-1)--(0.75,1) ;
\end{tric}
}

\def\neckcollare{
\begin{tric}
\draw[darkred,thick] (0,-0.2) rectangle (1.5,0.2)
                      node[pos=0.5,black,scale=0.7] {$(0,1)$}; 
\draw[double,thin] (0.75,-0.2)--(0.75,-1)  (0.75,0.2)--(0.75,1) ;
\end{tric}
}

\def\neckcollarf{
\begin{tric}
\draw[double,thin] (0.75,-1)--(0.75,1) ;
\end{tric}
}

\def\neckcollarg{
\begin{tric}
\draw(0.75,-2)--(0.75,2) node[midway,left,scale=0.7] {$(c,d)$} ;
\end{tric}
}

\def\neckcollarh{
\begin{tric}
 \draw (0.75,0.2)--(0.75,2) node[midway,left,scale=0.7] {$(c,d)$}
        (0.75,-0.2)--(0.75,-2) node[midway,left,scale=0.7] {$(c,d)$};
\draw[darkred,thick] (0,-0.2) rectangle (1.5,0.2)
                      node[pos=0.5,black,scale=0.7] {$(c,d)$}; 
\end{tric}
}

\def\neckcollari{
\begin{tric}
\draw [darkred,thick] (0.75,1.2) ellipse (1.4 and 0.5)
                      node[black] {${\  }^1 F_{(c,d)(e,f)}^{j{(e,f)}}$}
                      (0.75,-1.2) ellipse (1.4 and 0.5)
                      node[black] {${\  }^2 F_{(c,d)(e,f)}^{j{(e,f)}}$}
                     (0,0.2)rectangle(1.5,-0.2)node[pos=0.5,black,scale=0.7] {$(e,f)$};
\draw (0.75,0.2)--(0.75,0.7) node[midway,left,scale=0.7] {$(e,f)$}
       (0.75,-0.2)--(0.75,-0.7) node[midway,left,scale=0.7] {$(e,f)$}
      (0.75,-1.7)--(0.75,-2.2) node[midway,left,scale=0.7] {$(c,d)$}
       (0.75,1.7)--(0.75,2.2) node[midway,left,scale=0.7] {$(c,d)$};
\end{tric}
}

\def\neckcollarj{
\begin{tric}
\draw(0.75,-1)--(0.75,1) node[midway,left,scale=0.7] {$(c,d)$} ;
\end{tric}
}

\def\neckcollark{
\begin{tric}
\draw [darkred,thick] (0.75,1.2) ellipse (1.1 and 0.5)
                      node[black] {${\  }^1 D_{(c,d)}^{i{(c,d)}}$}
                      (0.75,-1.2) ellipse (1.1 and 0.5)
                      node[black] {${\  }^2 D_{(c,d)}^{i{(c,d)}}$}
                      (0,0.2)rectangle(1.5,-0.2)
                      node[black,pos=0.5,scale=0.7]{$(c,d)$};
\draw (0.75,0.2)--(0.75,0.7) node[midway,left,scale=0.7] {$(c,d)$}
      (0.75,-0.2)--(0.75,-0.7) node[midway,left,scale=0.7] {$(c,d)$}
      (0.75,-1.7)--(0.75,-2.2) node[midway,left,scale=0.7] {$(a,b)$}
       (0.75,1.7)--(0.75,2.2) node[midway,left,scale=0.7] {$(a,b)$};
\end{tric}
}

\def\neckcollarl{
\begin{tric}
\draw [darkred,thick] (0.75,1.2) ellipse (1.4 and 0.5)
                      node[black] {${\  }^1 F_{(c,d)(e,f)}^{j{(e,f)}}$}
                      (0.75,-1.2) ellipse (1.4 and 0.5)
                      node[black] {${\  }^2 F_{(c,d)(e,f)}^{j{(e,f)}}$}
                     (0,0.2)rectangle(1.5,-0.2)node[pos=0.5,black,scale=0.7] {$(e,f)$}
                      (0.75,2.7) ellipse (1.1 and 0.5)
                      node[black] {${\  }^1 D_{(c,d)}^{i{(c,d)}}$}
                      (0.75,-2.7) ellipse (1.1 and 0.5)
                      node[black] {${\  }^2 D_{(c,d)}^{i{(c,d)}}$};
\draw (0.75,0.2)--(0.75,0.7) node[midway,left,scale=0.7] {$(e,f)$}
       (0.75,-0.2)--(0.75,-0.7) node[midway,left,scale=0.7] {$(e,f)$}
      (0.75,-1.7)--(0.75,-2.2) node[midway,left,scale=0.7] {$(c,d)$}
       (0.75,1.7)--(0.75,2.2) node[midway,left,scale=0.7] {$(c,d)$}
       (0.75,-3.7)--(0.75,-3.2) node[midway,left,scale=0.7] {$(a,b)$}
       (0.75,3.7)--(0.75,3.2) node[midway,left,scale=0.7] {$(a,b)$};
\end{tric}
}

\def\defad{
\begin{tric}
\draw[darkred,thick] (0,-1.3) rectangle (1.5,-1.7)
                      node[pos=0.5,black,scale=0.7] {$(a,b)$}; 
\draw[darkred,thick] (0,1.3) rectangle (1.5,1.7) 
                      node[pos=0.5,black,scale=0.7] {$(c,d)$};
\draw[darkred,thick] (0.75,0) ellipse (0.75 and 0.2)
                      node[black,scale=0.7] {$D$};
\draw (0.75,-1.3)--(0.75,-0.2) node[left,midway,scale=0.7] {$(a,b)$};
\draw (0.75,1.3)--(0.75,0.2) node[left,midway,scale=0.7] {$(c,d)$};
\draw (0.75,-1.7)--(0.75,-2.5) node[left,midway,scale=0.7] {$(a,b)$};
\draw (0.75,1.7)--(0.75,2.5) node[left,midway,scale=0.7] {$(c,d)$};
\end{tric}
}

\def\diagram{
\begin{tric}
\draw[darkred,thick] (0.75,0) ellipse (0.75 and 0.2)
                      node[black,scale=0.7] {$D$};
\draw (0.75,0.2)--(0.75,0.6) node[left,midway,scale=0.7] {$(c,d)$};
\draw (0.75,-0.2)--(0.75,-0.6) node[left,midway,scale=0.7] {$(a,b)$};
\end{tric}
}

\def\lowclasp{
\begin{tric}
\draw[darkred,thick] (0,1.3) rectangle (1.5,1.7) 
                      node[pos=0.5,black,scale=0.7] {$(c,d)$};
\draw (0.75,1.3)--(0.75,0.9) node[left,midway,scale=0.7] {$(c,d)$};
\draw (0.75,1.7)--(0.75,2.1) node[left,midway,scale=0.7] {$(c,d)$};
\end{tric}
}

\def\samedefa{
\begin{tric}
\draw[darkred,thick] 
      (0.75,0.6)ellipse(0.75 and 0.3) node[black,scale=0.7]{${\ }^1 D$} 
      (0.75,-0.6)ellipse(0.75 and 0.3) node[black,scale=0.7]{${\ }^2 D$};
\draw (0.75,-0.3)--(0.75,0.3) node[black,scale=0.7,left,midway]{$(c,d)$} 
      (0.75,-0.9)--(0.75,-1.4) node[black,scale=0.7,left,midway]{$(a,b)$} 
      (0.75,0.9)--(0.75,1.4) node[black,scale=0.7,left,midway]{$(m,n)$} ;
\end{tric}
}

\def\samedefb{
\begin{tric}
\draw[darkred,thick] 
      (0.75,0.6)ellipse(0.75 and 0.3) node[black,scale=0.7]{${\ }^1 D$} 
      (0.75,-0.6)ellipse(0.75 and 0.3) node[black,scale=0.7]{${\ }^2 D$}
      (0,-1.4)rectangle(1.5,-1.8) node[pos=0.5,black,scale=0.7]{$(a,b)$};
\draw (0.75,-0.3)--(0.75,0.3) node[black,scale=0.7,left,midway]{$(c,d)$} 
      (0.75,-0.9)--(0.75,-1.4) node[black,scale=0.7,left,midway]{$(a,b)$} 
      (0.75,0.9)--(0.75,1.4) node[black,scale=0.7,left,midway]{$(m,n)$}
      (0.75,-1.8)--(0.75,-2.3)node[black,scale=0.7,left,midway]{$(a,b)$} ;
\end{tric}
}

\def\samedefc{
\begin{tric}
\draw [darkred,thick] (0.75,1.2) ellipse (0.75 and 0.3)
                      node[black,scale=0.7] {${\  }^1 D$}
                      (0.75,-1.2) ellipse (0.75 and 0.3)
                      node[black,scale=0.7] {${\  }^2 D$}
                      (0,0.2)rectangle(1.5,-0.2)
                      node[black,pos=0.5,scale=0.7]{$(c,d)$}
                      (0,-2.2)rectangle(1.5,-2.6)
                      node[black,pos=0.5,scale=0.7]{$(a,b)$};
\draw (0.75,0.2)--(0.75,0.9) node[midway,left,scale=0.7] {$(c,d)$}
      (0.75,-0.2)--(0.75,-0.9) node[midway,left,scale=0.7] {$(c,d)$}
      (0.75,-1.5)--(0.75,-2.2) node[midway,left,scale=0.7] {$(a,b)$}
       (0.75,1.5)--(0.75,2.2) node[midway,left,scale=0.7] {$(m,n)$}
       (0.75,-2.6)--(0.75,-3.1)node[midway,left,scale=0.7] {$(a,b)$};
\end{tric}
}

\def\samedefd{
\begin{tric}
\draw [darkred,thick] (0.75,1.2) ellipse (1.4 and 0.5)
                      node[black] {${\  }^1 F_{(c,d)(e,f)}^{j{(e,f)}}$}
                      (0.75,-1.2) ellipse (1.4 and 0.5)
                      node[black] {${\  }^2 F_{(c,d)(e,f)}^{j{(e,f)}}$}
                     (0,0.2)rectangle(1.5,-0.2)node[pos=0.5,black,scale=0.7] {$(e,f)$}
                      (0.75,2.5) ellipse (0.75 and 0.3)
                      node[black,scale=0.7] {${\  }^1 D$}
                      (0.75,-2.5) ellipse (0.75 and 0.3)
                      node[black,scale=0.7] {${\  }^2 D$}
                      (0,-3.3)rectangle(1.5,-3.7)node[pos=0.5,black,scale=0.7] {$(a,b)$} ;
\draw (0.75,0.2)--(0.75,0.7) node[midway,left,scale=0.7] {$(e,f)$}
       (0.75,-0.2)--(0.75,-0.7) node[midway,left,scale=0.7] {$(e,f)$}
      (0.75,-1.7)--(0.75,-2.2) node[midway,left,scale=0.7] {$(c,d)$}
       (0.75,1.7)--(0.75,2.2) node[midway,left,scale=0.7] {$(c,d)$}
       (0.75,-3.3)--(0.75,-2.8) node[midway,left,scale=0.7] {$(a,b)$}
       (0.75,3.3)--(0.75,2.8) node[midway,left,scale=0.7] {$(m,n)$}
       (0.75,-3.7)--(0.75,-4.2) node[midway,left,scale=0.7] {$(a,b)$};
\end{tric}
}

\begin{definition}\label{gclaspdef}
Let $\unw \in \DD$. A diagrammatic $\unw$-$\textbf{clasp}$ is a morphism $C_{\unw}\in \End_{\DD}(\unw)$ which satisfies the following conditions:
\begin{enumerate} 
\item $C_{\unw} \ne 0$,
\item $C_{\unw}\circ C_{\unw}= C_{\unw}$,
\item If $D\in \Hom_{\DD}(\unw, \unu)$ and $\wt \unu< \wt \unw$, then $D\circ C_{\unw} = 0$.
\end{enumerate}
\end{definition}

\begin{remark}
Note that we only use the terminology clasp to refer to idempotents. This is consistent with Kuperberg's original use of the term \cite{Kupe}, but less general than Elias's \cite[Definition 1.12]{elias2015light}. In Section \ref{graphicalclasp}, we will define \emph{generalized clasps}, which will agree with Elias's notion of clasp. 
\end{remark}

\begin{lemma}\label{uniqueessofclasps}
If the $\unw$ clasp exists, then it is unique and $N_{C_{\unw}}= 1$. 
\end{lemma}
\begin{proof}
Suppose $C_{\unw}$ and $C_{\unw}'$ are both $\unw$-clasps. By Lemma \ref{factorizationlemma} we can write $C_{\unw} = N_{C_{\unw}}\id + \sum_i A_i\circ B_i$ and $C_{\unw}' = N_{C_{\unw}'}\id + \sum_i A_i'\circ B_i'$. As a consequence of the definition of clasps, we find
\[
C_{\unw} = C_{\unw}\circ C_{\unw} = \left(N_{C_{\unw}}\id + \sum_i A_i\circ B_i\right)\circ C_{\unw} = N_{C_{\unw}}C_{\unw}
\]
and
\[
C_{\unw}' = C_{\unw}'\circ C_{\unw}' = \left(N_{C_{\unw}'}\id + \sum_i A_i'\circ B_i'\right)\circ C_{\unw}' = N_{C_{\unw}'}C_{\unw}'.
\]
Since $\unw$-clasps are non-zero elements of the vector space $\End_{\DD}(\unw)$, it follows that $N_{C_{\unw}}= 1= N_{C_{\unw}'}$. Thus,
\[
C_{\unw}' = N_{C_{\unw}} C_{\unw}'= \left(N_{C_{\unw}}\id + \sum_i A_i\circ B_i\right)C_{\unw}'=C_{\unw}\circ C_{\unw}'=C_{\unw}\left(N_{C_{\unw}'}\id + \sum_i A_i'\circ B_i'\right) = N_{C_{\unw}'}C_{\unw} = C_{\unw}. 
\]
\end{proof}

\begin{defn}
Let $\pi_{\unw}\in \End_{\Uq}(V(\unw))$ be the idempotent  endomorphism with image $V(\wt \unw)$. The endomorphism $\Phi^{-1}(\pi_{\unw})$ in $\End_{\DD}(\unw)$ is the \textbf{algebraic $\unw$-clasp}
\end{defn}
\begin{lemma}\label{algclaspisgraphclasp}
The algebraic clasp is a clasp, and $\Phi^{-1}(\pi_{\unw}) = C_{\unw}$.
\end{lemma}
\begin{proof}
Since the algebraic clasp is non-zero and idempotent, we just need to argue that the algebraic clasp satisfies the third condition in the definition of clasp. 
Fix $\unu$ such that $\wt \unu< \wt \unw$ and let $D\in \Hom_{\DD}(\unw, \unu)$. The module $V(\wt \unw)$ is not isomorphic to any summand of $V(\unu)$, so we know that $\Hom_{\Uq}(V(\wt \unw), V(\unu))= 0$. Therefore, $\Phi(D) \circ \pi_{\unw}= 0$, and by Theorem \ref{equivalencethm} we may conclude that $D\circ \left( \Phi^{-1}(\pi_{\unw}) \right)= 0$. 
\end{proof}

\begin{lemma}
Let $\unw$ be an object in $\DD$, then there is a unique clasp $C_{\unw}\in \End_{\DD}(\unw)$ and $\Phi(C_{\unw})$ is an idempotent endomorphism of $V(\unw)$ projecting to $V(\wt \unw)$. 
\end{lemma}
\begin{proof}
From Lemma \ref{algclaspisgraphclasp} we see that clasps exist and map under $\Phi$ to the projector for $V(\wt \unw)$. Uniqueness follows from Lemma \ref{uniqueessofclasps}.
\end{proof}

\begin{lemma}\label{newgclaspdef}
Let $E\in \End_{\DD}(\unw)$ be a non-zero endomorphism so that $E^2= E$. If $C_{\unu}\circ D\circ E= 0$ for all $D\in \Hom_{\DD}(\unw, \unu)$ such that $\wt \unu< \wt \unw$, then $E= C_{\unw}$. 
\end{lemma}
\begin{proof}
Since we assume $E$ is non-zero and idempotent, we just need to show that $E$ satisfies the third condition in Definition \ref{gclaspdef}. Fix $\unu$ such that $\wt \unu < \wt \unw$. Suppose inductively that $B\circ E= 0$ for all $B\in \Hom_{\DD}(\unw, \unv)$ where $\wt \unv < \wt \unu$. By Lemma \ref{factorizationlemma} we can write $C_{\unu} = \id + \sum_i A_i\circ B_i$ where each $A_i$ has domain $\unv_i$ such that $\wt \unv_i< \wt \unu$. If $D\in \Hom_{\DD}(\unw, \unu)$, then 
\begin{equation}\label{E:DE-equation}
D\circ E = \id \circ D \circ E = C_{\unu}\circ D \circ E - \sum_i A_i\circ B_i \circ D \circ E = - \sum_i A_i\circ B_i\circ D\circ E. 
\end{equation}
Since $B_i\in \Hom_{\DD}(\unu, \unv_i)$ and $D\in \Hom_{\DD}(\unw, \unu)$, we have $B_i\circ D\in \Hom_{\DD}(\unw, \unv_i)$. The induction hypothesis applies, so $(B_i\circ D)\circ E = 0$, and Equation \ref{E:DE-equation} implies $D\circ E =0$. 
\end{proof}

\begin{lemma}[Clasp Schur's Lemma]\label{claspschur}
Let $\unu, \unv\in \DD$ and let $D\in \Hom_{\DD}(\unu, \unv)$.
\begin{itemize}
\item If $\wt \unu \ne \wt \unv$, then $C_{\unv}\circ D\circ C_{\unu} =0$.
\item If $\unu = \unv$, then $C_{\unv}\circ D\circ C_{\unu} = N_D\cdot C_{\unu}$. 
\end{itemize}
\end{lemma}
\begin{proof}
By Corollary \ref{karoubi}, we find:
$\dim\Hom_{\Kar\DD}(C_{\unu}, C_{\unv})= \delta_{\wt\unu, \wt\unv}$. Thus, we can deduce the following.
\begin{itemize}
\item If $\wt\unu \ne \wt \unv$, then $\Hom_{\Kar\DD}(C_{\unu}, C_{\unv}) = 0$.
\item If $\unu= \unv$, then $\Hom_{\Kar\DD}(C_{\unu}, C_{\unv}) = \mathbb{C}(q)\cdot C_{\unu}$.
\end{itemize}
\end{proof}

\begin{lemma}\label{claspedwebspace}
Let $\unw, \unu$, and $\unv\in \DD$. If $V(\wt\unv)$ is not a direct summand of $V(\wt\unw)\otimes V(\wt\unu)$, then 
\[
C_{\unv}\circ D\circ (C_{\unw}\otimes C_{\unu}) = 0
\]
for all $D\in \Hom_{\DD}(\unw\otimes\unu,\unv)$. 
\end{lemma}

\begin{proof}
Corollary \ref{karoubi} implies that $\dim\Hom_{\Kar(\DD)}(C_{\unw}\otimes C_{\unu}, C_{\unv})=0$.
\end{proof}

\begin{lemma}[Clasp absorption]\label{claspabsorption}
Let $\unw= \unx \ot \uny \ot \unz$ in $\DD$, then
\[
(\id_{\unx}\otimes C_{\uny}\otimes\id_{\unz})\circ C_{\unw} = C_{\unw}=  C_{\unw}\circ (\id_{\unx}\otimes C_{\uny}\otimes \id_{\unz}).
\]
\end{lemma}
\begin{proof}
Since $V(\wt\unw)$ appears with multiplicity one in $V(\unw)$, it follows that $\pi_{\unw}$ is a central idempotent in $\End_{\Uq}(V(\unw))$. Therefore, $(\id_{\unx}\otimes C_{\uny}\otimes\id_{\unz})\circ C_{\unw}= C_{\unw}\circ (\id_{\unx}\otimes C_{\uny}\otimes \id_{\unz})$ is also an idempotent and $C_{\unu}\circ D \circ (\id_{\unx}\otimes C_{\uny}\otimes \id_{\unz})\circ C_{\unw} = 0$ for all $D\in \Hom_{\DD}(\unw, \unu)$ such that $\wt \unu< \wt \unw$. Thus, by Lemma \ref{newgclaspdef} it suffices to show that $(\id_{\unx}\otimes C_{\uny}\otimes \id_{\unz})\circ C_{\unw} \ne 0$. This is deduced from observing that the morphism $\Phi((\id_{\unx}\otimes C_{\uny}\otimes \id_{\unz})\circ C_{\unw})$ acts on $V(\unw)_{\wt \unw}$ as multiplication by $1$.  
\end{proof}

\subsection{Neutral diagrams and generalized clasps}
\label{graphicalclasp}

\begin{definition}
We will write $\Xa:= \Ha$ and $\Xb:= \Hb$. These are the \textbf{basic neutral diagrams}. 
\end{definition}

\begin{lemma}[Neutral absorption]\label{neturalabsorption}
Let $\unw= \unw_1 \varpi_1\varpi_2\unw_2$ and $\unw' = \unw_1\varpi_2\varpi_1\unw_2$. Then we have the following equality of morphisms.
\[
(\id_{\unw_1}\ot \Xa\ot \id_{\unw_2})\circ C_{\unw}\circ (\id_{\unw_1}\ot \Xb\ot \id_{\unw_2})= C_{\unw'}
\]
\end{lemma}
\begin{proof}
Write $\mathsf{H}_{\unw}^{\unw'}:= \id_{\unw_1}\ot \Xa\ot \id_{\unw_2}$ and $\mathsf{H}_{\unw'}^{\unw}:= \id_{\unw_1}\ot \Xb\ot \id_{\unw_2}$. By Lemma \ref{uniqueessofclasps} we only need to show that $(\mathsf{H}_{\unw}^{\unw'})\circ C_{\unw}\circ (\mathsf{H}_{\unw'}^{\unw})$ satisfies the defining properties of a clasp.

Let $D\in \Hom_{\DD}(\unw', \unu)$ where $\wt \unu < \wt \unw'$. Then $D\circ \mathsf{H}_{\unw}^{\unw'} \in \Hom_{\DD}(\unw, \unu)$ and $\wt \unw = \wt \unw'> \wt \unu$, so $D\circ \mathsf{H}_{\unw}^{\unw'} \circ C_{\unw} = 0$. So $\mathsf{H}_{\unw}^{\unw'} \circ C_{\unw}\circ \mathsf{H}_{\unw'}^{\unw}$ satisfies the third condition in the definition of clasps.   

The following calculation shows that $\mathsf{H}_{\unw'}^{\unw} \circ \mathsf{H}_{\unw}^{\unw'} \circ C_{\unw} = C_{\unw}$. 
\begin{align} \label{HinverseH}
  \Hmoved &= \frac{1}{[3]} \Hmovee + \frac{[4][6]}{[2]^2[12]}\Hmovef
         +\Hmoveg - \frac{1}{[2]}\Hmoveh \notag\\
           &= - \frac{1}{[2]}\Hmoveh \ \ \  = \ \ \  \Hmovei
\end{align}
So we have 
\[
(\mathsf{H}_{\unw}^{\unw'}\circ C_{\unw}\circ \mathsf{H}_{\unw'}^{\unw})\circ (\mathsf{H}_{\unw}^{\unw'}\circ C_{\unw}\circ \mathsf{H}_{\unw'}^{\unw})=
\mathsf{H}_{\unw}^{\unw'}\circ C_{\unw}\circ(\mathsf{H}_{\unw'}^{\unw} \circ \mathsf{H}_{\unw}^{\unw'} \circ C_{\unw})\circ \mathsf{H}_{\unw'}^{\unw}= \mathsf{H}_{\unw}^{\unw'}\circ C_{\unw}\circ C_{\unw} \circ \mathsf{H}_{\unw'}^{\unw}=\mathsf{H}_{\unw}^{\unw'}\circ C_{\unw} \circ \mathsf{H}_{\unw'}^{\unw}.
\]
This tells us that $\mathsf{H}_{\unw}^{\unw'}\circ C_{\unw} \circ \mathsf{H}_{\unw'}^{\unw}$ satisfies the second condition in the definition of clasps.

What's more, $\mathsf{H}_{\unw}^{\unw'}\circ C_{\unw} \circ \mathsf{H}_{\unw'}^{\unw}\ne 0$. Otherwise  $ C_{\unw} = \mathsf{H}_{\unw'}^{\unw} \circ ( \mathsf{H}_{\unw}^{\unw'}\circ C_{\unw} \circ \mathsf{H}_{\unw'}^{\unw}) \circ \mathsf{H}_{\unw}^{\unw'} = 0 $, which is a contradiction. 
\end{proof}

\begin{definition}
A \textbf{neutral diagram} $\mathsf{N}_{\unw}^{\unw'} \in \Hom_{\DD}(\unw,\unw') $ is a composition of tensor products of identity diagrams and basic neutral diagrams. A \textbf{reduced neutral diagram}, is a neutral diagram such that $\Xb\circ \Xa$ or $\Xa\circ \Xb$ do not occur as subdiagrams of $\mathsf{N}_{\unw}^{\unw'}$.
\end{definition}

\begin{lemma}\label{neutrallemma}
Fix $\unw$ and $\unw'$. 
\begin{enumerate}
\item There is a neutral diagram $\mathsf{N}_{\unw}^{\unw'}$ if and only if $\wt\unw=\wt\unw'$. 
\item If $\wt\unw = \wt \unw'$, then there is a reduced neutral diagram in $\Hom_{\DD}(\unw, \unw')$. 
\item Reduced neutral diagrams are unique.
\item Suppose ${^1\mathsf{N}}_{\unw}^{\unw'}$ and ${^2\mathsf{N}}_{\unw}^{\unw'}$ are two neutral diagrams. Then ${^1\mathsf{N}}_{\unw}^{\unw'}\circ C_{\unw} = {^2\mathsf{N}}_{\unw}^{\unw'}\circ C_{\unw}$.
\end{enumerate}
\end{lemma}
\begin{proof}
Omitted. 
\end{proof}

\begin{notation}
Suppose $\wt\unw= \wt\unw'$, then we will write $\mathsf{H}_{\unw}^{\unw'}$ for the reduced neutral diagram in $\Hom_{\DD}(\unw, \unw')$. 
\end{notation}

\def\neutralconnect{
\begin{tric}
\draw  (0.5,0)--(0.5,0.5)  (1.5,0)--(1.5,0.5) (2,0)--(2,0.5);
\draw[double,thin](0,0)--(0,0.5) (1,0)--(1,0.5);
\draw (0,1)--(0,0.5) (1,1)--(1,0.5) (2,1)--(2,0.5) ;
\draw[double,thin]  (0.5,1)--(0.5,0.5) (1.5,1)--(1.5,0.5) ;
\draw (0,1)--(0,1.5) (0.5,1)--(0.5,1.5)  (1.5,1)--(1.5,1.5)  ;
\draw[double,thin] (1,1)--(1,1.5) (2,1)--(2,1.5);
\draw (0,0.5)--(0.5,0.5) (0.5,1)--(1,1) (1,0.5)--(1.5,0.5) (1.5,1)--(2,1);
\end{tric}
}

\begin{example}
Consider $\unw=\varpi_2\varpi_1\varpi_2\varpi_1\varpi_1$ ,  $\unw'=\varpi_1\varpi_1\varpi_2\varpi_1\varpi_2$. We know that $\wt \unw = \wt \unw' =(3,2) $. The reduced neutral diagram is  
$\mathsf{H}_{\unw}^{\unw'}=\neutralconnect$.
\end{example}

\begin{definition}
Given a diagram $D$ in $\DD$ we will write $\mathbb{D}(D)$ for the diagram obtained by flipping $D$ upside down. Note that \   $\mathbb{D}\left(\mathsf{H}_{\unw}^{\unw'}\right)=\mathsf{H}_{\unw'}^{\unw}$.
\end{definition}

\begin{definition}
Given $\unx, \uny$ so that $\wt\unx= \wt\unw = \wt \uny$, we define the \textbf{generalized clasp} $C_{\unx}^{\uny}:= \mathsf{H}_{\unw}^{\uny}\circ C_{\unw} \circ \mathsf{H}_{\unx}^{\unw}$. From Lemma \ref{neutrallemma} it follows that if $\mathsf{N}_{\unw}^{\uny}$ and $\mathsf{N}_{\unx}^{\unw}$ are any neutral diagrams, then $C_{\unx}^{\uny}= \mathsf{N}_{\unw}^{\uny}\circ C_{\unw}\circ \mathsf{N}_{\unx}^{\unw}$.
\end{definition}

\begin{proposition}
The generalized clasps satisfy the following properties: 
\[
C_{\unx}^{\unx}=C_{\unx}, \quad
C_{\unx}^{\uny}\circ \mathsf{H}_{\unz}^{\unx} =C_{\unz}^{\uny}, \quad \mathsf{H}_{\uny}^{\unz} \circ C_{\unx}^{\uny}  =C_{\unx}^{\unz}, \quad C_{\uny}^{\unz} \circ C_{\unx}^{\uny} = C_{\unx}^{\unz}, \quad \text{and} \quad \mathbb{D}(C_{\unx}^{\uny})= C_{\uny}^{\unx}.
\]
\end{proposition}

\begin{proof}
Exercise for the reader. For hints, see \cite[Proposition $3.2$]{elias2015light}.
\end{proof}


\subsection{Elementary Light Ladders}
\label{subsec-elemLLs}

\begin{notation}
We write $f_1^{(k)}:=\cfrac{f_1^k}{[k]_q!}$ and $f_2^{(k)}:=\cfrac{f_2^k}{[k]_{q^3}!}$, where $[k]_q! := [k][k-1]\dots [2][1]$ and $[k]_{q^3}!:=[k]_{q^3}[k-1]_{q^3}\dots [2]_{q^3}[1]_{q^3}$. Note that $[k]_{q^3} = \frac{[3k]}{[3]}$.
\end{notation}

For each fundamental weight $\varpi \in \lbrace \varpi_1, \varpi_2\rbrace$ we choose a basis $\lbrace ^{i}v_{\mu, \varpi}\rbrace_{ i=1, \ldots, \dim V(\varpi)_{\mu}}$ for all weight spaces $V(\varpi)_{\mu}$. Our convention will be to not record the superscript $i$ in ${^iv}_{\mu, \varpi}$ when the weight space is multiplicity one. Explicitly, we choose the following basis of $V(\varpi_1)$:
\[
v_{(1, 0), \varpi_1}= v_1, \ \ \  v_{(-1, 1), \varpi_1}= f_1v_1, \ \ \  v_{(2, -1), \varpi_1}= f_2f_1v_1, \ \ \  v_{(0, 0), \varpi_1}= f_1f_2f_1v_1,
\]
\[
v_{(-2, 1), \varpi_1}= f_1^{(2)}f_2f_1v_1, \ \ \  v_{(1, -1), \varpi_1}= f_2f_1^{(2)}f_2f_1v_1, \ \ \ \text{and} \ \ \ v_{(-1, 0), \varpi_1}= f_1f_2f_1^{(2)}f_2f_1v_1,
\]
and the following basis of $V(\varpi_2)$:
\[
v_{(0, 1), \varpi_2}= v_2, \ \ \ v_{(3, -1), \varpi_2}=f_2v_2, \ \ \  v_{(1, 0), \varpi_2}= f_1f_2v_2, \ \ \  v_{(-1, 1), \varpi_2}= f_1^{(2)}f_2f_1v_2,
\]
\[
v_{(2, -1), \varpi_2} = f_2f_1^{(2)}f_2v_2, \ \ \  v_{(-3, 2), \varpi_2}= f_1^{(3)}f_2v_2, \ \ \ {^1v}_{(0, 0), \varpi_2}= f_1f_2f_1^{(2)}f_2v_2, \ \ \  {^2v}_{(0, 0), \varpi_2}= f_2f_1^{(3)}f_2v_2,
\]
\[
v_{(3, -2), \varpi_2}= f_2^{(2)}f_1^{(3)}f_2v_2, \ \ \  v_{(-2, 1), \varpi_2}= f_1^{(2)}f_2f_1^{(2)}f_2v_2, \ \ \ v_{(1, -1), \varpi_2}= f_1f_2^{(2)}f_1^{(3)}f_2v_2, \ \ \  v_{(-1, 0), \varpi_2}= f_1^{(2)}f_2^{(2)}f_1^{(3)}f_2v_2,
\]
\[
v_{(-3, 1), \varpi_2}= f_1^{(3)}f_2^{(2)}f_1^{(3)}f_2v_2, \ \ \ \text{and} \ \ \ v_{(0, -1), \varpi_2}= f_2f_1^{(3)}f_2^{(2)}f_1^{(3)}f_2v_2.
\]

\begin{rmk}
The following relation holds in $V(\varpi_2)$:
\[
f_1f_2^{(2)}f_1^{(3)}f_2v_2=f_2f_1^{(2)}f_2f_1^{(2)}f_2v_2.
\]
Thus, there are two ways to present the vector $v_{(1, -1), \varpi_2}$.
\end{rmk}

\begin{defn}
For each vector $^i v_{\mu, \varpi}\in V(\varpi)$, we associate a diagram in $\DD$ denoted ${^iL}_{\mu, \varpi}$. Our convention will be to not record the superscript $i$ in ${^iL}_{\mu, \varpi}$ when the weight space is multiplicity one.

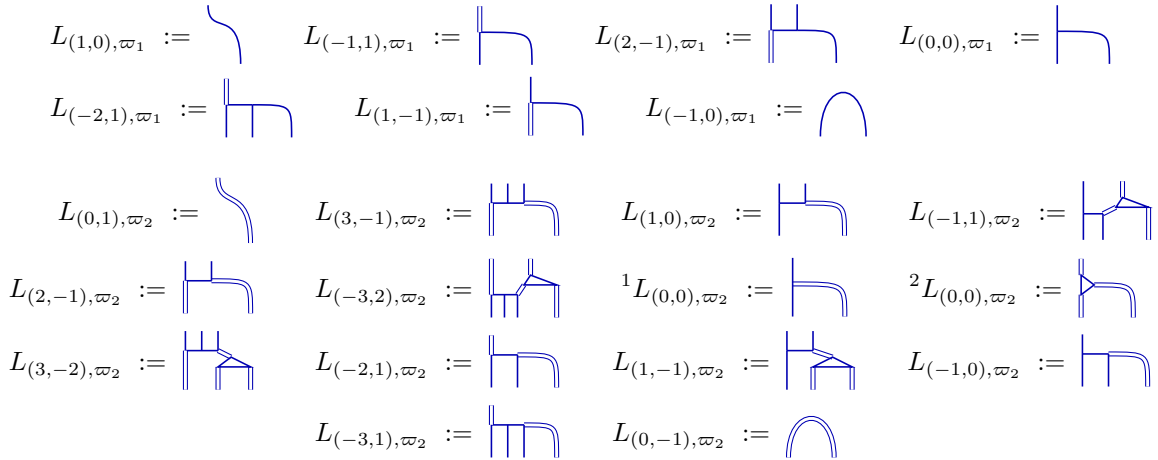
\def\elKb{
\begin{tric}
\draw (5.2,-1) ..controls (5.2,-0.1)and(5.9,-0.1) ..(5.9,-1);
\end{tric}}

\def\elKc{
\begin{tric}
\draw (9,-0.5) ..controls (9.8,-0.5) ..(9.8,-1);
\draw  (9,-0.1)--(9,-0.5);
\draw[double,thin] (9,-1)--(9,-0.5);
\end{tric}}
 
\def\elKd{
\begin{tric}
\draw (12.8,-0.5) --(12.8,-1);
\draw (13.2,-0.5) ..controls (13.8,-0.5) ..(13.8,-1);
\draw (13.2,-0.5) --(13.2,-1);
\draw (13.2,-0.5) -- (12.8,-0.5);
\draw[double,thin] (12.8,-0.5) --(12.8,-0.1);
\end{tric}}
 
\def\elKe{
\begin{tric}
\draw (1,-0.5) ..controls (1.8,-0.5) ..(1.8,-1);
\draw  (1,-0.1)--(1,-0.5);
\draw (1,-1)--(1,-0.5);
\end{tric}}

\def\elKf{
\begin{tric}
\draw (5.2,-0.5) ..controls (5.8,-0.5) ..(5.8,-1);
\draw (5.2,-0.5) --(5.2,-0.1);
\draw (5.2,-0.5) -- (4.8,-0.5);
\draw (4.8,-0.5) --(4.8,-0.1);
\draw[double,thin] (4.8,-0.5) --(4.8,-1);
\end{tric}}

\def\elKg {
\begin{tric}
\draw (9,-0.5) ..controls (9.8,-0.5) ..(9.8,-1);
\draw[double,thin]  (9,-0.1)--(9,-0.5);
\draw (9,-1)--(9,-0.5);
\end{tric}}

\def\elKh {
\begin{tric}
\draw (13.3,-0.1)..controls(13.3,-0.6)and(13.8,-0.1)..(13.8,-1);
\end{tric}}

\def\elRb{
\begin{tric}
\draw[double,thin] (1.2,-1) ..controls (1.2,-0.2)and(1.9,-0.2) ..(1.9,-1);
\end{tric}}

\def\elRc{
\begin{tric}
\draw [double,thin](1.3,-0.5)..controls (1.8,-0.5)..(1.8,-1);
\draw (1.3,-0.5) --(1.3,-1);
\draw (1.3,-0.5) -- (0.8,-0.5);
\draw [double,thin](0.8,-0.5) --(0.8,-0.2);
\draw (0.8,-0.5) --(0.8,-1);
\draw (1.05,-0.5)--(1.05,-1);
\end{tric}}

\def\elRd{
\begin{tric}
\draw [double,thin](1.2,-0.5)..controls (1.8,-0.5)..(1.8,-1);
\draw (1.2,-0.5) --(1.2,-1);
\draw (1.2,-0.5) -- (0.8,-0.5);
\draw (0.8,-0.5) --(0.8,-0.2);
\draw (0.8,-0.5) --(0.8,-1);
\end{tric}}

\def\elRe{
\begin{tric}
\draw [double,thin](1.2,-0.4)--(1.5,-0.5);
\draw [double,thin](1.2,-1)--(1.2,-0.65);
\draw [double,thin](1.8,-1)--(1.8,-0.65);
\draw (1.5,-0.5)--(1.2,-0.65)--(1.8,-0.65)--cycle;
\draw (1.2,-0.4) --(1.2,-0.1);
\draw (1.2,-0.4) -- (0.8,-0.4);
\draw (0.8,-0.4) --(0.8,-0.1);
\draw (0.8,-0.4) --(0.8,-1);
\end{tric}
}

\def\elRf{
\begin{tric}
\draw [double,thin](1.2,-0.5)..controls (1.8,-0.5)..(1.8,-1);
\draw (1.2,-0.5) --(1.2,-1);
\draw (1.2,-0.5) -- (0.8,-0.5);
\draw [double,thin](0.8,-0.5) --(0.8,-0.2);
\draw (0.8,-0.5) --(0.8,-1);
\end{tric}
}

\def\elRg{
\begin{tric}
\draw [double,thin](1.3,-0.4)--(1.5,-0.5);
\draw [double,thin](1.3,-1)--(1.3,-0.65);
\draw [double,thin](1.8,-1)--(1.8,-0.65);
\draw (1.5,-0.5)--(1.3,-0.65)--(1.8,-0.65)--cycle;
\draw (1.3,-0.4) --(1.3,-0.1);
\draw (1.3,-0.4) -- (0.8,-0.4);
\draw (0.8,-0.4) --(0.8,-0.1);
\draw [double,thin](0.8,-0.4) --(0.8,-1);
\draw (1.05,-0.4)--(1.05,-0.1);
\end{tric}
}

\def\elRh{
\begin{tric}
\draw[double,thin] (1,-0.5) ..controls (1.8,-0.5) ..(1.8,-1);
\draw  (1,-0.1)--(1,-0.5);
\draw (1,-1)--(1,-0.5);
\end{tric}
}

\def\elRi{
\begin{tric}
\draw[double,thin] (1.2,-0.5) ..controls (1.8,-0.5) ..(1.8,-1);
\draw [double,thin](1,-0.1)--(1,-0.35);
\draw [double,thin](1,-1)--(1,-0.65);
\draw (1,-0.35)--(1,-0.65)--(1.2,-0.5)--cycle;
\end{tric}
}

\def\elRl{
\begin{tric}
\draw [double,thin](1.2,-0.65)--(1.3,-0.5);
\draw [double,thin](1.4,-0.1)--(1.4,-0.35);
\draw [double,thin](1.8,-1)--(1.8,-0.5);
\draw (1.3,-0.5)--(1.4,-0.35)--(1.8,-0.5)--cycle;
\draw (1.2,-0.65) --(1.2,-1);
\draw (1.2,-0.65) -- (0.8,-0.65);
\draw (0.8,-0.65) --(0.8,-1);
\draw [double,thin](0.8,-0.65) --(0.8,-0.1);
\draw (1,-0.65)--(1,-1);
\end{tric}
}

\def\elRm{
\begin{tric}
\draw [double,thin](1.2,-0.5)..controls (1.8,-0.5)..(1.8,-1);
\draw (1.2,-0.5) --(1.2,-0.2);
\draw (1.2,-0.5) -- (0.8,-0.5);
\draw (0.8,-0.5) --(0.8,-0.2);
\draw [double,thin](0.8,-0.5) --(0.8,-1);
\end{tric}
}

\def\elRn{    
\begin{tric}
\draw [double,thin](1.1,-0.6)--(1.3,-0.5);
\draw [double,thin](1.4,-0.1)--(1.4,-0.35);
\draw [double,thin](1.8,-1)--(1.8,-0.5);
\draw (1.3,-0.5)--(1.4,-0.35)--(1.8,-0.5)--cycle;
\draw (1.1,-0.6) --(1.1,-1);
\draw (1.1,-0.6) -- (0.8,-0.6);
\draw (0.8,-0.6) --(0.8,-1);
\draw (0.8,-0.6) --(0.8,-0.1);
\end{tric}
}

\def\elRo{    
\begin{tric}
\draw [double,thin](1.2,-0.5)..controls (1.8,-0.5)..(1.8,-1);
\draw (1.2,-0.5) --(1.2,-0.2);
\draw (1.2,-0.5) -- (0.8,-0.5);
\draw (0.8,-0.5) --(0.8,-0.2);
\draw (0.8,-0.5) --(0.8,-1);
\end{tric}  
}

\def\elRp{    
\begin{tric}
\draw [double,thin](1.3,-0.5)..controls (1.8,-0.5)..(1.8,-1);
\draw (1.3,-0.5) --(1.3,-0.2);
\draw (1.3,-0.5) -- (0.8,-0.5);
\draw (0.8,-0.5) --(0.8,-0.2);
\draw (1.05,-0.5)--(1.05,-0.2);
\draw [double,thin](0.8,-0.5) --(0.8,-1);
\end{tric}
}

\def\elRq{    
\begin{tric}
\draw [double,thin](13.3,0)..controls(13.3,-0.5)and(13.8,-0.1)..(13.8,-1);
\end{tric}
}   
 
\begin{align*}
  &L_{(1,0), \varpi_1}  \ := \ \elKh \ \ \ \ \ \ \ 
   L_{(-1,1), \varpi_1}    \ := \  \elKg \ \ \ \ \ \ \ 
   L_{(2,-1), \varpi_1}    \ := \  \elKf \ \ \ \ \ \ \ 
   L_{(0,0), \varpi_1}    \ := \  \elKe\\
  &L_{(-2,1), \varpi_1}    \ := \  \elKd \ \ \ \ \ \ \ 
   L_{(1,-1), \varpi_1}    \ := \  \elKc \ \ \ \ \ \ \ 
   L_{(-1,0), \varpi_1}    \ := \  \elKb
\end{align*}
\begin{align*}
   L_{(0,1), \varpi_2}    \ := \  \elRq \ \ \ \ \ \ \
    &L_{(3,-1), \varpi_2}    \ := \  \elRp \ \ \ \ \ \ \
    L_{(1,0), \varpi_2}    \ := \  \elRo \ \ \ \ \ \ \
    L_{(-1,1), \varpi_2}    \ := \  \elRn \\
    L_{(2,-1), \varpi_2}    \ := \  \elRm \ \ \ \ \ \ \
   &L_{(-3,2), \varpi_2}    \ := \  \elRl \ \ \ \ \ \ \
    {^1L}_{(0,0),\varpi_2}    \ := \  \elRh \ \ \ \ \ \ \
    {^2L}_{(0,0),\varpi_2}    \ := \  \elRi \\
    L_{(3,-2), \varpi_2}    \ := \  \elRg \ \ \ \ \ \ \
    &L_{(-2,1), \varpi_2}    \ := \  \elRf \ \ \ \ \ \ 
    L_{(1,-1), \varpi_2}   \ := \  \elRe \ \ \ \ \ \ 
    L_{(-1,0), \varpi_2}    \ := \  \elRd \\
    &L_{(-3,1), \varpi_2}    \ := \  \elRc \ \ \ \ \ \ 
    L_{(0,-1), \varpi_2}    \ := \  \elRb 
\end{align*}

The diagram ${^iL}_{\mu, \varpi}$ is a morphism from $^i\unx_{\mu, \varpi}\otimes \varpi \rightarrow ^i\uny_{\mu, \varpi}$. We will refer to $^i\unx_{\mu, \varpi}$ as the \textbf{in strand} of ${^iL}_{\mu, \varpi}$ and $^i\uny_{\mu, \varpi}$ as the \textbf{out strand} of ${^iL}_{\mu, \varpi}$. Note that $\mu= \wt(^i\uny_{\mu, \varpi})-\wt(^i\unx_{\mu, \varpi})$. 
\end{defn}

\def\elRf{
\begin{tric}
\draw [double,thin](1.2,-0.5)..controls (1.8,-0.5)..(1.8,-1);
\draw (1.2,-0.5) --(1.2,-1);
\draw (1.2,-0.5) -- (0.8,-0.5);
\draw [double,thin](0.8,-0.5) --(0.8,-0.2);
\draw (0.8,-0.5) --(0.8,-1);
\end{tric}
}

\begin{example}
   Consider  $$ L_{(-2,1), \varpi_2}    \ := \  \elRf \ \  { \ }_,  $$  we know that \ $\unx_{(-2,1), \varpi_2}=\varpi_1\varpi_1$ , $\uny_{(-2,1), \varpi_2}=\varpi_2$, and  $(-2,1)= \wt(\uny_{(-2,1), \varpi_2})-\wt(\unx_{(-2,1), \varpi_2})$. 
\end{example}

\begin{notation}
If $W$ is a subspace of a $\Uq$ module, we will write 
\[
\Kab(W) := \lbrace w\in W \ : \ e_1^{a+1}w= 0 = e_2^{b+1}w\rbrace.
\]
\end{notation}

\begin{lemma}
Let $a,b\in \mathbb{Z}_{\ge 0}$. Fix a fundamental weight $\varpi$ and let $\mu\in \wt V(\varpi)$. Then
\begin{equation}\label{PRV}
[V(a, b)\otimes V(\varpi): V((a, b)+ \mu)] = \dim \Kab(V(\varpi)_{\mu}).
\end{equation} 

\end{lemma}
\begin{proof}
Follows from \cite[Theorem 2.1]{PRVformula}.
\end{proof} 

\begin{lemma}\label{vecleaf}
The following are equivalent: 
\begin{enumerate}
    \item $^i v_{\mu, \varpi}\in \Kab(V(\varpi)_{\mu})$, \quad \text{and}
    \item There is $(c, d)\in \mathbb{N}\times \mathbb{N}$ such that $\wt (^i\unx_{\mu, \varpi}) + (c, d) =  (a, b)$.
\end{enumerate}
\end{lemma}

\begin{proof}
The lemma can be deduced from the following claim: the weight of the in strand for ${^iL}_{\mu, \varpi}$, $\wt ( ^i\unx_{\mu, \varpi} )$, is equal to the minimal $(a, b)$ so that $^i v_{\mu, \varpi}\in \Kab(V(\varpi)_{\mu})$. The claim is verified from the vector to diagram correspondence $^i v_{\mu, \varpi} \mapsto {^iL}_{\mu, \varpi}$, along with Equation \eqref{PRV}, and the description of action of $e_1$ and $e_2$ on the vectors in each fundamental representation. Computing $e_k\cdot{^iv}_{\mu, \varpi}$ is left as an exercise, the most interesting case is the zero weight space for the second fundamental representation.
\end{proof}

\begin{example}
When $a\geq 3$ and $b\geq 2$ ,  $[V(a, b)\otimes V(\varpi_2): V((a, b)+ \mu)]= 1$ when $\mu \neq (0,0)$, and $[V(a, b)\otimes V(\varpi_2): V(a, b)]=2$. The reader should compare this with the observation that for each ${^iL}_{\mu, \varpi_2}$ the number of $\varpi_1$ colored in strands is less than or equal to $3$ and the number of $\varpi_2$ colored in strands is less than or equal to $2$.
\end{example}

For each dominant integral weight $\lambda = a\varpi_1 + b\varpi_2\in X_+$, we choose a distinguished object $\unu_{\lambda}\in \DD$ such that $\wt \unu_{\lambda} = \lambda$. 

\begin{example}
We must have $\unu_{(2, 0)}= \varpi_1\varpi_1$ and for $\unu_{(1, 1)}$ we choose one of $ \varpi_1\varpi_2$ or $\varpi_2\varpi_1$.
\end{example}

\begin{defn}
Let $\unw$ be an object in $\DD$ and let $\lambda= \wt \unw$. Suppose that $^i v_{\mu, \varpi} \in \text{Ker}_{\lambda}(V(\varpi)_{\mu})$, so in particular $\lambda+ \mu\in X_+$. Let $(c,d):=\lambda - \wt( ^i \unx_{\mu,\varpi})$, then by Lemma \ref{vecleaf} we have $(c,d) \in \mathbb{N}\times \mathbb{N}$, so there is a reduced neutral diagram 
\[
\mathsf{H}_{\unw}^{\underline{s}_{(c,d)}\ot (^i\unx_{\mu, \varpi})}: \unw \rightarrow {\underline{s}_{(c,d)}}\ot (^i\unx_{\mu, \varpi}).
\] 
Since $(c,d)+\wt(^i\uny_{\mu, \varpi})=\lambda+\mu$, there is also a reduced neutral diagram
\[
\mathsf{H}_{\underline{s}_{(c,d)} \ot (^i\uny_{\mu, \varpi})}^{\unu_{\lambda+\mu}}:
\underline{s}_{(c,d)}\ot  (^i\uny_{\mu, \varpi}) \rightarrow \unu_{\lambda+ \mu}.
\]
We define the \textbf{elementary light ladder diagram} to be
\[
^iELL_{\unw, \varpi}^{\unu_{\lambda+ \mu}}:= 
\mathsf{H}_{\underline{s}_{(c,d)} \ot (^i\uny_{\mu, \varpi})}^{\unu_{\lambda+\mu}}\circ
\Big( (\id_{\underline{s}_{(c,d)}}) \ot {^iL}_{\mu, \varpi} \Big)\circ \Big( \mathsf{H}_{\unw}^{\underline{s}_{(c,d)}\ot (^i\unx_{\mu, \varpi})}\ot \id_{\varpi} \Big).
\]
\end{defn}

\def\elexa{
\begin{tric}
\draw [double,thin](1.2,-0.5)..controls (1.8,-0.5)..(1.8,-1);
\draw (1.2,-0.5) --(1.2,-1);
\draw (1.2,-0.5) -- (0.8,-0.5);
\draw [double,thin](0.8,-0.5) --(0.8,0);
\draw (0.8,-0.5) --(0.8,-1);
\draw [double,thin] (0.5,-1)--(0.5,0);
\draw (0.2,-1)--(0.2,0);
\end{tric}
}

\def\elexb{
\begin{tric}
\draw[double,thin] (1,0) ..controls (1.5,0) ..(1.5,-1);
\draw  (0.7,1)--(0.7,0.4)--(1,0.4)--(1,-1);
\draw (0.7,-1)--(0.7,-0.4)--(0.4,-0.4)--(0.4,1);
\draw[double,thin] (0.7,-0.4)--(0.7,0.4) (1,0.4)--(1,1) (0.4,-0.4)--(0.4,-1);
\draw (0.1,-1)--(0.1,1);
\end{tric}
}

\def\elexc{
\begin{tric}
\draw[double,thin] (1.2,-0.5) ..controls (1.5,-0.5) ..(1.5,-1.5);
\draw [double,thin](1,0.5)--(1,-0.35);
\draw [double,thin](1,-0.9)--(1,-0.65);
\draw (1,-0.35)--(1,-0.65)--(1.2,-0.5)--cycle;
\draw (1,-1.5)--(1,-0.9)--(0.7,-0.9)--(0.7,0.5);
\draw (0.7,-1.5)--(0.7,-1.2)--(0.4,-1.2)--(0.4,0.5);
\draw [double,thin] (0.4,-1.2)--(0.4,-1.5) (0.7,-1.2)--(0.7,-0.9);
\draw (0.1,-1.5)--(0.1,0.5);
\end{tric}
}

\begin{example}
Consider $\unw=\varpi_1\varpi_2\varpi_1\varpi_1$ and  \ $\varpi=\varpi_2$. We know that $\lambda= \wt \unw = (3,1)$. When $\mu= (-2,1) $, so $\lambda+\mu=(1,2)$, choose $\unu_{\lambda+\mu}=\varpi_1\varpi_2\varpi_2$. Then 
\[
^1ELL_{\unw, \varpi}^{\unu_{\lambda+ \mu}}=\elexa \quad .
\]
When $\mu= (0,0)$, $\lambda+\mu=(3,1)$, choose $\unu_{\lambda+\mu}=\varpi_1\varpi_1\varpi_1\varpi_2$. Then
\[
^1ELL_{\unw, \varpi}^{\unu_{\lambda+ \mu}}=\elexb \qquad \text{and} \qquad ^2ELL_{\unw, \varpi}^{\unu_{\lambda+ \mu}}=\elexc \quad .
\]
\end{example}

\begin{defn}\label{claspedbasis}
Let $\unw \in \DD$. Write $\lambda= \wt\unw$, and suppose that $^i v_{\mu, \varpi}\in \text{Ker}_{\lambda}(V(\varpi)_{\mu})$. We define the \textbf{(clasped) light ladder diagram} to be the following diagram:
\[
^iLL_{\unw, \varpi}^{\unu_{\lambda+ \mu}}:= C_{\unu_{\lambda+ \mu}}\circ (^iELL_{\unw, \varpi}^{\unu_{\lambda+ \mu}})\circ (C_{\unw}\otimes \id_{\varpi}).
\]
\end{defn}

\begin{lemma}
Let $\mathsf{N}_{\unw}^{\unw'}:\unw\rightarrow \unw'$ be a neutral diagram. Then
\[
 {^{i}LL}_{\unw', \varpi}^{\unu_{\lambda+ \mu}}
\circ (N_{\unw}^{\unw'}\ot \id_{\varpi})\circ C_{\unw}\ot \id_{\varpi}= {^{i}LL}_{\unw, \varpi}^{\unu_{\lambda+ \mu}}.
\]
\end{lemma}

\begin{proof}
Follows from Lemma \ref{neutrallemma}.
\end{proof}

\begin{defn}
Suppose that $^i v_{\mu, \varpi}, ^j v_{\mu, \varpi}\in \text{Ker}_{\wt\unw}(V(\varpi)_{\mu})$ (we allow for $i= j$). We define the \textbf{(clasped) double ladder diagram} to be the following diagram in $\End_{\DD}(\unw\otimes \varpi)$:
\[
^{ij}\mathbb{LL}_{\unw, \varpi}^{\unu_{\lambda+ \mu}} := (\mathbb{D}(^iLL_{\unw, \varpi}^{\unu_{\lambda+ \mu}})) \circ (^jLL_{\unw, \varpi}^{\unu_{\lambda+ \mu}}).
\]

In the case that $V(\varpi)_{\mu}$ is one dimensional, we will drop the superscripts $ij$ in double ladders and drop superscript $i$ in (clasped) light ladders.
\end{defn}

\def\dlexa{
\begin{tric}
\draw [double,thin](1.2,-0.8)..controls (1.8,-0.8)..(1.8,-2.2);
\draw (1.2,-0.8) --(1.2,-1.3);
\draw (1.2,-0.8) -- (0.8,-0.8);
\draw [double,thin](0.8,-0.8) --(0.8,-0.3);
\draw (0.8,-0.8) --(0.8,-1.3);
\draw [double,thin] (0.5,-1.3)--(0.5,-0.3);
\draw (0.2,-1.3)--(0.2,-0.3);
\draw (0.2,-1.9)--(0.2,-2.2) (0.8,-1.9)--(0.8,-2.2) (1.2,-1.9)--(1.2,-2.2);
\draw [double,thin] (0.5,-1.9)--(0.5,-2.2);
\draw (0.2,0.3)--(0.2,0.6);
\draw[double,thin](0.5,0.3)--(0.5,0.6) (0.8,0.3)--(0.8,0.6);
\draw [darkred,thick] (-0.1,-1.3) rectangle (1.5,-1.9) node[midway,black]{$C_{\unw}$};
\draw [darkred,thick] (-0.1,-0.3) rectangle (1.5,0.3) node[midway,black]{$C_{\unu_{\lambda+\mu}}$};
\end{tric}
}

\def\dlexb{
\begin{tric}
\draw [double,thin](1.2,-0.8)..controls (1.8,-0.8)..(1.8,-2.2);
\draw (1.2,-0.8) --(1.2,-1.3);
\draw (1.2,-0.8) -- (0.8,-0.8);
\draw [double,thin](0.8,-0.8) --(0.8,-0.3);
\draw (0.8,-0.8) --(0.8,-1.3);
\draw [double,thin] (0.5,-1.3)--(0.5,-0.3);
\draw (0.2,-1.3)--(0.2,-0.3);
\draw (0.2,-1.9)--(0.2,-2.2) (0.8,-1.9)--(0.8,-2.2) (1.2,-1.9)--(1.2,-2.2);
\draw [double,thin] (0.5,-1.9)--(0.5,-2.2);
\draw [darkred,thick] (-0.1,-1.3) rectangle (1.5,-1.9) node[midway,black]{$C_{\unw}$};

\draw [double,thin](1.2,0.8)..controls (1.8,0.8)..(1.8,2.2);
\draw (1.2,0.8) --(1.2,1.3);
\draw (1.2,0.8) -- (0.8,0.8);
\draw [double,thin](0.8,0.8) --(0.8,0.3);
\draw (0.8,0.8) --(0.8,1.3);
\draw [double,thin] (0.5,1.3)--(0.5,0.3);
\draw (0.2,1.3)--(0.2,0.3);
\draw (0.2,1.9)--(0.2,2.2) (0.8,1.9)--(0.8,2.2) (1.2,1.9)--(1.2,2.2);
\draw [double,thin] (0.5,1.9)--(0.5,2.2);
\draw [darkred,thick] (-0.1,1.3) rectangle (1.5,1.9) node[midway,black]{$C_{\unw}$};
\draw [darkred,thick] (-0.1,-0.3) rectangle (1.5,0.3) node[midway,black]{$C_{\unu_{\lambda+\mu}}$};
\end{tric}
}

\begin{example}
Consider $\unw=\varpi_1\varpi_2\varpi_1\varpi_1  ,\ \varpi=\varpi_2.$ When $\mu= (-2,1) $ , choose $\unu_{\lambda+\mu}=\varpi_1\varpi_2\varpi_2 $, then
\[
LL_{\unw, \varpi}^{\unu_{\lambda+ \mu}}=\dlexa
\quad \ \ \text{and} \ \ \quad \mathbb{LL}_{\unw, \varpi}^{\unu_{\lambda+ \mu}}=\dlexb \quad .
\]
\end{example}

\begin{remark}
Using the definition of the elementary light ladder, and basic properties of clasps, we can expand the clasped light ladder: 
\begin{align*}
^iLL_{\unw, \varpi}^{\unu_{\lambda+ \mu}}:&= C_{\unu_{\lambda+ \mu}}\circ (^iELL_{\unw, \varpi}^{\unu_{\lambda+ \mu}})\circ (C_{\unw}\otimes \id_{\varpi}) \\
&=  C_ {\underline{s}_{(c,d)} \ot (^i\uny_{\mu, \varpi})}
^{\unu_{\lambda+ \mu}} 
\circ \Big( (\id_{\underline{s}_{(c,d)}}) \ot {^iL}_{\mu, \varpi}\Big) 
\circ \Big( C_{\unw}
^{\underline{s}_{(c,d)}\ot (^i\unx_{\mu, \varpi})}
\ot \id_{\varpi} \Big).
\end{align*}
We can similarly expand the clasped double ladder:
\begin{align*}
^{ij}\mathbb{LL}_{\unw, \varpi}^{\unu_{\lambda+ \mu}} :&= (\mathbb{D}(^iLL_{\unw, \varpi}^{\unu_{\lambda+ \mu}})) \circ (^jLL_{\unw, \varpi}^{\unu_{\lambda+ \mu}}) \ \  =
( C_{\unw}\ot \id_{\varpi} ) \circ (\mathbb{D}(^iELL_{\unw, \varpi}^{\unu_{\lambda+ \mu}}))\circ C_{\unu_{\lambda+ \mu}} \circ (^jELL_{\unw, \varpi}^{\unu_{\lambda+ \mu}})\circ ( C_{\unw}\ot \id_{\varpi} ) \\
&=\Big( C^{\unw}
_{\underline{s}_{(c,d)}\ot (^i\unx_{\mu, \varpi})}
\ot \id_{\varpi} \Big)
\circ \Big( (\id_{\underline{s}_{(c,d)}}) \ot \mathbb{D}({^iL}_{\mu, \varpi})\Big) 
\circ C_ {\underline{s}_{(e,f)}\ot (^j\uny_{\mu, \varpi})}
^{\underline{s}_{(c,d)}\ot (^i\uny_{\mu, \varpi})} \\
&\ \ \ \ \ \ \ \ \ \ \ \ \ \ \ \ \ \ \ \ \ \ \ \ \ \ \ \ \ \ \ \ \ \ \ \ \ \ \ \ \ \ \ \   \circ \Big( (\id_{\underline{s}_{(e,f)}}) \ot {^jL}_{\mu, \varpi}\Big) 
\circ \Big( C_{\unw}
^{\underline{s}_{(e,f)}\ot (^j\unx_{\mu, \varpi})}
\ot \id_{\varpi} \Big).
\end{align*}
This more complicated looking expanded formula, is actually simpler when viewed in terms of the graphical calculus, as we illustrate in Example \ref{exampleladders}.
\end{remark}

\def\dlexa{
\begin{tric}
\draw [double,thin](1.2,-0.8)..controls (1.8,-0.8)..(1.8,-2.2);
\draw (1.2,-0.8) --(1.2,-1.3);
\draw (1.2,-0.8) -- (0.8,-0.8);
\draw [double,thin](0.8,-0.8) --(0.8,-0.3);
\draw (0.8,-0.8) --(0.8,-1.3);
\draw [double,thin] (0.5,-1.3)--(0.5,-0.3);
\draw (0.2,-1.3)--(0.2,-0.3);
\draw (0.2,-1.9)--(0.2,-2.2) (0.8,-1.9)--(0.8,-2.2) (1.2,-1.9)--(1.2,-2.2);
\draw [double,thin] (0.5,-1.9)--(0.5,-2.2);
\draw (0.2,0.3)--(0.2,0.6);
\draw[double,thin](0.5,0.3)--(0.5,0.6) (0.8,0.3)--(0.8,0.6);
\draw [darkred,thick] (-0.1,-1.3) rectangle (1.5,-1.9) node[midway,black]{$C_{\unw}$};
\draw [darkred,thick] (-0.1,-0.3) rectangle (1.5,0.3) node[midway,black]{$C_{\unu_{\lambda+\mu}}$};
\end{tric}
}

\def\dlexaa{
\begin{tric}
\draw [double,thin](1.2,-0.8)..controls (1.8,-0.8)..(1.8,-2.2);
\draw (1.2,-0.8) --(1.2,-1.3);
\draw (1.2,-0.8) -- (0.8,-0.8);
\draw [double,thin](0.8,-0.8) --(0.8,-0.3);
\draw (0.8,-0.8) --(0.8,-1.3);
\draw (0.4,-1.3)--(0.4,-0.3) node[left, midway, scale=0.7]{$(1,1)$};
\draw (0.2,-1.9)--(0.2,-2.2) (0.8,-1.9)--(0.8,-2.2) (1.2,-1.9)--(1.2,-2.2);
\draw [double,thin] (0.5,-1.9)--(0.5,-2.2);
\draw (0.2,0.3)--(0.2,0.6);
\draw[double,thin](0.5,0.3)--(0.5,0.6) (0.8,0.3)--(0.8,0.6);
\draw [darkred,thick] (-0.1,-1.3) rectangle (1.5,-1.9) node[midway,black]{$(3,1)$};
\draw [darkred,thick] (-0.1,-0.3) rectangle (1.5,0.3) node[midway,black]{$(1,2)$};
\end{tric}
}

\def\dlexc{
\begin{tric}
\draw [double,thin](1.2,-0.8)..controls (1.8,-0.8)..(1.8,-2.2);
\draw (1.2,-0.8) --(1.2,-1.3);
\draw (1.2,-0.8) -- (0.8,-0.8);
\draw [double,thin](0.8,-0.8) --(0.8,-0.3);
\draw (0.8,-0.8) --(0.8,-1.3);
\draw (0.4,-1.3)--(0.4,-0.3)node[left,midway,scale=0.7]{$(1,1)$};
\draw (0.2,-1.9)--(0.2,-2.2) (0.8,-1.9)--(0.8,-2.2) (1.2,-1.9)--(1.2,-2.2);
\draw [double,thin] (0.5,-1.9)--(0.5,-2.2);
\draw [darkred,thick] (-0.1,-1.3) rectangle (1.5,-1.9) node[midway,black]{$(3,1)$};

\draw [double,thin](1.2,0.8)..controls (1.8,0.8)..(1.8,2.2);
\draw (1.2,0.8) --(1.2,1.3);
\draw (1.2,0.8) -- (0.8,0.8);
\draw [double,thin](0.8,0.8) --(0.8,0.3);
\draw (0.8,0.8) --(0.8,1.3);
\draw (0.4,1.3)--(0.4,0.3)node[left,midway,scale=0.7]{$(1,1)$};
\draw (0.2,1.9)--(0.2,2.2) (0.8,1.9)--(0.8,2.2) (1.2,1.9)--(1.2,2.2);
\draw [double,thin] (0.5,1.9)--(0.5,2.2);
\draw [darkred,thick] (-0.1,1.3) rectangle (1.5,1.9) node[midway,black]{$(3,1)$};
\draw [darkred,thick] (-0.1,-0.3) rectangle (1.5,0.3) node[midway,black]{$(1,2)$};
\end{tric}
}

\def\dlexd{
\begin{tric}
\draw[double,thin] (1,-1.3) ..controls (1.8,-1.3) ..(1.8,-3.2);
\draw  (0.7,-0.3)--(0.7,-0.9)--(1,-0.9)--(1,-2.3);
\draw (0.7,-2.3)--(0.7,-1.7)--(0.4,-1.7)--(0.4,-0.3);
\draw[double,thin] (0.7,-1.7)--(0.7,-0.9) (1,-0.9)--(1,-0.3) (0.4,-1.7)--(0.4,-2.3);
\draw (0.1,-2.3)--(0.1,-0.3);
\draw (0.1,-2.9)--(0.1,-3.2) (0.7,-2.9)--(0.7,-3.2) (1,-2.9)--(1,-3.2);
\draw[double,thin] (0.4,-2.9)--(0.4,-3.2);

\draw (0.1,0.3)--(0.1,0.6) (0.4,0.3)--(0.4,0.6) (0.7,0.3)--(0.7,0.6) ;
\draw[double,thin] (1,0.3)--(1,0.6);

\draw [darkred,thick] (-0.3,-2.3) rectangle (1.4,-2.9) node[midway,black]{$C_{\unw}$};
\draw [darkred,thick] (-0.3,-0.3) rectangle (1.4,0.3) node[midway,black]{$C_{\unu_{\lambda+\mu}}$};
\end{tric}
}

\def\dlexda{
\begin{tric}
\draw[double,thin] (1,-1.3) ..controls (1.8,-1.3) ..(1.8,-3.2);
\draw  (1,-0.3)--(1,-2.3);
\draw (0.3,-0.3)--(0.3,-2.3) node[left,midway,scale=0.7]{$(2,1)$};

\draw (0.1,-2.9)--(0.1,-3.2) (0.7,-2.9)--(0.7,-3.2) (1,-2.9)--(1,-3.2);
\draw[double,thin] (0.4,-2.9)--(0.4,-3.2);

\draw (0.1,0.3)--(0.1,0.6) (0.4,0.3)--(0.4,0.6) (0.7,0.3)--(0.7,0.6) ;
\draw[double,thin] (1,0.3)--(1,0.6);

\draw [darkred,thick] (-0.3,-2.3) rectangle (1.4,-2.9) node[midway,black]{$(3,1)$};
\draw [darkred,thick] (-0.3,-0.3) rectangle (1.4,0.3) node[midway,black]{$(3,1)$};
\end{tric}
}

\def\dlexe{
\begin{tric}
\draw[double,thin] (1.2,-1.3) ..controls (1.8,-1.3) ..(1.8,-3.2);
\draw [double,thin](1,-0.3)--(1,-1.15);
\draw [double,thin](1,-1.7)--(1,-1.45);
\draw (1,-1.15)--(1,-1.45)--(1.2,-1.3)--cycle;
\draw (1,-2.3)--(1,-1.7)--(0.7,-1.7)--(0.7,-0.3);
\draw (0.7,-2.3)--(0.7,-2)--(0.4,-2)--(0.4,-0.3);
\draw [double,thin] (0.4,-2)--(0.4,-2.3) (0.7,-2)--(0.7,-1.7);
\draw (0.1,-2.3)--(0.1,-0.3);
\draw (0.1,-2.9)--(0.1,-3.2) (0.7,-2.9)--(0.7,-3.2) (1,-2.9)--(1,-3.2);
\draw[double,thin] (0.4,-2.9)--(0.4,-3.2);

\draw (0.1,0.3)--(0.1,0.6) (0.4,0.3)--(0.4,0.6) (0.7,0.3)--(0.7,0.6) ;
\draw[double,thin] (1,0.3)--(1,0.6);

\draw [darkred,thick] (-0.3,-2.3) rectangle (1.4,-2.9) node[midway,black]{$C_{\unw}$};
\draw [darkred,thick] (-0.3,-0.3) rectangle (1.4,0.3) node[midway,black]{$C_{\unu_{\lambda+\mu}}$};
\end{tric}
}

\def\dlexea{
\begin{tric}
\draw[double,thin] (1.2,-1.3) ..controls (1.8,-1.3) ..(1.8,-3.2);
\draw (1,-1.15)--(1,-1.45)--(1.2,-1.3)--cycle;
\draw[double,thin] (1,-1.15)--(1,-0.3) (1,-1.45)--(1,-2.3);
\draw (0.3,-0.3)--(0.3,-2.3) node[left,midway,scale=0.7]{$(3,0)$};

\draw (0.1,-2.9)--(0.1,-3.2) (0.7,-2.9)--(0.7,-3.2) (1,-2.9)--(1,-3.2);
\draw[double,thin] (0.4,-2.9)--(0.4,-3.2);

\draw (0.1,0.3)--(0.1,0.6) (0.4,0.3)--(0.4,0.6) (0.7,0.3)--(0.7,0.6) ;
\draw[double,thin] (1,0.3)--(1,0.6);

\draw [darkred,thick] (-0.3,-2.3) rectangle (1.4,-2.9) node[midway,black]{$(3,1)$};
\draw [darkred,thick] (-0.3,-0.3) rectangle (1.4,0.3) node[midway,black]{$(3,1)$};
\end{tric}
}

\def\dlexf{
\begin{tric}
\draw[double,thin] (1,-1.3) ..controls (1.8,-1.3) ..(1.8,-3.2);
\draw  (0.7,-0.3)--(0.7,-0.9)--(1,-0.9)--(1,-2.3);
\draw (0.7,-2.3)--(0.7,-1.7)--(0.4,-1.7)--(0.4,-0.3);
\draw[double,thin] (0.7,-1.7)--(0.7,-0.9) (1,-0.9)--(1,-0.3) (0.4,-1.7)--(0.4,-2.3);
\draw (0.1,-2.3)--(0.1,-0.3);
\draw (0.1,-2.9)--(0.1,-3.2) (0.7,-2.9)--(0.7,-3.2) (1,-2.9)--(1,-3.2);
\draw[double,thin] (0.4,-2.9)--(0.4,-3.2);

\draw[double,thin] (1,1.3) ..controls (1.8,1.3) ..(1.8,3.2);
\draw  (0.7,0.3)--(0.7,0.9)--(1,0.9)--(1,2.3);
\draw (0.7,2.3)--(0.7,1.7)--(0.4,1.7)--(0.4,0.3);
\draw[double,thin] (0.7,1.7)--(0.7,0.9) (1,0.9)--(1,0.3) (0.4,1.7)--(0.4,2.3);
\draw (0.1,2.3)--(0.1,0.3);
\draw (0.1,2.9)--(0.1,3.2) (0.7,2.9)--(0.7,3.2) (1,2.9)--(1,3.2);
\draw[double,thin] (0.4,2.9)--(0.4,3.2);

\draw [darkred,thick] (-0.3,-2.3) rectangle (1.4,-2.9) node[midway,black]{$C_{\unw}$};
\draw [darkred,thick] (-0.3,-0.3) rectangle (1.4,0.3) node[midway,black]{$C_{\unu_{\lambda+\mu}}$};
\draw [darkred,thick] (-0.3,2.3) rectangle (1.4,2.9) node[midway,black]{$C_{\unw}$};
\end{tric}
}

\def\dlexg{
\begin{tric}
\draw[double,thin] (1,-1.3) ..controls (1.8,-1.3) ..(1.8,-3.2);
\draw (1,-2.3)--(1,-0.3);
\draw (0.4,-2.3)--(0.4,-0.3) node[left,midway,scale=0.7]{$(2,1)$};
\draw (0.1,-2.9)--(0.1,-3.2) (0.7,-2.9)--(0.7,-3.2) (1,-2.9)--(1,-3.2);
\draw[double,thin] (0.4,-2.9)--(0.4,-3.2);

\draw[double,thin] (1,1.3) ..controls (1.8,1.3) ..(1.8,3.2);
\draw (1,2.3)--(1,0.3);
\draw (0.4,2.3)--(0.4,0.3) node[left,midway,scale=0.7]{$(2,1)$};
\draw (0.1,2.9)--(0.1,3.2) (0.7,2.9)--(0.7,3.2) (1,2.9)--(1,3.2);
\draw[double,thin] (0.4,2.9)--(0.4,3.2);

\draw [darkred,thick] (-0.3,-2.3) rectangle (1.4,-2.9) node[midway,black]{$(3,1)$};
\draw [darkred,thick] (-0.3,-0.3) rectangle (1.4,0.3) node[midway,black]{$(3,1)$};
\draw [darkred,thick] (-0.3,2.3) rectangle (1.4,2.9) node[midway,black]{$(3,1)$};
\end{tric}
}

\def\dlexh{
\begin{tric}
\draw[double,thin] (1.2,-1.3) ..controls (1.8,-1.3) ..(1.8,-3.2);
\draw [double,thin](1,-0.3)--(1,-1.15);
\draw [double,thin](1,-1.7)--(1,-1.45);
\draw (1,-1.15)--(1,-1.45)--(1.2,-1.3)--cycle;
\draw (1,-2.3)--(1,-1.7)--(0.7,-1.7)--(0.7,-0.3);
\draw (0.7,-2.3)--(0.7,-2)--(0.4,-2)--(0.4,-0.3);
\draw [double,thin] (0.4,-2)--(0.4,-2.3) (0.7,-2)--(0.7,-1.7);
\draw (0.1,-2.3)--(0.1,-0.3);
\draw (0.1,-2.9)--(0.1,-3.2) (0.7,-2.9)--(0.7,-3.2) (1,-2.9)--(1,-3.2);
\draw[double,thin] (0.4,-2.9)--(0.4,-3.2);

\draw[double,thin] (1,1.3) ..controls (1.8,1.3) ..(1.8,3.2);
\draw  (0.7,0.3)--(0.7,0.9)--(1,0.9)--(1,2.3);
\draw (0.7,2.3)--(0.7,1.7)--(0.4,1.7)--(0.4,0.3);
\draw[double,thin] (0.7,1.7)--(0.7,0.9) (1,0.9)--(1,0.3) (0.4,1.7)--(0.4,2.3);
\draw (0.1,2.3)--(0.1,0.3);
\draw (0.1,2.9)--(0.1,3.2) (0.7,2.9)--(0.7,3.2) (1,2.9)--(1,3.2);
\draw[double,thin] (0.4,2.9)--(0.4,3.2);

\draw [darkred,thick] (-0.3,-2.3) rectangle (1.4,-2.9) node[midway,black]{$C_{\unw}$};
\draw [darkred,thick] (-0.3,-0.3) rectangle (1.4,0.3) node[midway,black]{$C_{\unu_{\lambda+\mu}}$};
\draw [darkred,thick] (-0.3,2.3) rectangle (1.4,2.9) node[midway,black]{$C_{\unw}$};
\end{tric}
}

\def\dlexi{
\begin{tric}
\draw[double,thin] (1.2,-1.3) ..controls (1.8,-1.3) ..(1.8,-3.2);
\draw [double,thin](1,-0.3)--(1,-1.15);
\draw [double,thin](1,-2.3)--(1,-1.45);
\draw (1,-1.15)--(1,-1.45)--(1.2,-1.3)--cycle;
\draw (0.1,-2.9)--(0.1,-3.2) (0.7,-2.9)--(0.7,-3.2) (1,-2.9)--(1,-3.2);
\draw[double,thin] (0.4,-2.9)--(0.4,-3.2);
\draw (0.4,-2.3)--(0.4,-0.3) node[left,midway,scale=0.7]{$(3,0)$};

\draw[double,thin] (1,1.3) ..controls (1.8,1.3) ..(1.8,3.2);
\draw (1,2.3)--(1,0.3);
\draw (0.1,2.9)--(0.1,3.2) (0.7,2.9)--(0.7,3.2) (1,2.9)--(1,3.2);
\draw[double,thin] (0.4,2.9)--(0.4,3.2);
\draw (0.4,2.3)--(0.4,0.3) node[left,midway,scale=0.7]{$(2,1)$};

\draw [darkred,thick] (-0.3,-2.3) rectangle (1.4,-2.9) node[midway,black]{$(3,1)$};
\draw [darkred,thick] (-0.3,-0.3) rectangle (1.4,0.3) node[midway,black]{$(3,1)$};
\draw [darkred,thick] (-0.3,2.3) rectangle (1.4,2.9) node[midway,black]{$(3,1)$};
\end{tric}
}

\def\dlexj{
\begin{tric}
\draw[double,thin] (1,-1.3) ..controls (1.8,-1.3) ..(1.8,-3.2);
\draw  (0.7,-0.3)--(0.7,-0.9)--(1,-0.9)--(1,-2.3);
\draw (0.7,-2.3)--(0.7,-1.7)--(0.4,-1.7)--(0.4,-0.3);
\draw[double,thin] (0.7,-1.7)--(0.7,-0.9) (1,-0.9)--(1,-0.3) (0.4,-1.7)--(0.4,-2.3);
\draw (0.1,-2.3)--(0.1,-0.3);
\draw (0.1,-2.9)--(0.1,-3.2) (0.7,-2.9)--(0.7,-3.2) (1,-2.9)--(1,-3.2);
\draw[double,thin] (0.4,-2.9)--(0.4,-3.2);

\draw[double,thin] (1.2,1.3) ..controls (1.8,1.3) ..(1.8,3.2);
\draw [double,thin](1,0.3)--(1,1.15);
\draw [double,thin](1,1.7)--(1,1.45);
\draw (1,1.15)--(1,1.45)--(1.2,1.3)--cycle;
\draw (1,2.3)--(1,1.7)--(0.7,1.7)--(0.7,0.3);
\draw (0.7,2.3)--(0.7,2)--(0.4,2)--(0.4,0.3);
\draw [double,thin] (0.4,2)--(0.4,2.3) (0.7,2)--(0.7,1.7);
\draw (0.1,2.3)--(0.1,0.3);
\draw (0.1,2.9)--(0.1,3.2) (0.7,2.9)--(0.7,3.2) (1,2.9)--(1,3.2);
\draw[double,thin] (0.4,2.9)--(0.4,3.2);

\draw [darkred,thick] (-0.3,-2.3) rectangle (1.4,-2.9) node[midway,black]{$C_{\unw}$};
\draw [darkred,thick] (-0.3,-0.3) rectangle (1.4,0.3) node[midway,black]{$C_{\unu_{\lambda+\mu}}$};
\draw [darkred,thick] (-0.3,2.3) rectangle (1.4,2.9) node[midway,black]{$C_{\unw}$};
\end{tric}
}

\def\dlexk{
\begin{tric}
\draw[double,thin] (1.2,1.3) ..controls (1.8,1.3) ..(1.8,3.2);
\draw [double,thin](1,0.3)--(1,1.15);
\draw [double,thin](1,2.3)--(1,1.45);
\draw (1,1.15)--(1,1.45)--(1.2,1.3)--cycle;
\draw (0.1,-2.9)--(0.1,-3.2) (0.7,-2.9)--(0.7,-3.2) (1,-2.9)--(1,-3.2);
\draw[double,thin] (0.4,-2.9)--(0.4,-3.2);
\draw (0.4,2.3)--(0.4,0.3) node[left,midway,scale=0.7]{$(3,0)$};

\draw[double,thin] (1,-1.3) ..controls (1.8,-1.3) ..(1.8,-3.2);
\draw (1,-2.3)--(1,-0.3);
\draw (0.1,2.9)--(0.1,3.2) (0.7,2.9)--(0.7,3.2) (1,2.9)--(1,3.2);
\draw[double,thin] (0.4,2.9)--(0.4,3.2);
\draw (0.4,-2.3)--(0.4,-0.3) node[left,midway,scale=0.7]{$(2,1)$};

\draw [darkred,thick] (-0.3,-2.3) rectangle (1.4,-2.9) node[midway,black]{$(3,1)$};
\draw [darkred,thick] (-0.3,-0.3) rectangle (1.4,0.3) node[midway,black]{$(3,1)$};
\draw [darkred,thick] (-0.3,2.3) rectangle (1.4,2.9) node[midway,black]{$(3,1)$};
\end{tric}
}

\def\dlexl{
\begin{tric}
\draw[double,thin] (1.2,-1.3) ..controls (1.8,-1.3) ..(1.8,-3.2);
\draw [double,thin](1,-0.3)--(1,-1.15);
\draw [double,thin](1,-1.7)--(1,-1.45);
\draw (1,-1.15)--(1,-1.45)--(1.2,-1.3)--cycle;
\draw (1,-2.3)--(1,-1.7)--(0.7,-1.7)--(0.7,-0.3);
\draw (0.7,-2.3)--(0.7,-2)--(0.4,-2)--(0.4,-0.3);
\draw [double,thin] (0.4,-2)--(0.4,-2.3) (0.7,-2)--(0.7,-1.7);
\draw (0.1,-2.3)--(0.1,-0.3);
\draw (0.1,-2.9)--(0.1,-3.2) (0.7,-2.9)--(0.7,-3.2) (1,-2.9)--(1,-3.2);
\draw[double,thin] (0.4,-2.9)--(0.4,-3.2);

\draw[double,thin] (1.2,1.3) ..controls (1.8,1.3) ..(1.8,3.2);
\draw [double,thin](1,0.3)--(1,1.15);
\draw [double,thin](1,1.7)--(1,1.45);
\draw (1,1.15)--(1,1.45)--(1.2,1.3)--cycle;
\draw (1,2.3)--(1,1.7)--(0.7,1.7)--(0.7,0.3);
\draw (0.7,2.3)--(0.7,2)--(0.4,2)--(0.4,0.3);
\draw [double,thin] (0.4,2)--(0.4,2.3) (0.7,2)--(0.7,1.7);
\draw (0.1,2.3)--(0.1,0.3);
\draw (0.1,2.9)--(0.1,3.2) (0.7,2.9)--(0.7,3.2) (1,2.9)--(1,3.2);
\draw[double,thin] (0.4,2.9)--(0.4,3.2);

\draw [darkred,thick] (-0.3,-2.3) rectangle (1.4,-2.9) node[midway,black]{$C_{\unw}$};
\draw [darkred,thick] (-0.3,-0.3) rectangle (1.4,0.3) node[midway,black]{$C_{\unu_{\lambda+\mu}}$};
\draw [darkred,thick] (-0.3,2.3) rectangle (1.4,2.9) node[midway,black]{$C_{\unw}$};
\end{tric}
}

\def\dlexm{
\begin{tric}
\draw[double,thin] (1.2,1.3) ..controls (1.8,1.3) ..(1.8,3.2);
\draw [double,thin](1,0.3)--(1,1.15);
\draw [double,thin](1,2.3)--(1,1.45);
\draw (1,1.15)--(1,1.45)--(1.2,1.3)--cycle;
\draw (0.1,2.9)--(0.1,3.2) (0.7,2.9)--(0.7,3.2) (1,2.9)--(1,3.2);
\draw[double,thin] (0.4,2.9)--(0.4,3.2);
\draw (0.4,2.3)--(0.4,0.3) node[left,midway,scale=0.7]{$(3,0)$};

\draw[double,thin] (1.2,-1.3) ..controls (1.8,-1.3) ..(1.8,-3.2);
\draw [double,thin](1,-0.3)--(1,-1.15);
\draw [double,thin](1,-2.3)--(1,-1.45);
\draw (1,-1.15)--(1,-1.45)--(1.2,-1.3)--cycle;
\draw (0.1,-2.9)--(0.1,-3.2) (0.7,-2.9)--(0.7,-3.2) (1,-2.9)--(1,-3.2);
\draw[double,thin] (0.4,-2.9)--(0.4,-3.2);
\draw (0.4,-2.3)--(0.4,-0.3) node[left,midway,scale=0.7]{$(3,0)$};

\draw [darkred,thick] (-0.3,-2.3) rectangle (1.4,-2.9) node[midway,black]{$(3,1)$};
\draw [darkred,thick] (-0.3,-0.3) rectangle (1.4,0.3) node[midway,black]{$(3,1)$};
\draw [darkred,thick] (-0.3,2.3) rectangle (1.4,2.9) node[midway,black]{$(3,1)$};
\end{tric}
}

\begin{notation}
When $\wt\unx=\wt\uny=(a,b)$, we will use an $(a, b)$ labelled box in $\Hom_{\DD}(\unx, \uny)$ to denote $C_{\unx}^{\uny}$.
\end{notation}

\begin{example}\label{exampleladders}
Let $\unw=\varpi_1\varpi_2\varpi_1\varpi_1  ,\ \varpi=\varpi_2$. When $\mu= (0,0)$,  choose $\unu_{\lambda+\mu}=\varpi_1\varpi_1\varpi_1\varpi_2 $, then (using Lemma \ref{neturalabsorption}) we have
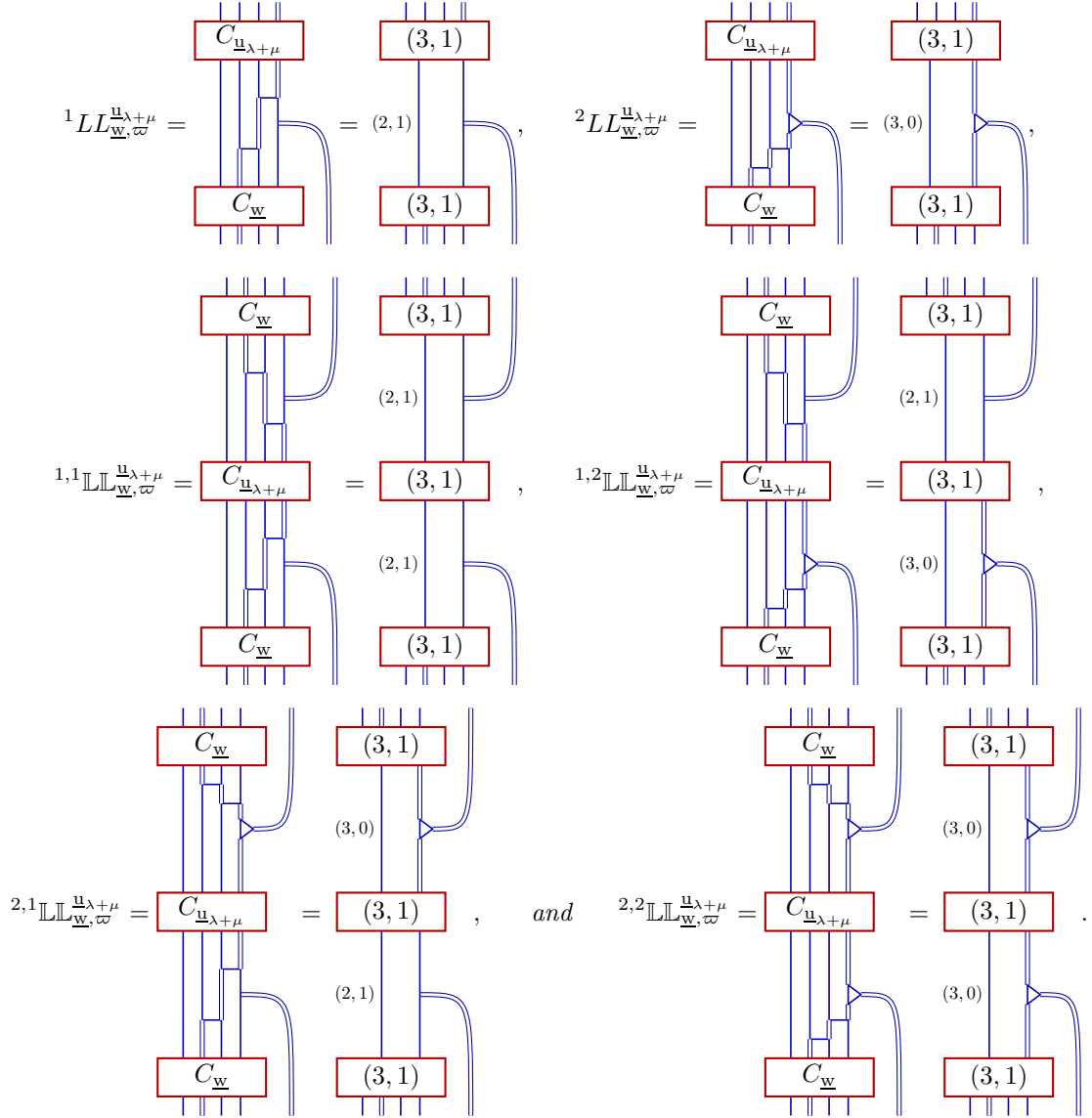
\[
^1LL_{\unw, \varpi}^{\unu_{\lambda+ \mu}}=\dlexd=\dlexda, \ \ \ \ \ ^2LL_{\unw, \varpi}^{\unu_{\lambda+ \mu}}=\dlexe=\dlexea,
\]
\[
^{1,1}\mathbb{LL}_{\unw, \varpi}^{\unu_{\lambda+ \mu}}= \dlexf=\dlexg, \ \ \ \ \ ^{1,2}\mathbb{LL}_{\unw, \varpi}^{\unu_{\lambda+ \mu}}= \dlexh=\dlexi,
\]
\[
^{2,1}\mathbb{LL}_{\unw, \varpi}^{\unu_{\lambda+ \mu}}= \dlexj=\dlexk, \ \ \ \ \ \text{and} \ \ \ \ \  ^{2,2}\mathbb{LL}_{\unw, \varpi}^{\unu_{\lambda+ \mu}}= \dlexl=\dlexm.
\]
\end{example}

\begin{remark}
As $(a, b)$ varies, so does $\Kab(V(\varpi))$. However, the vector $v_{\varpi} = v_{\varpi, 1}$ is always contained in $\Kab(V(\varpi))$. Moreover, the associated elementary light ladder is just a composition of neutral diagrams, so the associated (clasped) double ladder can be simplified to
\[
\mathbb{LL}_{\unw, \varpi}^{\unu_{\lambda+ \varpi}}= C_{\unw\otimes \varpi}.
\]
Thus, the (clasped) double ladder associated to the highest weight vector in $V(\varpi)$ is itself a clasp.
\end{remark}

\newpage
\enlargethispage{30pt}

\section{Triple Clasp Formula}
\label{FormulaProof}
\subsection{Formulas}
\label{results}

\def\gcalcaa{
\begin{tric}
\draw (0.75,0.2)--(0.75,2.5) node[left,midway,scale=0.7] {$(a,b)$}
       (0.75,-0.2)--(0.75,-2.5) node[left,midway,scale=0.7] {$(a,b)$}
       (1.8,-2.5)--(1.8,2.5);

\draw[darkred,thick](0,0.2)rectangle(1.3,-0.2) 
                     node[pos=0.5,black,scale=0.7]{$(a,b)$};
\end{tric}
}

\def\gcalcab{
\begin{tric}
\draw (0.9,0.2)--(0.9,0.6) node[left,midway,scale=0.7] {$(a+p,b+\ell)$}
       (0.9,-0.2)--(0.9,-0.6) node[left,midway,scale=0.7] {$(a+p,b+\ell)$}
       (0.3,1.4)--(0.3,1.8) node[left,midway,scale=0.7] {$(a,b)$}
       (0.3,2.2)--(0.3,2.7) node[left,midway,scale=0.7] {$(a,b)$}
       (0.3,-1.4)--(0.3,-1.8) node[left,midway,scale=0.7] {$(a,b)$}
       (0.3,-2.2)--(0.3,-2.7) node[left,midway,scale=0.7] {$(a,b)$}
       (1.4,2.7)--(1.4,1.4) (1.4,-2.7)--(1.4,-1.4);
       
\draw[darkred,thick](-0.5,0.2)rectangle(1.7,-0.2) 
                     node[pos=0.5,black,scale=0.7]{$(a+p,b+\ell)$}
                     (-0.3,-1.8)rectangle(1,-2.2) 
                     node[pos=0.5,black,scale=0.7]{$(a,b)$}
                     (-0.3,1.8)rectangle(1,2.2) 
                     node[pos=0.5,black,scale=0.7]{$(a,b)$}
                     (-0.3,1.4)--(1.7,1.4)--(1.4,0.6)--(0,0.6)--cycle
                     (0.7,1) node[black]{$\overline{L}_{(a,b)}^{(p,\ell)}$}
                     (-0.3,-1.4)--(1.7,-1.4)--(1.4,-0.6)--(0,-0.6)--cycle
                     (0.7,-1)node[black]{$L_{(a,b)}^{(p,\ell)}$};
\end{tric}
}

\def\gcalcac{
\begin{tric}
\draw (0.75,0.2)--(0.75,2.5) node[left,midway,scale=0.7] {$(a,b)$}
       (0.75,-0.2)--(0.75,-2.5) node[left,midway,scale=0.7] {$(a,b)$}; 
\draw[double,thin] (1.8,-2.5)--(1.8,2.5);

\draw[darkred,thick](0,0.2)rectangle(1.3,-0.2) 
                     node[pos=0.5,black,scale=0.7]{$(a,b)$};
\end{tric}
}

\def\gcalcad{
\begin{tric}
\draw (0.9,0.2)--(0.9,0.6) node[left,midway,scale=0.7] {$(a+p,b+\ell)$}
       (0.9,-0.2)--(0.9,-0.6) node[left,midway,scale=0.7] {$(a+p,b+\ell)$}
       (0.3,1.4)--(0.3,1.8) node[left,midway,scale=0.7] {$(a,b)$}
       (0.3,2.2)--(0.3,2.7) node[left,midway,scale=0.7] {$(a,b)$}
       (0.3,-1.4)--(0.3,-1.8) node[left,midway,scale=0.7] {$(a,b)$}
       (0.3,-2.2)--(0.3,-2.7) node[left,midway,scale=0.7] {$(a,b)$};
\draw[double,thin] (1.4,2.7)--(1.4,1.4) (1.4,-2.7)--(1.4,-1.4);

\draw[darkred,thick](-0.5,0.2)rectangle(1.7,-0.2) 
                     node[pos=0.5,black,scale=0.7]{$(a+p,b+\ell)$}
                     (-0.3,-1.8)rectangle(1,-2.2) 
                     node[pos=0.5,black,scale=0.7]{$(a,b)$}
                     (-0.3,1.8)rectangle(1,2.2) 
                     node[pos=0.5,black,scale=0.7]{$(a,b)$}
                     (-0.3,1.4)--(1.7,1.4)--(1.4,0.6)--(0,0.6)--cycle
                (0.7,1) node[black]{$\overline{T}_{(a,b)}^{(p,\ell)}$}
                     (-0.3,-1.4)--(1.7,-1.4)--(1.4,-0.6)--(0,-0.6)--cycle
                     (0.7,-1)node[black]{$T_{(a,b)}^{(p,\ell)}$};
\end{tric}
}

\def\gcalcae{
\begin{tric}
\draw[double,thin] (1,0.6) ..controls (1.8,0.6) ..(1.8,2)
                   (1,-0.6) ..controls (1.8,-0.6) ..(1.8,-2);
\draw  (1,0.2)--(1,0.6)--(1,1) (1,-0.2)--(1,-0.6)--(1,-1);
\draw (0.75,2)--(0.75,1.4) node[left,midway,black,scale=0.7] {$(a,b)$};
\draw (0.75,-2)--(0.75,-1.4) node[left,midway,black,scale=0.7] {$(a,b)$};;
\draw (0.5,1)--(0.5,0.2) node[left,midway,black,scale=0.7] {$(a-1,b)$};
\draw (0.5,-1)--(0.5,-0.2) node[left,midway,black,scale=0.7] {$(a-1,b)$};
\draw[darkred,thick] (-0.3,-0.2) rectangle (1.5,0.2) node[pos=0.5,scale=0.7,black] {$(a,b)$};
\draw[darkred,thick] (-0.3,1) rectangle (1.5,1.4) 
node[pos=0.5,scale=0.7,black] {$(a,b)$};
\draw[darkred, thick] (-0.3,-1.4) rectangle (1.5,-1) node[pos=0.5,scale=0.7,black] {$(a,b)$};
\end{tric}
}

\def\gcalcaf{
\begin{tric}
\draw[double,thin] (1,0.6) ..controls (1.8,0.6) ..(1.8,2)
                    (1.2,-0.6) ..controls (1.8,-0.6) ..(1.8,-2);
\draw  (1,0.2)--(1,0.6)--(1,1);
\draw [double,thin](1,-0.2)--(1,-0.45) (1,-1)--(1,-0.75);
\draw (1,-0.45)--(1,-0.75)--(1.2,-0.6)--cycle;
\draw (0.75,2)--(0.75,1.4)node[left,midway,black,scale=0.7] {$(a,b)$};
\draw (0.75,-2)--(0.75,-1.4)node[left,midway,black,scale=0.7] {$(a,b)$};
\draw (0.5,1)--(0.5,0.2) node[left,midway,black,scale=0.7] {$(a-1,b)$};
\draw (0.5,-1)--(0.5,-0.2)node[left,midway,black,scale=0.7] {$(a,b-1)$};
\draw[darkred,thick] (-0.3,-0.2) rectangle (1.5,0.2) node[pos=0.5,scale=0.7,black] {$(a,b)$};
\draw[darkred,thick] (-0.3,1) rectangle (1.5,1.4) node[pos=0.5,scale=0.7,black] {$(a,b)$};
\draw[darkred, thick] (-0.3,-1.4) rectangle (1.5,-1) node[pos=0.5,scale=0.7,black] {$(a,b)$};
\end{tric}
}

\def\gcalcag{
\begin{tric}
\draw[double,thin] (1,-0.6) ..controls (1.8,-0.6) ..(1.8,-2)
                    (1.2,0.6) ..controls (1.8,0.6) ..(1.8,2);
\draw  (1,-0.2)--(1,-0.6)--(1,-1);
\draw [double,thin](1,0.2)--(1,0.45) (1,1)--(1,0.75);
\draw (1,0.45)--(1,0.75)--(1.2,0.6)--cycle;
\draw (0.75,2)--(0.75,1.4)node[left,midway,black,scale=0.7] {$(a,b)$};
\draw (0.75,-2)--(0.75,-1.4)node[left,midway,black,scale=0.7] {$(a,b)$};
\draw (0.5,-1)--(0.5,-0.2) node[left,midway,black,scale=0.7] {$(a-1,b)$};
\draw (0.5,1)--(0.5,0.2)node[left,midway,black,scale=0.7] {$(a,b-1)$};
\draw[darkred,thick] (-0.3,-0.2) rectangle (1.5,0.2) node[pos=0.5,scale=0.7,black] {$(a,b)$};
\draw[darkred,thick] (-0.3,1) rectangle (1.5,1.4) node[pos=0.5,scale=0.7,black] {$(a,b)$};
\draw[darkred, thick] (-0.3,-1.4) rectangle (1.5,-1) node[pos=0.5,scale=0.7,black] {$(a,b)$};
\end{tric}
}

\def\gcalcah{
\begin{tric}
\draw[double,thin] (1.2,-0.6) ..controls (1.8,-0.6) ..(1.8,-2)
                    (1.2,0.6) ..controls (1.8,0.6) ..(1.8,2);
\draw [double,thin](1,0.2)--(1,0.45) (1,1)--(1,0.75)
                    (1,-0.2)--(1,-0.45) (1,-1)--(1,-0.75);
\draw (1,0.45)--(1,0.75)--(1.2,0.6)--cycle
      (1,-0.45)--(1,-0.75)--(1.2,-0.6)--cycle;
\draw (0.75,2)--(0.75,1.4)node[left,midway,black,scale=0.7] {$(a,b)$};
\draw (0.75,-2)--(0.75,-1.4)node[left,midway,black,scale=0.7] {$(a,b)$};
\draw (0.5,-1)--(0.5,-0.2) node[left,midway,black,scale=0.7] {$(a,b-1)$};
\draw (0.5,1)--(0.5,0.2)node[left,midway,black,scale=0.7] {$(a,b-1)$};
\draw[darkred,thick] (-0.3,-0.2) rectangle (1.5,0.2) node[pos=0.5,scale=0.7,black] {$(a,b)$};
\draw[darkred,thick] (-0.3,1) rectangle (1.5,1.4) node[pos=0.5,scale=0.7,black] {$(a,b)$};
\draw[darkred, thick] (-0.3,-1.4) rectangle (1.5,-1) node[pos=0.5,scale=0.7,black] {$(a,b)$};
\end{tric}
}

\def\implicitLLa{
\begin{tric}
\draw (0.9,0.1)--(0.9,0.6) node[left,midway,scale=0.7] {$(a+p,b+\ell)$}
       (0.3,1.4)--(0.3,1.9) node[left,midway,scale=0.7] {$(a,b)$}
       (1.4,1.9)--(1.4,1.4);
       
\draw[darkred,thick] (-0.3,1.4)--(1.7,1.4)--(1.4,0.6)--(0,0.6)--cycle
                     (0.7,1) node[black]{$\overline{L}_{(a,b)}^{(p,\ell)}$};
\end{tric}
}

\def\implicitLLb{
\begin{tric}
\draw  (0.9,-0.1)--(0.9,-0.6) node[left,midway,scale=0.7] {$(a+p,b+\ell)$}  
       (0.3,-1.4)--(0.3,-1.9) node[left,midway,scale=0.7] {$(a,b)$}
       (1.4,-1.9)--(1.4,-1.4);
       
\draw[darkred,thick] (-0.3,-1.4)--(1.7,-1.4)--(1.4,-0.6)--(0,-0.6)--cycle
                     (0.7,-1)node[black]{$L_{(a,b)}^{(p,\ell)}$};
\end{tric}
}

\def\implicitLLc{
\begin{tric}
\draw (0.9,0.1)--(0.9,0.6) node[left,midway,scale=0.7] {$(a+p,b+\ell)$}
       (0.3,1.4)--(0.3,1.9) node[left,midway,scale=0.7] {$(a,b)$};
\draw[double,thin]       (1.4,1.9)--(1.4,1.4);
       
\draw[darkred,thick] (-0.3,1.4)--(1.7,1.4)--(1.4,0.6)--(0,0.6)--cycle
                     (0.7,1) node[black]{$\overline{T}_{(a,b)}^{(p,\ell)}$};
\end{tric}
}

\def\implicitLLd{
\begin{tric}
\draw  (0.9,-0.1)--(0.9,-0.6) node[left,midway,scale=0.7] {$(a+p,b+\ell)$}  
       (0.3,-1.4)--(0.3,-1.9) node[left,midway,scale=0.7] {$(a,b)$};
\draw[double,thin]       (1.4,-1.9)--(1.4,-1.4);
       
\draw[darkred,thick] (-0.3,-1.4)--(1.7,-1.4)--(1.4,-0.6)--(0,-0.6)--cycle
                     (0.7,-1)node[black]{$T_{(a,b)}^{(p,\ell)}$};
\end{tric}
}

\def\Ka{
\begin{tric}
\draw (1.5,2)--(1.5,-2);
\draw (0.6,2)--(0.6,0.2);
\draw (0.6,1) node[anchor=east,black,scale=0.7]{$(a,b)$};
\draw (0.6,-2)--(0.6,-0.2);
\draw (0.6,-1) node[anchor=east,black,scale=0.7]{$(a,b)$};
\draw[darkred,thick] (0,-0.2) rectangle (1.2,0.2) node[pos=0.5,scale=0.7,black] {$(a,b)$};
\end{tric}}

\def\Kb{
\begin{tric}
\draw (5.2,0.8) ..controls (5.2,0.5)and(5.9,-0.2) ..(5.9,2);
\draw (5.2,-0.8) ..controls (5.2,-0.5)and(5.9,0.2) ..(5.9,-2);
\draw (4.75,1.2)--(4.75,2);
\draw (4.5,0.8)--(4.5,0.2);
\draw (4.75,-1.2)--(4.75,-2);
\draw (4.5,-0.8)--(4.5,-0.2);
\draw (4.75,1.6) node[anchor=east,black,scale=0.7]  {$(a,b)$};
\draw (4.5,0.5) node[anchor=east,black,scale=0.7] {$(a-1,b)$};
\draw (4.75,-1.6) node[anchor=east,black,scale=0.7] {$(a,b)$};
\draw (4.5,-0.5) node[anchor=east,black,scale=0.7] {$(a-1,b)$};
\draw[darkred,thick](3.7,-0.2) rectangle (5.5,0.2)node[pos=0.5,scale=0.7,black] {$(a-1,b)$};
\draw[darkred,thick](3.7,0.8) rectangle (5.5,1.2)node[pos=0.5,scale=0.7,black] {$(a,b)$};
\draw[darkred,thick](3.7,-1.2) rectangle
(5.5,-0.8)node[pos=0.5,scale=0.7,black] {$(a,b)$};
\end{tric}}

\def\Kc{
\begin{tric}
\draw (9,0.5) ..controls (9.8,0.5) ..(9.8,2);
\draw  (9,0.2)--(9,0.5);
\draw[double,thin] (9,0.8)--(9,0.5);
\draw (9,-0.5) ..controls (9.8,-0.5) ..(9.8,-2);
\draw  (9,-0.2)--(9,-0.5);
\draw[double,thin] (9,-0.8)--(9,-0.5);
\draw (8.75,2)--(8.75,1.2);
\draw (8.75,-2)--(8.75,-1.2);
\draw (8.5,0.8)--(8.5,0.2);
\draw (8.5,-0.8)--(8.5,-0.2);
\draw (8.75,1.6) node[anchor=east,black,scale=0.7]  {$(a,b)$};
\draw (8.5,0.5) node[anchor=east,black,scale=0.7] {$(a,b-1)$};
\draw (8.75,-1.6) node[anchor=east,black,scale=0.7] {$(a,b)$};
\draw (8.5,-0.5) node[anchor=east,black,scale=0.7] {$(a,b-1)$};

\draw[darkred, thick](7.7,-0.2) rectangle (9.5,0.2)node[pos=0.5,scale=0.7,black] {$(a+1,b-1)$};;
\draw[darkred, thick](7.7,0.8) rectangle (9.5,1.2)node[pos=0.5,scale=0.7,black] {$(a,b)$};;
\draw[darkred, thick](7.7,-1.2) rectangle (9.5,-0.8)node[pos=0.5,scale=0.7,black] {$(a,b)$};;
\end{tric}}

\def\Kd{
\begin{tric}
\draw (13.2,0.5) ..controls (13.8,0.5) ..(13.8,2);
\draw (13.2,0.5) --(13.2,0.8);
\draw (13.2,0.5) -- (12.8,0.5);
\draw[double,thin] (12.8,0.5) --(12.8,0.2);
\draw (12.8,0.5) --(12.8,0.8);
\draw (13.2,-0.5) ..controls (13.8,-0.5) ..(13.8,-2);
\draw (13.2,-0.5) --(13.2,-0.8);
\draw (13.2,-0.5) -- (12.8,-0.5);
\draw[double,thin] (12.8,-0.5) --(12.8,-0.2);
\draw (12.8,-0.5) --(12.8,-0.8);
\draw (12.7,2)--(12.7,1.2);
\draw (12.7,-2)--(12.7,-1.2);
\draw (12.3,0.8)--(12.3,0.2);
\draw (12.3,-0.8)--(12.3,-0.2);
\draw (12.7,1.6) node[anchor=east,black,scale=0.7]  {$(a,b)$};
\draw (12.3,0.5) node[anchor=east,black,scale=0.7] {$(a-2,b)$};
\draw (12.7,-1.6) node[anchor=east,black,scale=0.7] {$(a,b)$};
\draw (12.3,-0.5) node[anchor=east,black,scale=0.7] {$(a-2,b)$};

\draw[darkred,thick](11.7,-0.2) rectangle (13.5,0.2)node[pos=0.5,scale=0.7,black] {$(a-2,b+1)$};
\draw[darkred,thick](11.7,0.8) rectangle (13.5,1.2)node[pos=0.5,scale=0.7,black] {$(a,b)$};
\draw[darkred,thick](11.7,-1.2) rectangle (13.5,-0.8)node[pos=0.5,scale=0.7,black] {$(a,b)$};
\end{tric}}
 
\def\Ke{
\begin{tric}
\draw (1,0.5) ..controls (1.8,0.5) ..(1.8,2);
\draw  (1,0.2)--(1,0.5);
\draw (1,0.8)--(1,0.5);
\draw (1,-0.5) ..controls (1.8,-0.5) ..(1.8,-2);
\draw  (1,-0.2)--(1,-0.5);
\draw (1,-0.8)--(1,-0.5);
\draw (0.75,2)--(0.75,1.2);
\draw (0.75,-2)--(0.75,-1.2);
\draw (0.5,0.8)--(0.5,0.2);
\draw (0.5,-0.8)--(0.5,-0.2);
\draw (0.75,1.6) node[anchor=east,black,scale=0.7] {$(a,b)$};
\draw (0.5,0.5) node[anchor=east,black,scale=0.7] {$(a-1,b)$};
\draw (0.75,-1.6) node[anchor=east,black,scale=0.7] {$(a,b)$};
\draw (0.5,-0.5) node[anchor=east,black,scale=0.7] {$(a-1,b)$};

\draw[darkred,thick] (-0.3,-0.2) rectangle (1.5,0.2) node[pos=0.5,scale=0.7,black] {$(a,b)$};
\draw[darkred,thick] (-0.3,0.8) rectangle (1.5,1.2) node[pos=0.5,scale=0.7,black] {$(a,b)$};
\draw[darkred, thick] (-0.3,-1.2) rectangle (1.5,-0.8) node[pos=0.5,scale=0.7,black] {$(a,b)$};
\end{tric}}

\def\Kf{
\begin{tric}
\draw (5.2,0.5) ..controls (5.8,0.5) ..(5.8,2);
\draw (5.2,0.5) --(5.2,0.2);
\draw (5.2,0.5) -- (4.8,0.5);
\draw (4.8,0.5) --(4.8,0.2);
\draw[double,thin] (4.8,0.5) --(4.8,0.8);
\draw (5.2,-0.5) ..controls (5.8,-0.5) ..(5.8,-2);
\draw (5.2,-0.5) --(5.2,-0.2);
\draw (5.2,-0.5) -- (4.8,-0.5);
\draw (4.8,-0.5) --(4.8,-0.2);
\draw[double,thin] (4.8,-0.5) --(4.8,-0.8);
\draw (4.7,2)--(4.7,1.2);
\draw (4.7,-2)--(4.7,-1.2);
\draw (4.3,0.8)--(4.3,0.2);
\draw (4.3,-0.8)--(4.3,-0.2);
\draw (4.7,1.6) node[anchor=east,black,scale=0.7]  {$(a,b)$};
\draw (4.3,0.5) node[anchor=east,black,scale=0.7] {$(a,b-1)$};
\draw (4.7,-1.6) node[anchor=east,black,scale=0.7] {$(a,b)$};
\draw (4.3,-0.5) node[anchor=east,black,scale=0.7] {$(a,b-1)$};

\draw[darkred,thick](3.7,-0.2) rectangle (5.5,0.2) node[pos=0.5,scale=0.7,black] {$(a+2,b-1)$};
\draw[darkred,thick](3.7,0.8) rectangle (5.5,1.2)
node[pos=0.5,scale=0.7,black] {$(a,b)$};
\draw[darkred,thick](3.7,-1.2) rectangle (5.5,-0.8) node[pos=0.5,scale=0.7,black] {$(a,b)$};
\end{tric}}

\def\Kg {
\begin{tric}
\draw (9,0.5) ..controls (9.8,0.5) ..(9.8,2);
\draw[double,thin]  (9,0.2)--(9,0.5);
\draw (9,0.8)--(9,0.5);
\draw (9,-0.5) ..controls (9.8,-0.5) ..(9.8,-2);
\draw[double,thin]  (9,-0.2)--(9,-0.5);
\draw (9,-0.8)--(9,-0.5);
\draw (8.75,2)--(8.75,1.2);
\draw (8.75,-2)--(8.75,-1.2);
\draw (8.5,0.8)--(8.5,0.2);
\draw (8.5,-0.8)--(8.5,-0.2);
\draw (8.75,1.6) node[anchor=east,black,scale=0.7] {$(a,b)$};
\draw (8.5,0.5) node[anchor=east,black,scale=0.7] {$(a-1,b)$};
\draw (8.75,-1.6) node[anchor=east,black,scale=0.7] {$(a,b)$};
\draw (8.5,-0.5) node[anchor=east,black,scale=0.7] {$(a-1,b)$};

\draw[darkred, thick](7.7,-0.2) rectangle (9.5,0.2) node[pos=0.5,scale=0.7,black] {$(a-1,b+1)$};
\draw[darkred, thick](7.7,0.8) rectangle (9.5,1.2) node[pos=0.5,scale=0.7,black] {$(a,b)$};
\draw[darkred,thick](7.7,-1.2) rectangle (9.5,-0.8) node[pos=0.5,scale=0.7,black] {$(a,b)$};
\end{tric}}

\def\Kh {
\begin{tric}
\draw (13.3,0.2)..controls(13.3,0.6)and(13.8,0.1)..(13.8,2);
\draw (13.3,-0.2)..controls(13.3,-0.6)and(13.8,-0.1)..(13.8,-2);
\draw (12.75,2)--(12.75,1.2);
\draw (12.75,-2)--(12.75,-1.2);
\draw (12.5,0.8)--(12.5,0.2);
\draw (12.5,-0.8)--(12.5,-0.2);
\draw (12.75,1.6) node[anchor=east,black,scale=0.7] {$(a,b)$};
\draw (12.5,0.5) node[anchor=east,black,scale=0.7] {$(a,b)$};
\draw (12.75,-1.6) node[anchor=east,black,scale=0.7] {$(a,b)$};
\draw (12.5,-0.5) node[anchor=east,black,scale=0.7] {$(a,b)$};

\draw[darkred,thick](11.7,-0.2) rectangle (13.5,0.2) node[pos=0.5,scale=0.7,black] {$(a+1,b)$};
\draw[darkred,thick](11.7,0.8) rectangle (13.5,1.2) node[pos=0.5,scale=0.7,black] {$(a,b)$};
\draw[darkred, thick](11.7,-1.2) rectangle (13.5,-0.8) node[pos=0.5,scale=0.7,black] {$(a,b)$};
\end{tric}}

\def\Ra{
\begin{tric}
\draw [double,thin](1.5,2)--(1.5,-2);
\draw (0.6,2)--(0.6,0.2);
\draw (0.6,1) node[anchor=east,black,scale=0.7]{$(a,b)$};
\draw (0.6,-2)--(0.6,-0.2);
\draw (0.6,-1) node[anchor=east,black,scale=0.7]{$(a,b)$};
\draw[darkred,thick] (0,-0.2) rectangle (1.2,0.2) node[pos=0.5,scale=0.7,black] {$(a,b)$};
\end{tric}}

\def\Rb{
\begin{tric}
\draw[double,thin] (1.2,0.8) ..controls (1.2,0.5)and(1.9,-0.2) ..(1.9,2);
\draw[double,thin] (1.2,-0.8) ..controls (1.2,-0.5)and(1.9,0.2) ..(1.9,-2);
\draw (0.75,1.2)--(0.75,2);
\draw (0.5,0.8)--(0.5,0.2);
\draw (0.75,-1.2)--(0.75,-2);
\draw (0.5,-0.8)--(0.5,-0.2);
\draw (0.75,1.6) node[anchor=east,black,scale=0.7]  {$(a,b)$};
\draw (0.5,0.5) node[anchor=east,black,scale=0.7] {$(a,b-1)$};
\draw (0.75,-1.6) node[anchor=east,black,scale=0.7] {$(a,b)$};
\draw (0.5,-0.5) node[anchor=east,black,scale=0.7] {$(a,b-1)$};
\draw[darkred,thick](0,-0.2) rectangle (1.5,0.2)node[pos=0.5,scale=0.7,black] {$(a,b-1)$};
\draw[darkred,thick](0,0.8) rectangle (1.5,1.2)node[pos=0.5,scale=0.7,black] {$(a,b)$};
\draw[darkred,thick](0,-1.2) rectangle
(1.5,-0.8)node[pos=0.5,scale=0.7,black] {$(a,b)$};
\end{tric}}

\def\Rc{
\begin{tric}
\draw [double,thin](1.3,0.5) ..controls (1.8,0.5) ..(1.8,2);
\draw (1.3,0.5) --(1.3,0.8);
\draw (1.3,0.5) -- (0.8,0.5);
\draw [double,thin](0.8,0.5) --(0.8,0.2);
\draw (0.8,0.5) --(0.8,0.8);
\draw [double,thin](1.3,-0.5)..controls (1.8,-0.5)..(1.8,-2);
\draw (1.3,-0.5) --(1.3,-0.8);
\draw (1.3,-0.5) -- (0.8,-0.5);
\draw [double,thin](0.8,-0.5) --(0.8,-0.2);
\draw (0.8,-0.5) --(0.8,-0.8);
\draw (1.05,0.5)--(1.05,0.8);
\draw (1.05,-0.5)--(1.05,-0.8);
\draw (0.7,2)--(0.7,1.2);
\draw (0.7,-2)--(0.7,-1.2);
\draw (0.3,0.8)--(0.3,0.2);
\draw (0.3,-0.8)--(0.3,-0.2);
\draw (0.7,1.6) node[anchor=east,black,scale=0.7]  {$(a,b)$};
\draw (0.3,0.5) node[anchor=east,black,scale=0.7] {$(a-3,b)$};
\draw (0.7,-1.6) node[anchor=east,black,scale=0.7] {$(a,b)$};
\draw (0.3,-0.5) node[anchor=east,black,scale=0.7] {$(a-3,b)$};
\draw[darkred,thick](-0.3,-0.2) rectangle (1.5,0.2)node[pos=0.5,scale=0.7,black] {$(a-3,b+1)$};
\draw[darkred,thick](-0.3,0.8) rectangle (1.5,1.2)node[pos=0.5,scale=0.7,black] {$(a,b)$};
\draw[darkred,thick](-0.3,-1.2) rectangle (1.5,-0.8)node[pos=0.5,scale=0.7,black] {$(a,b)$};
\end{tric}}

\def\Rd{
\begin{tric}
\draw [double,thin](1.2,0.5) ..controls (1.8,0.5) ..(1.8,2);
\draw (1.2,0.5) --(1.2,0.8);
\draw (1.2,0.5) -- (0.8,0.5);
\draw (0.8,0.5) --(0.8,0.2);
\draw (0.8,0.5) --(0.8,0.8);
\draw [double,thin](1.2,-0.5)..controls (1.8,-0.5)..(1.8,-2);
\draw (1.2,-0.5) --(1.2,-0.8);
\draw (1.2,-0.5) -- (0.8,-0.5);
\draw (0.8,-0.5) --(0.8,-0.2);
\draw (0.8,-0.5) --(0.8,-0.8);
\draw (0.7,2)--(0.7,1.2);
\draw (0.7,-2)--(0.7,-1.2);
\draw (0.3,0.8)--(0.3,0.2);
\draw (0.3,-0.8)--(0.3,-0.2);
\draw (0.7,1.6) node[anchor=east,black,scale=0.7]  {$(a,b)$};
\draw (0.3,0.5) node[anchor=east,black,scale=0.7] {$(a-2,b)$};
\draw (0.7,-1.6) node[anchor=east,black,scale=0.7] {$(a,b)$};
\draw (0.3,-0.5) node[anchor=east,black,scale=0.7] {$(a-2,b)$};
\draw[darkred,thick](-0.3,-0.2) rectangle (1.5,0.2)node[pos=0.5,scale=0.7,black] {$(a-1,b)$};
\draw[darkred,thick](-0.3,0.8) rectangle (1.5,1.2)node[pos=0.5,scale=0.7,black] {$(a,b)$};
\draw[darkred,thick](-0.3,-1.2) rectangle (1.5,-0.8)node[pos=0.5,scale=0.7,black] {$(a,b)$};
\end{tric}}

\def\Re{
\begin{tric}
\draw [double,thin](1.2,0.4)--(1.5,0.5);
\draw [double,thin](1.2,0.8)--(1.2,0.65);
\draw [double,thin](1.8,2)--(1.8,0.65);
\draw (1.5,0.5)--(1.2,0.65)--(1.8,0.65)--cycle;
\draw (1.2,0.4) --(1.2,0.2);
\draw (1.2,0.4) -- (0.8,0.4);
\draw (0.8,0.4) --(0.8,0.2);
\draw (0.8,0.4) --(0.8,0.8);
\draw [double,thin](1.2,-0.4)--(1.5,-0.5);
\draw [double,thin](1.2,-0.8)--(1.2,-0.65);
\draw [double,thin](1.8,-2)--(1.8,-0.65);
\draw (1.5,-0.5)--(1.2,-0.65)--(1.8,-0.65)--cycle;
\draw (1.2,-0.4) --(1.2,-0.2);
\draw (1.2,-0.4) -- (0.8,-0.4);
\draw (0.8,-0.4) --(0.8,-0.2);
\draw (0.8,-0.4) --(0.8,-0.8);
\draw (0.7,2)--(0.7,1.2);
\draw (0.7,-2)--(0.7,-1.2);
\draw (0.3,0.8)--(0.3,0.2);
\draw (0.3,-0.8)--(0.3,-0.2);
\draw (0.7,1.6) node[anchor=east,black,scale=0.7]  {$(a,b)$};
\draw (0.3,0.5) node[anchor=east,black,scale=0.7] {$(a-1,b-1)$};
\draw (0.7,-1.6) node[anchor=east,black,scale=0.7] {$(a,b)$};
\draw (0.3,-0.5) node[anchor=east,black,scale=0.7] {$(a-1,b-1)$};
\draw[darkred,thick](-0.3,-0.2) rectangle (1.5,0.2)node[pos=0.5,scale=0.7,black] {$(a+1,b-1)$};
\draw[darkred,thick](-0.3,0.8) rectangle (1.5,1.2)node[pos=0.5,scale=0.7,black] {$(a,b)$};
\draw[darkred,thick](-0.3,-1.2) rectangle (1.5,-0.8)node[pos=0.5,scale=0.7,black] {$(a,b)$};
\end{tric}
}

\def\Rf{
\begin{tric}
\draw [double,thin](1.2,0.5) ..controls (1.8,0.5) ..(1.8,2);
\draw (1.2,0.5) --(1.2,0.8);
\draw (1.2,0.5) -- (0.8,0.5);
\draw [double,thin](0.8,0.5) --(0.8,0.2);
\draw (0.8,0.5) --(0.8,0.8);
\draw [double,thin](1.2,-0.5)..controls (1.8,-0.5)..(1.8,-2);
\draw (1.2,-0.5) --(1.2,-0.8);
\draw (1.2,-0.5) -- (0.8,-0.5);
\draw [double,thin](0.8,-0.5) --(0.8,-0.2);
\draw (0.8,-0.5) --(0.8,-0.8);
\draw (0.7,2)--(0.7,1.2);
\draw (0.7,-2)--(0.7,-1.2);
\draw (0.3,0.8)--(0.3,0.2);
\draw (0.3,-0.8)--(0.3,-0.2);
\draw (0.7,1.6) node[anchor=east,black,scale=0.7]  {$(a,b)$};
\draw (0.3,0.5) node[anchor=east,black,scale=0.7] {$(a-2,b)$};
\draw (0.7,-1.6) node[anchor=east,black,scale=0.7] {$(a,b)$};
\draw (0.3,-0.5) node[anchor=east,black,scale=0.7] {$(a-2,b)$};
\draw[darkred,thick](-0.3,-0.2) rectangle (1.5,0.2)node[pos=0.5,scale=0.7,black] {$(a-2,b+1)$};
\draw[darkred,thick](-0.3,0.8) rectangle (1.5,1.2)node[pos=0.5,scale=0.7,black] {$(a,b)$};
\draw[darkred,thick](-0.3,-1.2) rectangle (1.5,-0.8)node[pos=0.5,scale=0.7,black] {$(a,b)$};
\end{tric}
}

\def\Rg{
\begin{tric}
\draw [double,thin](1.3,0.4)--(1.5,0.5);
\draw [double,thin](1.3,0.8)--(1.3,0.65);
\draw [double,thin](1.8,2)--(1.8,0.65);
\draw (1.5,0.5)--(1.3,0.65)--(1.8,0.65)--cycle;
\draw (1.3,0.4) --(1.3,0.2);
\draw (1.3,0.4) -- (0.8,0.4);
\draw (0.8,0.4) --(0.8,0.2);
\draw [double,thin](0.8,0.4) --(0.8,0.8);
\draw (1.05,0.4)--(1.05,0.2);
\draw [double,thin](1.3,-0.4)--(1.5,-0.5);
\draw [double,thin](1.3,-0.8)--(1.3,-0.65);
\draw [double,thin](1.8,-2)--(1.8,-0.65);
\draw (1.5,-0.5)--(1.3,-0.65)--(1.8,-0.65)--cycle;
\draw (1.3,-0.4) --(1.3,-0.2);
\draw (1.3,-0.4) -- (0.8,-0.4);
\draw (0.8,-0.4) --(0.8,-0.2);
\draw [double,thin](0.8,-0.4) --(0.8,-0.8);
\draw (0.7,2)--(0.7,1.2);
\draw (0.7,-2)--(0.7,-1.2);
\draw (0.3,0.8)--(0.3,0.2);
\draw (0.3,-0.8)--(0.3,-0.2);
\draw (1.05,-0.4)--(1.05,-0.2);
\draw (0.7,1.6) node[anchor=east,black,scale=0.7]  {$(a,b)$};
\draw (0.3,0.5) node[anchor=east,black,scale=0.7] {$(a,b-2)$};
\draw (0.7,-1.6) node[anchor=east,black,scale=0.7] {$(a,b)$};
\draw (0.3,-0.5) node[anchor=east,black,scale=0.7] {$(a,b-2)$};
\draw[darkred,thick](-0.3,-0.2) rectangle (1.5,0.2)node[pos=0.5,scale=0.7,black] {$(a+3,b-2)$};
\draw[darkred,thick](-0.3,0.8) rectangle (1.5,1.2)node[pos=0.5,scale=0.7,black] {$(a,b)$};
\draw[darkred,thick](-0.3,-1.2) rectangle (1.5,-0.8)node[pos=0.5,scale=0.7,black] {$(a,b)$};
\end{tric}
}

\def\Rh{
\begin{tric}
\draw[double,thin] (1,0.5) ..controls (1.8,0.5) ..(1.8,2);
\draw  (1,0.2)--(1,0.5);
\draw (1,0.8)--(1,0.5);
\draw[double,thin] (1,-0.5) ..controls (1.8,-0.5) ..(1.8,-2);
\draw  (1,-0.2)--(1,-0.5);
\draw (1,-0.8)--(1,-0.5);
\draw (0.75,2)--(0.75,1.2);
\draw (0.75,-2)--(0.75,-1.2);
\draw (0.5,0.8)--(0.5,0.2);
\draw (0.5,-0.8)--(0.5,-0.2);
\draw (0.75,1.6) node[anchor=east,black,scale=0.7] {$(a,b)$};
\draw (0.5,0.5) node[anchor=east,black,scale=0.7] {$(a-1,b)$};
\draw (0.75,-1.6) node[anchor=east,black,scale=0.7] {$(a,b)$};
\draw (0.5,-0.5) node[anchor=east,black,scale=0.7] {$(a-1,b)$};
\draw[darkred,thick] (-0.3,-0.2) rectangle (1.5,0.2) node[pos=0.5,scale=0.7,black] {$(a,b)$};
\draw[darkred,thick] (-0.3,0.8) rectangle (1.5,1.2) node[pos=0.5,scale=0.7,black] {$(a,b)$};
\draw[darkred, thick] (-0.3,-1.2) rectangle (1.5,-0.8) node[pos=0.5,scale=0.7,black] {$(a,b)$};
\end{tric}
}

\def\Ri{
\begin{tric}
\draw[double,thin] (1,0.5) ..controls (1.8,0.5) ..(1.8,2);
\draw  (1,0.2)--(1,0.5);
\draw (1,0.8)--(1,0.5);
\draw[double,thin] (1.2,-0.5) ..controls (1.8,-0.5) ..(1.8,-2);
\draw [double,thin](1,-0.2)--(1,-0.35);
\draw [double,thin](1,-0.8)--(1,-0.65);
\draw (1,-0.35)--(1,-0.65)--(1.2,-0.5)--cycle;
\draw (0.75,2)--(0.75,1.2);
\draw (0.75,-2)--(0.75,-1.2);
\draw (0.5,0.8)--(0.5,0.2);
\draw (0.5,-0.8)--(0.5,-0.2);
\draw (0.75,1.6) node[anchor=east,black,scale=0.7] {$(a,b)$};
\draw (0.5,0.5) node[anchor=east,black,scale=0.7] {$(a-1,b)$};
\draw (0.75,-1.6) node[anchor=east,black,scale=0.7] {$(a,b)$};
\draw (0.5,-0.5) node[anchor=east,black,scale=0.7] {$(a,b-1)$};
\draw[darkred,thick] (-0.3,-0.2) rectangle (1.5,0.2) node[pos=0.5,scale=0.7,black] {$(a,b)$};
\draw[darkred,thick] (-0.3,0.8) rectangle (1.5,1.2) node[pos=0.5,scale=0.7,black] {$(a,b)$};
\draw[darkred, thick] (-0.3,-1.2) rectangle (1.5,-0.8) node[pos=0.5,scale=0.7,black] {$(a,b)$};
\end{tric}
}

\def\Rj{
\begin{tric}
\draw[double,thin] (1.2,0.5) ..controls (1.8,0.5) ..(1.8,2);
\draw [double,thin](1,0.2)--(1,0.35);
\draw [double,thin](1,0.8)--(1,0.65);
\draw (1,0.35)--(1,0.65)--(1.2,0.5)--cycle;
\draw[double,thin] (1,-0.5) ..controls (1.8,-0.5) ..(1.8,-2);
\draw  (1,-0.2)--(1,-0.5);
\draw (1,-0.8)--(1,-0.5);
\draw (0.75,2)--(0.75,1.2);
\draw (0.75,-2)--(0.75,-1.2);
\draw (0.5,0.8)--(0.5,0.2);
\draw (0.5,-0.8)--(0.5,-0.2);
\draw (0.75,1.6) node[anchor=east,black,scale=0.7] {$(a,b)$};
\draw (0.5,0.5) node[anchor=east,black,scale=0.7] {$(a,b-1)$};
\draw (0.75,-1.6) node[anchor=east,black,scale=0.7] {$(a,b)$};
\draw (0.5,-0.5) node[anchor=east,black,scale=0.7] {$(a-1,b)$};
\draw[darkred,thick] (-0.3,-0.2) rectangle (1.5,0.2) node[pos=0.5,scale=0.7,black] {$(a,b)$};
\draw[darkred,thick] (-0.3,0.8) rectangle (1.5,1.2) node[pos=0.5,scale=0.7,black] {$(a,b)$};
\draw[darkred, thick] (-0.3,-1.2) rectangle (1.5,-0.8) node[pos=0.5,scale=0.7,black] {$(a,b)$};
\end{tric}
}

\def\Rk{
\begin{tric}
\draw[double,thin] (1.2,0.5) ..controls (1.8,0.5) ..(1.8,2);
\draw [double,thin](1,0.2)--(1,0.35);
\draw [double,thin](1,0.8)--(1,0.65);
\draw (1,0.35)--(1,0.65)--(1.2,0.5)--cycle;
\draw[double,thin] (1.2,-0.5) ..controls (1.8,-0.5) ..(1.8,-2);
\draw [double,thin](1,-0.2)--(1,-0.35);
\draw [double,thin](1,-0.8)--(1,-0.65);
\draw (1,-0.35)--(1,-0.65)--(1.2,-0.5)--cycle;
\draw (0.75,2)--(0.75,1.2);
\draw (0.75,-2)--(0.75,-1.2);
\draw (0.5,0.8)--(0.5,0.2);
\draw (0.5,-0.8)--(0.5,-0.2);
\draw (0.75,1.6) node[anchor=east,black,scale=0.7] {$(a,b)$};
\draw (0.5,0.5) node[anchor=east,black,scale=0.7] {$(a,b-1)$};
\draw (0.75,-1.6) node[anchor=east,black,scale=0.7] {$(a,b)$};
\draw (0.5,-0.5) node[anchor=east,black,scale=0.7] {$(a,b-1)$};
\draw[darkred,thick] (-0.3,-0.2) rectangle (1.5,0.2) node[pos=0.5,scale=0.7,black] {$(a,b)$};
\draw[darkred,thick] (-0.3,0.8) rectangle (1.5,1.2) node[pos=0.5,scale=0.7,black] {$(a,b)$};
\draw[darkred, thick] (-0.3,-1.2) rectangle (1.5,-0.8) node[pos=0.5,scale=0.7,black] {$(a,b)$};
\end{tric}
}

\def\Rl{
\begin{tric}
\draw [double,thin](1.2,0.65)--(1.3,0.5);
\draw [double,thin](1.4,0.2)--(1.4,0.35);
\draw [double,thin](1.8,2)--(1.8,0.5);
\draw (1.3,0.5)--(1.4,0.35)--(1.8,0.5)--cycle;
\draw (1.2,0.65) --(1.2,0.8);
\draw (1.2,0.65) -- (0.8,0.65);
\draw (0.8,0.65) --(0.8,0.8);
\draw [double,thin](0.8,0.65) --(0.8,0.2);
\draw (1,0.65)--(1,0.8);
\draw [double,thin](1.2,-0.65)--(1.3,-0.5);
\draw [double,thin](1.4,-0.2)--(1.4,-0.35);
\draw [double,thin](1.8,-2)--(1.8,-0.5);
\draw (1.3,-0.5)--(1.4,-0.35)--(1.8,-0.5)--cycle;
\draw (1.2,-0.65) --(1.2,-0.8);
\draw (1.2,-0.65) -- (0.8,-0.65);
\draw (0.8,-0.65) --(0.8,-0.8);
\draw [double,thin](0.8,-0.65) --(0.8,-0.2);
\draw (1,-0.65)--(1,-0.8);
\draw (0.7,2)--(0.7,1.2);
\draw (0.7,-2)--(0.7,-1.2);
\draw (0.3,0.8)--(0.3,0.2);
\draw (0.3,-0.8)--(0.3,-0.2);
\draw (0.7,1.6) node[anchor=east,black,scale=0.7]  {$(a,b)$};
\draw (0.3,0.5) node[anchor=east,black,scale=0.7] {$(a-3,b)$};
\draw (0.7,-1.6) node[anchor=east,black,scale=0.7] {$(a,b)$};
\draw (0.3,-0.5) node[anchor=east,black,scale=0.7] {$(a-3,b)$};
\draw[darkred,thick](-0.3,-0.2) rectangle (1.5,0.2)node[pos=0.5,scale=0.7,black] {$(a-3,b+2)$};
\draw[darkred,thick](-0.3,0.8) rectangle (1.5,1.2)node[pos=0.5,scale=0.7,black] {$(a,b)$};
\draw[darkred,thick](-0.3,-1.2) rectangle (1.5,-0.8)node[pos=0.5,scale=0.7,black] {$(a,b)$};
\end{tric}
}

\def\Rm{
\begin{tric}
\draw [double,thin](1.2,0.5) ..controls (1.8,0.5) ..(1.8,2);
\draw (1.2,0.5) --(1.2,0.2);
\draw (1.2,0.5) -- (0.8,0.5);
\draw (0.8,0.5) --(0.8,0.2);
\draw [double,thin](0.8,0.5) --(0.8,0.8);
\draw [double,thin](1.2,-0.5)..controls (1.8,-0.5)..(1.8,-2);
\draw (1.2,-0.5) --(1.2,-0.2);
\draw (1.2,-0.5) -- (0.8,-0.5);
\draw (0.8,-0.5) --(0.8,-0.2);
\draw [double,thin](0.8,-0.5) --(0.8,-0.8);
\draw (0.7,2)--(0.7,1.2);
\draw (0.7,-2)--(0.7,-1.2);
\draw (0.3,0.8)--(0.3,0.2);
\draw (0.3,-0.8)--(0.3,-0.2);
\draw (0.7,1.6) node[anchor=east,black,scale=0.7]  {$(a,b)$};
\draw (0.3,0.5) node[anchor=east,black,scale=0.7] {$(a,b-1)$};
\draw (0.7,-1.6) node[anchor=east,black,scale=0.7] {$(a,b)$};
\draw (0.3,-0.5) node[anchor=east,black,scale=0.7] {$(a,b-1)$};
\draw[darkred,thick](-0.3,-0.2) rectangle (1.5,0.2)node[pos=0.5,scale=0.7,black] {$(a+2,b-1)$};
\draw[darkred,thick](-0.3,0.8) rectangle (1.5,1.2)node[pos=0.5,scale=0.7,black] {$(a,b)$};
\draw[darkred,thick](-0.3,-1.2) rectangle (1.5,-0.8)node[pos=0.5,scale=0.7,black] {$(a,b)$};
\end{tric}
}

\def\Rn{    
\begin{tric}
\draw [double,thin](1.1,0.6)--(1.3,0.5);
\draw [double,thin](1.4,0.2)--(1.4,0.35);
\draw [double,thin](1.8,2)--(1.8,0.5);
\draw (1.3,0.5)--(1.4,0.35)--(1.8,0.5)--cycle;
\draw (1.1,0.6) --(1.1,0.8);
\draw (1.1,0.6) -- (0.8,0.6);
\draw (0.8,0.6) --(0.8,0.8);
\draw (0.8,0.6) --(0.8,0.2);
\draw [double,thin](1.1,-0.6)--(1.3,-0.5);
\draw [double,thin](1.4,-0.2)--(1.4,-0.35);
\draw [double,thin](1.8,-2)--(1.8,-0.5);
\draw (1.3,-0.5)--(1.4,-0.35)--(1.8,-0.5)--cycle;
\draw (1.1,-0.6) --(1.1,-0.8);
\draw (1.1,-0.6) -- (0.8,-0.6);
\draw (0.8,-0.6) --(0.8,-0.8);
\draw (0.8,-0.6) --(0.8,-0.2);
\draw (0.7,2)--(0.7,1.2);
\draw (0.7,-2)--(0.7,-1.2);
\draw (0.3,0.8)--(0.3,0.2);
\draw (0.3,-0.8)--(0.3,-0.2);
\draw (0.7,1.6) node[anchor=east,black,scale=0.7]  {$(a,b)$};
\draw (0.3,0.5) node[anchor=east,black,scale=0.7] {$(a-2,b)$};
\draw (0.7,-1.6) node[anchor=east,black,scale=0.7] {$(a,b)$};
\draw (0.3,-0.5) node[anchor=east,black,scale=0.7] {$(a-2,b)$};
\draw[darkred,thick](-0.3,-0.2) rectangle (1.5,0.2)node[pos=0.5,scale=0.7,black] {$(a-1,b+1)$};
\draw[darkred,thick](-0.3,0.8) rectangle (1.5,1.2)node[pos=0.5,scale=0.7,black] {$(a,b)$};
\draw[darkred,thick](-0.3,-1.2) rectangle (1.5,-0.8)node[pos=0.5,scale=0.7,black] {$(a,b)$};
\end{tric}
}

\def\Ro{    
\begin{tric}
\draw [double,thin](1.2,0.5) ..controls (1.8,0.5) ..(1.8,2);
\draw (1.2,0.5) --(1.2,0.2);
\draw (1.2,0.5) -- (0.8,0.5);
\draw (0.8,0.5) --(0.8,0.2);
\draw (0.8,0.5) --(0.8,0.8);
\draw [double,thin](1.2,-0.5)..controls (1.8,-0.5)..(1.8,-2);
\draw (1.2,-0.5) --(1.2,-0.2);
\draw (1.2,-0.5) -- (0.8,-0.5);
\draw (0.8,-0.5) --(0.8,-0.2);
\draw (0.8,-0.5) --(0.8,-0.8);
\draw (0.7,2)--(0.7,1.2);
\draw (0.7,-2)--(0.7,-1.2);
\draw (0.3,0.8)--(0.3,0.2);
\draw (0.3,-0.8)--(0.3,-0.2);
\draw (0.7,1.6) node[anchor=east,black,scale=0.7]  {$(a,b)$};
\draw (0.3,0.5) node[anchor=east,black,scale=0.7] {$(a-1,b)$};
\draw (0.7,-1.6) node[anchor=east,black,scale=0.7] {$(a,b)$};
\draw (0.3,-0.5) node[anchor=east,black,scale=0.7] {$(a-1,b)$};
\draw[darkred,thick](-0.3,-0.2) rectangle (1.5,0.2)node[pos=0.5,scale=0.7,black] {$(a+1,b)$};
\draw[darkred,thick](-0.3,0.8) rectangle (1.5,1.2)node[pos=0.5,scale=0.7,black] {$(a,b)$};
\draw[darkred,thick](-0.3,-1.2) rectangle (1.5,-0.8)node[pos=0.5,scale=0.7,black] {$(a,b)$};
\end{tric}  
}

\def\Rp{    
\begin{tric}
\draw [double,thin](1.3,0.5) ..controls (1.8,0.5) ..(1.8,2);
\draw (1.3,0.5) --(1.3,0.2);
\draw (1.3,0.5) -- (0.8,0.5);
\draw (0.8,0.5) --(0.8,0.2);
\draw (1.05,0.5)--(1.05,0.2);
\draw [double,thin](0.8,0.5) --(0.8,0.8);
\draw [double,thin](1.3,-0.5)..controls (1.8,-0.5)..(1.8,-2);
\draw (1.3,-0.5) --(1.3,-0.2);
\draw (1.3,-0.5) -- (0.8,-0.5);
\draw (0.8,-0.5) --(0.8,-0.2);
\draw (1.05,-0.5)--(1.05,-0.2);
\draw [double,thin](0.8,-0.5) --(0.8,-0.8);
\draw (0.7,2)--(0.7,1.2);
\draw (0.7,-2)--(0.7,-1.2);
\draw (0.3,0.8)--(0.3,0.2);
\draw (0.3,-0.8)--(0.3,-0.2);
\draw (0.7,1.6) node[anchor=east,black,scale=0.7]  {$(a,b)$};
\draw (0.3,0.5) node[anchor=east,black,scale=0.7] {$(a,b-1)$};
\draw (0.7,-1.6) node[anchor=east,black,scale=0.7] {$(a,b)$};
\draw (0.3,-0.5) node[anchor=east,black,scale=0.7] {$(a,b-1)$};
\draw[darkred,thick](-0.3,-0.2) rectangle (1.5,0.2)node[pos=0.5,scale=0.7,black] {$(a+3,b-1)$};
\draw[darkred,thick](-0.3,0.8) rectangle (1.5,1.2)node[pos=0.5,scale=0.7,black] {$(a,b)$};
\draw[darkred,thick](-0.3,-1.2) rectangle (1.5,-0.8)node[pos=0.5,scale=0.7,black] {$(a,b)$};
\end{tric}
}

\def\Rq{    
\begin{tric}
\draw [double,thin](13.3,0.2)..controls(13.3,0.5)and(13.8,0.1)..(13.8,2);
\draw [double,thin](13.3,-0.2)..controls(13.3,-0.5)and(13.8,-0.1)..(13.8,-2);
\draw (12.75,2)--(12.75,1.2);
\draw (12.75,-2)--(12.75,-1.2);
\draw (12.5,0.8)--(12.5,0.2);
\draw (12.5,-0.8)--(12.5,-0.2);
\draw (12.75,1.6) node[anchor=east,black,scale=0.7] {$(a,b)$};
\draw (12.5,0.5) node[anchor=east,black,scale=0.7] {$(a,b)$};
\draw (12.75,-1.6) node[anchor=east,black,scale=0.7] {$(a,b)$};
\draw (12.5,-0.5) node[anchor=east,black,scale=0.7] {$(a,b)$};
\draw[darkred,thick](11.7,-0.2) rectangle (13.5,0.2) node[pos=0.5,scale=0.7,black] {$(a,b+1)$};
\draw[darkred,thick](11.7,0.8) rectangle (13.5,1.2) node[pos=0.5,scale=0.7,black] {$(a,b)$};
\draw[darkred, thick](11.7,-1.2) rectangle (13.5,-0.8) node[pos=0.5,scale=0.7,black] {$(a,b)$};
\end{tric}
}   
 
\begin{align}
&\Ka =\frac{1}{t_{(a,b), \varpi_1}^{(1,0)}} \Kh 
    + \frac{1}{t_{(a,b), \varpi_1}^{(-1,1)}} \Kg
    +\frac{1}{t_{(a,b), \varpi_1}^{(2,-1)}} \Kf
    +\frac{1}{t_{(a,b), \varpi_1}^{(0,0)}} \Ke \notag \\
    &+\frac{1}{t_{(a,b),\varpi_1}^{(-2,1)}} \Kd
    +\frac{1}{t_{(a,b), \varpi_1}^{(1,-1)}} \Kc
    +\frac{1}{t_{(a,b), \varpi_1}^{(-1,0)}} \Kb \label{FirstTripleClaspsExpansionGraph}\\ 
     &\Ra=\frac{1}{t_{(a,b), \varpi_2}^{(0,1)}}\Rq +     
       \frac{1}{t_{(a,b), \varpi_2}^{(3,-1)}}\Rp  + 
      \frac{1}{t_{(a,b), \varpi_2}^{(1,0)}}\Ro + 
      \frac{1}{t_{(a,b), \varpi_2}^{(-1,1)}}\Rn \notag \\     
      &+ \frac{1}{t_{(a,b), \varpi_2}^{(2,-1)}}\Rm
      + \frac{1}{t_{(a,b), \varpi_2}^{(-3,2)}}\Rl
  + \frac{ {\ }^{2,2} t_{(a,b), \varpi_2}^{(0,0)}}{D_{(a,b)}}\Rh 
  - \frac{ {\ }^{1,2} t_{(a,b), \varpi_2}^{(0,0)}}{D_{(a,b)}} \Ri \notag\\
 &-\frac{ {\ }^{2,1} t_{(a,b), \varpi_2}^{(0,0)}}{D_{(a,b)}}\Rj 
+ \frac{ {\ }^{1,1} t_{(a,b), \varpi_2}^{(0,0)}}{D_{(a,b)}}\Rk
     + \frac{1}{t_{(a,b), \varpi_2}^{(3,-2)}}\Rg  
     + \frac{1}{t_{(a,b), \varpi_2}^{(-2,1)}}\Rf \notag\\  
     &+ \frac{1}{t_{(a,b), \varpi_2}^{(1,-1)}}\Re
     + \frac{1}{t_{(a,b), \varpi_2}^{(-1,0)}}\Rd
     + \frac{1}{t_{(a,b), \varpi_2}^{(-3,1)}}\Rc
     + \frac{1}{t_{(a,b), \varpi_2}^{(0,-1)}}\Rb \label{SecondTripleClaspsExpansionGraph}
 \end{align}

\ \\
\begin{align}
 t^{(1,0)}_{(a,b), \varpi_1}&=1 \label{FirstExplicit} \\
   t^{(-1,1)}_{(a,b),\varpi_1}&=-\frac{[a+1]}{[a]} \label{texampleb}\\
   t^{(2,-1)}_{(a,b), \varpi_1}&=\frac{[3b+3][a+3b+4]}{[3b][a+3b+3]} \label{texampled}\\
   t^{(0,0)}_{(a,b), \varpi_1}&=-\frac{[a+2][a+3b+5][2a+3b+6]}{[2][a][a+3b+3][2a+3b+4]} \\
   t^{(-2,1)}_{(a,b), \varpi_1}&=\frac{[a+1][2a+3b+5][3a+3b+6]}{[a-1][2a+3b+4][3a+3b+3]}\\ 
   t^{(1,-1)}_{(a,b), \varpi_1}&=-\frac{[3b+3][a+3b+4][2a+3b+5][3a+6b+9]}{[3b][a+3b+2][2a+3b+4][3a+6b+6]}\\ 
   t^{(-1,0)}_{(a,b), \varpi_1}&=\frac{[a+1][a+3b+4][2a+3b+5][3a+3b+6][3a+6b+9]}{[a][a+3b+3][2a+3b+3][3a+3b+3][3a+6b+6]} \\ \notag\\
   t^{(0,1)}_{(a,b), \varpi_2}&=1\\
   t^{(3,-1)}_{(a,b), \varpi_2}&=-\frac{[3b+3]}{[3b]}\label{texamplec}\\ 
   t^{(1,0)}_{(a,b), \varpi_2}&=\frac{[a+3][a+3b+6]}{[3][a][a+3b+3]} \label{texamplea}\\
  t^{(-1,1)}_{(a,b), \varpi_2}&=-\frac{[a+1][a+2][2a+3b+7]}{[3][a-1][a][2a+3b+4]} \\
   t^{(2,-1)}_{(a,b), \varpi_2}&=\frac{[3b+3][a+3b+4][a+3b+5][2a+3b+7]}{[3][3b][a+3b+2][a+3b+3][2a+3b+4]} \\
   t^{(-3,2)}_{(a,b), \varpi_2}&=\frac{[a+1][3a+3b+6]}{[a-2][3a+3b+3]} \\ \notag\\
   t^{(3,-2)}_{(a,b), \varpi_2}&=\frac{[3b+3][a+3b+4][3a+6b+9]}{[3b-3][a+3b+1][3a+6b+6]} \label{tcoeffexample}\\
   t^{(-2,1)}_{(a,b), \varpi_2}&=\frac{[a+1][a+3b+6][2a+3b+5][2a+3b+6][3a+3b+6]}{[3][a-1][a+3b+3][2a+3b+3][2a+3b+4][3a+3b+3]}\\ 
   t^{(1,-1)}_{(a,b), \varpi_2}&=-\frac{[a+3][3b+3][a+3b+4][2a+3b+5][2a+3b+6][3a+6b+9]}{[3][a][3b][a+3b+2][2a+3b+3][2a+3b+4][3a+6b+6]}\\ 
    t^{(-1,0)}_{(a,b), \varpi_2}&=\frac{[a+1][a+2][a+3b+4][a+3b+5][2a+3b+5][3a+3b+6][3a+6b+9]}{[3][a-1][a][a+3b+2][a+3b+3][2a+3b+3][3a+3b+3][3a+6b+6]}\\ 
     t^{(-3,1)}_{(a,b), \varpi_2}&=-\frac{[a+1][2a+3b+5][3a+3b+6][3a+6b+9]}{[a-2][2a+3b+2][3a+3b][3a+6b+6]} \\
     t^{(0,-1)}_{(a,b), \varpi_2}&=\frac{[3b+3][a+3b+4][2a+3b+5][3a+3b+6][3a+6b+9]}{[3b][a+3b+1][2a+3b+2][3a+3b+3][3a+6b+3]} \label{LastRound}\\ \notag\\
\left( {\ }^{p,\ell} t_{(a,b),\varpi_2}^{(0,0)} \right)&=
\begin{pmatrix}
  {\ }^{1,1} t_{(a,b), \varpi_2}^{(0,0)}
 &  {\ }^{1,2} t_{(a,b), \varpi_2}^{(0,0)}  \\
  {\ }^{2,1} t_{(a,b), \varpi_2}^{(0,0)} 
 & {\ }^{2,2} t_{(a,b), \varpi_2}^{(0,0)}   
\end{pmatrix}\\
 \mathcal{D}_{(a,b)}:&=
\text{det} \left( {\ }^{p,\ell} t_{(a,b), \varpi_2}^{(0,0)}  \right)
\end{align}

\newpage

Moreover,
\begin{align}
 &{\ }^{1,2} t_{(a,b), \varpi_2}^{(0,0)}=  
 {\ }^{2,1} t_{(a,b), \varpi_2}^{(0,0)}, \\
&\mathcal{D}_{(a,b)}= \dfrac{[4][6][a+2][3b+6][a+3b+5][2a+3b+6][3a+3b+9][3a+6b+12]}
{[2][3][12][a][3b][a+3b+3][2a+3b+4][3a+3b+3][3a+6b+6]}, \label{DeterminantR00} 
\end{align}
and the entries of the matrix can be computed from the relations in Appendix \ref{Relations}.

\subsection{Verifying the clasp conjecture}
\label{subsec-proof-claspconj}

Before proving our main theorem, we will prove the following, which implies the clasp conjecture in type $G_2$.  

\begin{corollary}\label{cor-claspconj}
Fix $\lambda = a\varpi_1+ b\varpi_2$  with $a,b\in \mathbb{Z}_{\ge 0}$. Let $\varpi\in \{\varpi_1, \varpi_2\}$ be a fundamental weight, and let $\mu \in W\cdot \varpi$ be a weight in the Weyl group orbit of $\varpi$. Then 
\[
t_{\lambda, \varpi}^{\mu} = \pm\prod_{\alpha\in {\Phi_{\mu}}}\frac{[(\alpha^{\vee}, \lambda + \rho)]_{\ell^{\ell(\alpha)}}}{[(\alpha^{\vee}, \lambda +\varpi + \rho)]_{\ell^{\ell(\alpha)}}}.
\]
\end{corollary}

\begin{proof}
In type $G_2$, The $W$ invariant bilinear pairing on $\mathbb{Z}\Phi$ is determined by $(\alpha_1, \alpha_1) = 2$ and $(\alpha_2, \alpha_2) = 6$. In particular, $l(\alpha_1) = 1$ and $l(\alpha_2) = 3$. We set $\alpha^{\vee} = 2\alpha/(\alpha, \alpha)$. The positive roots are
\begin{equation}
\alpha_1, 3\alpha_1 + \alpha_2, 2\alpha_1 + \alpha_2, 3\alpha_1 + 2\alpha_2, \alpha_1+ \alpha_2, \alpha_2,
\end{equation}
the corresponding coroots are
\begin{equation}
\alpha_1^{\vee}, \alpha_1^{\vee} + \alpha_2^{\vee}, 2\alpha_1^{\vee} + 3\alpha_2^{\vee}, \alpha_1^{\vee} + 2\alpha_2^{\vee} , \alpha_1^{\vee} + 3\alpha_2^{\vee}, \alpha_2^{\vee}.
\end{equation}
To simplify notation, we will write $s_i:=s_{\alpha_i}$. It is not hard to see that
\begin{equation}
d_{(1, 0)} = 1, \ \ \ \ \ \ \ \ \ \ d_{(0, 1)} = 1, \ \ \ \ \ \ \ \ \ \ d_{(-2, 1)} = s_1s_2s_1, \ \ \ \ \ \ \ \ \ \ d_{(3, -2)} = s_2s_1s_2,  
\end{equation}
\begin{equation}
d_{(-1, 1)} = s_1, \ \ \ \ \ \ \ \ \ \ d_{(3, -1)}= s_2, \ \ \ \ \ \ \ \ \ \ d_{(1, -1)}= s_1s_2s_1s_2, \ \ \ \ \ \ \ \ \ \ d_{(-3, 1)} = s_2s_1s_2s_1,
\end{equation}
\begin{equation}
d_{(2, -1)} = s_1s_2, \ \ \ \ \ \ \ \ \ \ d_{(-3, 2)} = s_2s_1, \ \ \ \ \ \ \ \ \ \ d_{(-1, 0)} = s_1s_2s_1s_2s_1, \ \ \ \ \ \text{and} \ \ \ \ \ d_{(0, -1)} = s_2s_1s_2s_1s_2.
\end{equation}

The claim then follows from the formulas for $t_{\lambda, \varpi}^{\mu}$ in Section \ref{results}. One verifies this by using that if $w= s_{\beta_1}s_{\beta_2}\ldots s_{\beta_n}$, then 
\[
\lbrace \alpha\in \Phi_+ \ : \ w\alpha\in \Phi_-\rbrace = \lbrace \beta_n, s_{\beta_n}(\beta_{n-1}), s_{\beta_n}s_{\beta_{n-1}}(\beta_{n-2}), \ldots s_{\beta_n}s_{\beta_{n-1}}\ldots s_{\beta_2}(\beta_1)\rbrace,
\]
along with the quantum number identity $[n]_{q^{3}} = [3n]/[3]$.
\end{proof}

\subsection{Proof of Triple Clasp Formula}
\label{subsec-proof}

Suppose that $V(\wt\unu)$ is a summand of $V(\wt\unw)\ot V(\varpi)$, and that $ (a,b)=\lambda = \wt \unw $ and $(m, n)=\mu=\wt \unu- \wt \unw$. Then we will write 
\[
{^{p\ell}t}_{\unw, \varpi}^{\unu}:= {^{p\ell}t}_{(a,b), \varpi}^{(m,n)}.
\]
Our convention is that the $p\ell$ superscript is neglected when $\dim V(\varpi)_{\mu}=1$. By definition the elements ${^{p\ell}t}_{\unw, \varpi}^{\unu}$ only depend on the weights $\wt \unw$ and $\wt \unu$, not the words $\unw$ and $\unu$.

We will write $(^{p\ell}t_{\unw, \varpi}^{\unu})$ to denote the matrix of scalars $^{p\ell}t_{\unw, \varpi}^{\unu}$, for $v_p, v_\ell \in \text{Ker}_{\wt \unw}(V(\varpi)_{\mu})$.

\begin{lemma}\label{tsinvertible}
The matrix $(^{p\ell}t_{\unw, \varpi}^{\unu})$ is invertible over $\mathbb{C}(q)$.  
\end{lemma}
\begin{proof}
From Section \ref{results}, using Equation \eqref{FirstExplicit} to Equation \eqref{LastRound}, and Equation \eqref{DeterminantR00}, one can check that the determinant of this matrix is invertible in $\mathbb{C}(q)$.
\end{proof}

\begin{defn}
Let $\unw\in \DD$, write $\lambda= \wt \unw$, and let $\varpi \in \lbrace \varpi_1, \varpi_2\rbrace$. We define the \textbf{triple clasp} to be the following inductively defined diagram:
\[
T_{\unw\otimes\varpi} := T_{\unw}\otimes \id_{\varpi} - \sum_{\substack{\mu \in \wt V(\varpi)\backslash \lbrace\varpi\rbrace \\ v_i, v_j\in \text{Ker}_{\wt \unw}(V(\varpi)_{\mu})}} ((^{p\ell}t_{\unw, \varpi}^{\unu_{\lambda+ \mu}})^{-1})_{ij}\cdot{^{ij}\mathbb{TT}}_{\unw, \varpi}^{\unu_{\lambda+ \mu}},
\]
where
\[
{^{ij}\mathbb{TT}}_{\unw, \varpi}^{\unu_{\lambda+ \mu}}:=( T_{\unw}\ot \id_{\varpi} ) \circ (\mathbb{D}(^iELL_{\unw, \varpi}^{\unu_{\lambda+ \mu}}))\circ T_{\unu_{\lambda+ \mu}} \circ (^jELL_{\unw, \varpi}^{\unu_{\lambda + \mu}})\circ ( T_{\unw}\ot \id_{\varpi} ).
\]
\end{defn}

\begin{defn}
By Lemma \ref{claspschur} there is a scalar ${^{p\ell}\kappa_{\unw, \varpi}^{\unu_{\lambda+ \mu}}}$ such that
\[
{^pLL_{\unw, \varpi}^{\unu_{\lambda+ \mu}}}\circ \mathbb{D}({^\ell LL_{\unw,\varpi}^{\unu_{\lambda+ \mu}}}) = {^{p\ell}\kappa_{\unw, \varpi}^{\unu_{\lambda+ \mu}}}\cdot C_{\unu_{\lambda+ \mu}}.
\]
We call the matrix 
\[
(^{p\ell}\kappa_{\unw, \varpi}^{\unu_{\lambda + \mu}}),
\]
such that $v_p, v_\ell \in \text{Ker}_{\wt \unw}(V(\varpi)_{\mu})$, a \textbf{local intersection form matrix}.
\end{defn}

\begin{lemma}
If $\wt\unw= \wt \unw'$, then
\[
^{p\ell}\kappa_{\unw, \varpi}^{\unu_{\lambda + \mu}} = {^{p\ell}\kappa_{\unw', \varpi}^{\unu_{\lambda + \mu}}}.
\]
\end{lemma}
\begin{proof}
First observe that
\begin{align*}
^{p\ell}\kappa_{\unw, \varpi}^{\unu_{\lambda + \mu}}\cdot C_{\unu_{\lambda+ \mu}} &= {^pLL_{\unw, \varpi}^{\unu_{\lambda+ \mu}}}\circ \mathbb{D}({^\ell LL_{\unw,\varpi}^{\unu_{\lambda+ \mu}}}) \\
&= {^pLL_{\unw, \varpi}^{\unu_{\lambda+ \mu}}}\circ (C_{\unw}\ot \id_{\varpi}) \circ \mathbb{D}({^\ell LL_{\unw,\varpi}^{\unu_{\lambda+ \mu}}}) \\
&= {^pLL_{\unw, \varpi}^{\unu_{\lambda+ \mu}}}\circ (\mathsf{H}_{\unw'}^{\unw}\ot \id_{\varpi})\circ (C_{\unw'}\ot \id_{\varpi})\circ (\mathsf{H}_{\unw}^{\unw'}\ot \id_{\varpi})\circ \mathbb{D}({^\ell LL_{\unw,\varpi}^{\unu_{\lambda+ \mu}}}) \\
&= {^pLL_{\unw', \varpi}^{\unu_{\lambda+ \mu}}}\circ  (C_{\unw'}\ot \id_{\varpi})\circ \mathbb{D}({^\ell LL_{\unw',\varpi}^{\unu_{\lambda+ \mu}}}) \\
&= {^pLL_{\unw', \varpi}^{\unu_{\lambda+ \mu}}}\circ \mathbb{D}({^\ell LL_{\unw',\varpi}^{\unu_{\lambda+ \mu}}}) \\
&={^{p\ell}\kappa_{\unw', \varpi}^{\unu_{\lambda + \mu}}\cdot C_{\unu_{\lambda+ \mu}}}.
\end{align*}
The claim follows from comparing neutral coefficients. 
\end{proof}

\begin{notation}
We will write 
\[
^{p\ell}\kappa_{\lambda, \varpi}^{\mu}:={^{p\ell}\kappa_{\unw, \varpi}^{\unu_{\lambda+ \mu}},}
\]
where $\wt\unw= (a,b)$.
\end{notation}

\def\kexamplea{
\begin{tric}
\draw (0.3,-0.3)--(0.3,-2.2);
\draw [double,thin]  (0.7,-0.3)--(0.7,-2.2) (1.1,-0.3)--(1.1,-2.2);

\draw (0.3,0.3)--(0.3,2.2);
\draw [double,thin]  (0.7,0.3)--(0.7,2.2) (1.1,0.3)--(1.1,2.2);

\draw [darkred,thick] (-0.1,-0.3) rectangle (1.5,0.3) node[midway,black]{$C_{\unu_{\lambda+\mu}}$};
\end{tric}
}

\def\kexampleb{
\begin{tric}
\draw [double,thin](1.2,0.8)..controls (1.8,0.8)..(1.8,0) ..controls(1.8,-0.8)..(1.2,-0.8);
\draw (1.2,0.8) --(1.2,0.3);
\draw (1.2,0.8) -- (0.8,0.8);
\draw [double,thin](0.8,0.8) --(0.8,1.3);
\draw (0.8,0.8) --(0.8,0.3);
\draw [double,thin] (0.5,0.3)--(0.5,1.3);
\draw (0.2,0.3)--(0.2,1.3);

\draw (0.2,-1.9)--(0.2,-2.2)  ;
\draw [double,thin] (0.5,-1.9)--(0.5,-2.2) (0.8,-1.9)--(0.8,-2.2);

\draw (1.2,-0.8) --(1.2,-0.3);
\draw (1.2,-0.8) -- (0.8,-0.8);
\draw [double,thin](0.8,-0.8) --(0.8,-1.3);
\draw (0.8,-0.8) --(0.8,-0.3);
\draw [double,thin] (0.5,-0.3)--(0.5,-1.3);
\draw (0.2,-0.3)--(0.2,-1.3);

\draw (0.2,1.9)--(0.2,2.2) ;
\draw [double,thin] (0.5,1.9)--(0.5,2.2)  (0.8,1.9)--(0.8,2.2);

\draw [darkred,thick] (-0.1,-1.3) rectangle (1.5,-1.9) node[midway,black]{$C_{\unu_{\lambda+\mu}}$};
\draw [darkred,thick] (-0.1,1.3) rectangle (1.5,1.9) node[midway,black]{$C_{\unu_{\lambda+\mu}}$};
\draw [darkred,thick] (-0.1,-0.3) rectangle (1.5,0.3) node[midway,black]{$C_{\unw}$};
\end{tric}
}

\def\kexamplec{
\begin{tric}
\draw [double,thin](1.2,-0.8)..controls (1.8,-0.8)..(1.8,0)
..controls (1.8,0.8)..(1.2,0.8);
\draw (1.2,0.8) --(1.2,0.3);
\draw (1.2,0.8) -- (0.8,0.8);
\draw [double,thin](0.8,0.8) --(0.8,1.3);
\draw (0.8,0.8) --(0.8,0.3);
\draw (0.4,0.3)--(0.4,1.3)node[left,midway,scale=0.7]{$(1,1)$};

\draw (0.7,-1.9)--(0.7,-2.4) node[left,midway,scale=0.7]{$(1,2)$};

\draw (1.2,-0.8) --(1.2,-0.3);
\draw (1.2,-0.8) -- (0.8,-0.8);
\draw [double,thin](0.8,-0.8) --(0.8,-1.3);
\draw (0.8,-0.8) --(0.8,-0.3);
\draw (0.4,-0.3)--(0.4,-1.3)node[left,midway,scale=0.7]{$(1,1)$};

\draw (0.7,1.9)--(0.7,2.4) node[left,midway,scale=0.7]{$(1,2)$};

\draw [darkred,thick] (-0.1,-1.3) rectangle (1.5,-1.9) node[midway,black]{$(1,2)$};
\draw [darkred,thick] (-0.1,1.3) rectangle (1.5,1.9) node[midway,black]{$(1,2)$};
\draw [darkred,thick] (-0.1,-0.3) rectangle (1.5,0.3) node[midway,black]{$(3,1)$};
\end{tric}
}

\def\kexampled{
\begin{tric}
\draw (0.1,-0.3)--(0.1,-3.2) (0.4,-0.3)--(0.4,-3.2) (0.7,-0.3)--(0.7,-3.2) ;
\draw[double,thin] (1,-0.3)--(1,-3.2);

\draw (0.1,0.3)--(0.1,3.2) (0.4,0.3)--(0.4,3.2) (0.7,0.3)--(0.7,3.2) ;
\draw[double,thin] (1,0.3)--(1,3.2);

\draw [darkred,thick] (-0.3,-0.3) rectangle (1.4,0.3) node[midway,black]{$C_{\unu_{\lambda+\mu}}$};
\end{tric}
}

\def\kexamplee{
\begin{tric}
\draw[double,thin] (1.2,-1.3) ..controls (1.8,-1.3) ..(1.8,0) ..controls (1.8,1.3) ..(1,1.3);
\draw [double,thin](1,-2.3)--(1,-1.45);
\draw [double,thin](1,-0.9)--(1,-1.15);
\draw (1,-1.45)--(1,-1.15)--(1.2,-1.3)--cycle;
\draw (1,-0.3)--(1,-0.9)--(0.7,-0.9)--(0.7,-2.3);
\draw (0.7,-0.3)--(0.7,-0.6)--(0.4,-0.6)--(0.4,-2.3);
\draw [double,thin] (0.4,-0.6)--(0.4,-0.3) (0.7,-0.6)--(0.7,-0.9);
\draw (0.1,-0.3)--(0.1,-2.3);

\draw (0.1,-2.9)--(0.1,-3.2) (0.4,-2.9)--(0.4,-3.2) (0.7,-2.9)--(0.7,-3.2) ;
\draw[double,thin] (1,-2.9)--(1,-3.2) ;

\draw  (0.7,2.3)--(0.7,1.7)--(1,1.7)--(1,0.3);
\draw (0.7,0.3)--(0.7,0.9)--(0.4,0.9)--(0.4,2.3);
\draw[double,thin] (0.7,0.9)--(0.7,1.7) (1,1.7)--(1,2.3) (0.4,0.9)--(0.4,0.3);
\draw (0.1,0.3)--(0.1,2.3);

\draw (0.1,2.9)--(0.1,3.2) (0.4,2.9)--(0.4,3.2) (0.7,2.9)--(0.7,3.2) ;
\draw[double,thin] (1,2.9)--(1,3.2);

\draw [darkred,thick] (-0.3,-2.3) rectangle (1.4,-2.9) node[midway,black]{$C_{\unu_{\lambda+\mu}}$};
\draw [darkred,thick] (-0.3,-0.3) rectangle (1.4,0.3) node[midway,black]{$C_{\unw}$};
\draw [darkred,thick] (-0.3,2.3) rectangle (1.4,2.9) node[midway,black]{$C_{\unu_{\lambda+\mu}}$};
\end{tric}
}

\def\kexamplef{
\begin{tric}
\draw[double,thin] (1.2,-1.3) ..controls (1.8,-1.3) ..(1.8,0)..controls (1.8,1.3) ..(1,1.3);
\draw [double,thin](1,-2.3)--(1,-1.45);
\draw [double,thin](1,-0.3)--(1,-1.15);
\draw (1,-1.45)--(1,-1.15)--(1.2,-1.3)--cycle;

\draw (0.6,2.9)--(0.6,3.4) node[left,midway,scale=0.7]{$(3,1)$};
\draw (0.4,-0.3)--(0.4,-2.3) node[left,midway,scale=0.7]{$(3,0)$};

\draw (1,0.3)--(1,2.3);

\draw (0.6,-2.9)--(0.6,-3.4) node[left,midway,scale=0.7]{$(3,1)$};
\draw (0.4,0.3)--(0.4,2.3) node[left,midway,scale=0.7]{$(2,1)$};

\draw [darkred,thick] (-0.3,-2.3) rectangle (1.4,-2.9) node[midway,black]{$(3,1)$};
\draw [darkred,thick] (-0.3,-0.3) rectangle (1.4,0.3) node[midway,black]{$(3,1)$};
\draw [darkred,thick] (-0.3,2.3) rectangle (1.4,2.9) node[midway,black]{$(3,1)$};
\end{tric}
}

\def\kexampleg{
\begin{tric}
\draw (0.1,-0.3)--(0.1,-3.2) (0.4,-0.3)--(0.4,-3.2) (0.7,-0.3)--(0.7,-3.2) ;
\draw[double,thin] (1,-0.3)--(1,-3.2);

\draw (0.1,0.3)--(0.1,3.2) (0.4,0.3)--(0.4,3.2) (0.7,0.3)--(0.7,3.2) ;
\draw[double,thin] (1,0.3)--(1,3.2);

\draw [darkred,thick] (-0.3,-0.3) rectangle (1.4,0.3) node[midway,black]{$C_{\unu_{\lambda+\mu}}$};
\end{tric}
}

\def\kexampleh{
\begin{tric}
\draw[double,thin] (1.2,1.3) ..controls (1.8,1.3) ..(1.8,0) ..controls (1.8,-1.3) ..(1,-1.3);
\draw [double,thin](1,2.3)--(1,1.45);
\draw [double,thin](1,0.9)--(1,1.15);
\draw (1,1.45)--(1,1.15)--(1.2,1.3)--cycle;
\draw (1,0.3)--(1,0.9)--(0.7,0.9)--(0.7,2.3);
\draw (0.7,0.3)--(0.7,0.6)--(0.4,0.6)--(0.4,2.3);
\draw [double,thin] (0.4,0.6)--(0.4,0.3) (0.7,0.6)--(0.7,0.9);
\draw (0.1,0.3)--(0.1,2.3);

\draw (0.1,-2.9)--(0.1,-3.2) (0.4,-2.9)--(0.4,-3.2) (0.7,-2.9)--(0.7,-3.2) ;
\draw[double,thin] (1,-2.9)--(1,-3.2) ;

\draw  (0.7,-2.3)--(0.7,-1.7)--(1,-1.7)--(1,-0.3);
\draw (0.7,-0.3)--(0.7,-0.9)--(0.4,-0.9)--(0.4,-2.3);
\draw[double,thin] (0.7,-0.9)--(0.7,-1.7) (1,-1.7)--(1,-2.3) (0.4,-0.9)--(0.4,-0.3);
\draw (0.1,-0.3)--(0.1,-2.3);

\draw (0.1,2.9)--(0.1,3.2) (0.4,2.9)--(0.4,3.2) (0.7,2.9)--(0.7,3.2) ;
\draw[double,thin] (1,2.9)--(1,3.2);

\draw [darkred,thick] (-0.3,-2.3) rectangle (1.4,-2.9) node[midway,black]{$C_{\unu_{\lambda+\mu}}$};
\draw [darkred,thick] (-0.3,-0.3) rectangle (1.4,0.3) node[midway,black]{$C_{\unw}$};
\draw [darkred,thick] (-0.3,2.3) rectangle (1.4,2.9) node[midway,black]{$C_{\unu_{\lambda+\mu}}$};
\end{tric}
}

\def\kexamplei{
\begin{tric}
\draw[double,thin] (1.2,1.3) ..controls (1.8,1.3) ..(1.8,0)..controls (1.8,-1.3) ..(1,-1.3);
\draw [double,thin](1,2.3)--(1,1.45);
\draw [double,thin](1,0.3)--(1,1.15);
\draw (1,1.45)--(1,1.15)--(1.2,1.3)--cycle;
\draw (0.6,-2.9)--(0.6,-3.4) node[left,midway,scale=0.7]{$(3,1)$};
\draw (0.4,0.3)--(0.4,2.3) node[left,midway,scale=0.7]{$(3,0)$};

\draw (1,-0.3)--(1,-2.3);
\draw (0.6,2.9)--(0.6,3.4) node[left,midway,scale=0.7]{$(3,1)$};
\draw (0.4,-0.3)--(0.4,-2.3) node[left,midway,scale=0.7]{$(2,1)$};

\draw [darkred,thick] (-0.3,-2.3) rectangle (1.4,-2.9) node[midway,black]{$(3,1)$};
\draw [darkred,thick] (-0.3,-0.3) rectangle (1.4,0.3) node[midway,black]{$(3,1)$};
\draw [darkred,thick] (-0.3,2.3) rectangle (1.4,2.9) node[midway,black]{$(3,1)$};
\end{tric}
}

\begin{example}
Consider $\unw=\varpi_1\varpi_2\varpi_1\varpi_1$ and $\varpi=\varpi_2$. When $\mu= (-2,1)$, choose $\unu_{\lambda+\mu}=\varpi_1\varpi_2\varpi_2$. Then
\[
^{1,1}\kappa_{\unw, \varpi}^{\unu_{\lambda+ \mu}} \ \kexamplea
=\kexampleb= \kexamplec \quad .
\]
When $\mu= (0,0)$,  choose $\unu_{\lambda+\mu}=\varpi_1\varpi_1\varpi_1\varpi_2$. Then we have the following.
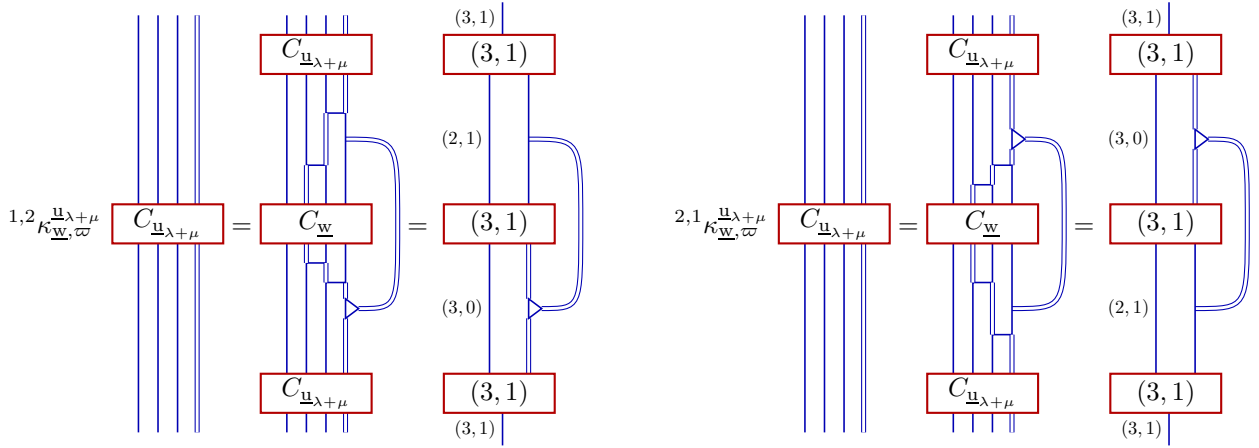
\[
^{1,2}\kappa_{\unw, \varpi}^{\unu_{\lambda+ \mu}} \ \kexampled
=\kexamplee= \kexamplef \ \ \ \ \ \ \ \ \ \ ^{2,1}\kappa_{\unw, \varpi}^{\unu_{\lambda+ \mu}} \ \kexampleg
=\kexampleh= \kexamplei
\]
\end{example}

\begin{remark}
When $\mu=(0,0)$, by taking the quantum trace, we know that $^{1,2}\kappa_{\unw, \varpi}^{\unu_{\lambda+ \mu}} 
= {\ }^{2,1}\kappa_{\unw, \varpi}^{\unu_{\lambda+ \mu}}$. So the local intersection form matrix  $\left(^{p\ell}\kappa_{\unw, \varpi}^{\unu_{\lambda+ \mu}}\right)$ is symmetric.\\
\end{remark}

\begin{notation}
Fix the following set of formal variables $\mathcal{X}:= \big\{ x_{(a,b), \varpi}^{(c,d)} \ | \ a,b,c,d\in \mathbb{Z}_{\ge 0} \ \text{and} \ \varpi\in \{\varpi_1, \varpi_2\}\big\}$. We will consider elements in the ring $\mathcal{A}:=\mathbb{C}(q)[x^{\pm 1} \ | \ x\in \mathcal{X}]$. 

Suppose that $ (a,b)=\lambda = \wt \unw $ and $(m, n)=\mu=\wt \unu_{\lambda+ \mu}- \wt \unw$,  then we will write \[
{^{p\ell}\rho}_{\unw, \varpi}^{\unu_{\lambda+ \mu}}:= {^{p\ell}\rho}_{(a,b), \varpi}^{(m, n)}\in \mathcal{A},
\]
where ${^{p\ell}\rho}_{(a,b), \varpi}^{(m, n)}$ is the recursive relation described by Equation \eqref{Recursion(1)} to Equation \eqref{Recursion(22)}, in Appendix \ref{Recursions}. 

We also write ${^{p\ell}\rho}_{(a,b), \varpi}^{(m, n)}(\kappa)$ to denote the right hand side of the recursive relation with each $^{p\ell}x_{\lambda, \varpi}^{\mu}$ replaced with $^{p\ell}\kappa_{\lambda, \varpi}^{\mu}$. Similarly, we write ${^{p\ell}\rho}_{(a,b), \varpi}^{(m, n)}(t)$ to denote the right hand side of the recursive relation with each $^{p\ell}x_{\lambda, \varpi}^{\mu}$ replaced by $^{p\ell}t_{\lambda, \varpi}^{\mu}$. 

Our convention is that the $p\ell$ superscript is neglected when $\dim V(\varpi)_{\mu}=1$.
\end{notation}

\begin{example}
Consider $\unw$ and $\unu_{\lambda+ \mu}$ in $\DD$ such that $\wt \unw=(a,b)$ and $\wt \unu_{\lambda+\mu}=(a+1,b)$. Also, let $\varpi= \varpi_2$. By Equation \eqref{Recursion(10)}:
\begin{align*}
  {\rho}_{\unw, \varpi_2}^{\unu_{\lambda+ \mu}}= {\rho}_{(a,b), \varpi_2}^{(1,0)}= \frac{[7]}{[3]}
-\frac{1}{x_{(a-1,b), \varpi_1}^{(-1,1)}} x_{(a-2,b+1), \varpi_2}^{(3,-1)}
          -\frac{1}{x_{(a-1,b), \varpi_1}^{(2,-1)}}.
\end{align*}
\end{example}

\begin{remark}
Since the recursive relations in Appendix \ref{Recursions} are elements of $\mathcal{A}$, there is no question whether a particular element in $\mathcal{X}$ appearing in a relation is invertible or not. Thanks to Lemma \ref{tsinvertible} the elements ${^{p\ell}\rho}_{(a,b), \varpi}^{(m,n)}(t)$ are also always well defined. This is not obviously true for ${^{p\ell}\rho}_{(a,b), \varpi}^{(m,n)}(\kappa)$. However we will prove that $\kappa_{\unw, \varpi}^{\unu_{\lambda+ \mu}} = t_{\unw, \varpi}^{\unu_{\lambda+ \mu}}$.
\end{remark}

\begin{thm}\label{mainthmtoprove}
If $\unw\in \DD$, $\lambda= \wt \unw$, and $\varpi\in \lbrace 1, 2\rbrace$, then
\[
T_{\unw} = C_{\unw} \quad \text{and} \quad \kappa_{\unw, \varpi}^{\unu_{\lambda+ \mu}} = t_{\unw, \varpi}^{\unu_{\lambda+ \mu}}. 
\]
\end{thm}

We will prove Theorem \ref{mainthmtoprove} by induction. To simplify the arguments, we will break the various steps of the proof into smaller lemmas about the following statements. In what follows we write $\lambda=\wt\unw$.
\begin{align*}
S_1(\unw):&= \bigg(T_{\unw} = C_{\unw}\bigg) \\
S_1'(\unw, \varpi):&= \bigg({^{p\ell}\kappa}_{\unw, \varpi}^{\unu_{\lambda+ \mu}} = {^{p\ell}t}_{\unw, \varpi}^{\unu_{\lambda+ \mu}},\text{for all} \ \mu\in \wt V(\varpi) \ \text{and} \  \text{for all} \ v_{\mu, p}, v_{\mu, \ell}\in \text{Ker}_{\lambda}(V(\varpi)_{\mu})\bigg) \\
S_2(\unw, \varpi):&= \bigg(C_{\unw}\otimes \id_{\varpi} \in \text{span} \ \bigcup_{\substack{\mu\in \wt V(\varpi) \\ v_p, v_\ell\in \text{Ker}_{\lambda}(V(\varpi)_{\mu})}}\lbrace {^{p\ell}\mathbb{LL}}_{\unw, \varpi}^{\unu_{\lambda+ \mu}}\rbrace \bigg) \\
S_2'(\unw, \varpi) :&= \bigg(\lbrace {^{p}LL}_{\unw, \varpi}^{\unu_{\lambda+ \mu}}\rbrace_{v_p\in \text{Ker}_{\lambda}(V(\varpi)_{\mu})} \ \text{is a linearly independent set, for all} \ \mu\in \wt V(\varpi) \bigg) \\
S_3(\unw, \varpi):&= \bigg(\lbrace {^{p}LL}_{\unw, \varpi}^{\unu_{\lambda+ \mu}}\rbrace_{v_p\in \text{Ker}_{\lambda}(V(\varpi)_{\mu})} \ \text{is a basis for} \ \Hom_{\Kar\DD}(C_{\unw}\otimes \id_{\varpi}, C_{\unu_{\lambda+ \mu}}), \text{for all} \ \mu\in \wt V(\varpi) \bigg) \\
S_4(\unw, \varpi):&= \bigg({^{p\ell}\kappa}_{\unw, \varpi}^{\unu_{\lambda+ \mu}} = {^{p\ell}\rho}_{\unw, \varpi}^{\unu_{\lambda+ \mu}}(\kappa), \text{for all} \ \mu\in \wt V(\varpi)\bigg) \\
S_5(\unw, \varpi):&= \bigg(T_{\unw\otimes \varpi}= C_{\unw}\otimes \id_{\varpi} - \sum_{\substack{\mu\in \wt V(\varpi)\backslash\lbrace{\varpi}\rbrace \\ v_i, v_j\in \text{Ker}_{\lambda}(V(\varpi)_{\mu})}} ({^{p\ell}\kappa}_{\unw, \varpi}^{\unu_{\lambda+ \mu}})^{-1}_{ij}\cdot {^{ij}\mathbb{LL}}_{\unw, \varpi}^{\unu_{\lambda+ \mu}}\bigg) \\
S_6(\unw):&= \bigg(T_{\unw}^2= T_{\unw}\bigg) \\
S_6'(\unw, \varpi):&= \bigg(C_{\unu}\circ D\circ T_{\unw\otimes \varpi} = 0, \text{for all $\unu$ such that $(\wt \unu- \lambda)\in \wt V(\varpi)\backslash\lbrace{\varpi}\rbrace$ for all possible diagrams $D$}\bigg)
\end{align*}

\begin{lemma}
If $V(\lambda+ \mu)$ is a summand of $V(\lambda)\ot V(\varpi)$, then $\lambda+ \mu\le \lambda + \varpi$.
\end{lemma}
\begin{proof}
Suppose $V(\lambda+ \mu)$ is a summand of $V(\lambda)\ot V(\varpi)$. Then $V(\varpi)_{\mu}\ne 0$. It follows that $\mu\in \varpi + \mathbb{Z}_{\le 0}\Phi_+$, so $\varpi - \mu \ge 0$. 
\end{proof}

\begin{lemma}\label{lem:1imp2}
If $S_1(\unx)$ for all $\unx$ such that $\wt \unx \le \wt (\unw\otimes \varpi)$, then $S_2(\unw, \varpi)$.
\end{lemma}
\begin{proof}
Write $\lambda= \wt\unw$. Since $T_{\unw\otimes \varpi} = C_{\unw\otimes\varpi}$,  
\[
C_{\unw\otimes\varpi}= T_{\unw\otimes\varpi}= T_{\unw}\otimes \id_{\varpi} - \sum_{\substack{\mu \in \wt V(\varpi)\backslash\lbrace{\varpi}\rbrace \\ v_i, v_j\in \text{Ker}_{\lambda}(V(\varpi)_{\mu})}} (^{p\ell}t_{\unw, \varpi}^{\unu_{\lambda+ \mu}})^{-1}_{ij}\cdot{^{ij}\mathbb{TT}}_{\unw, \varpi}^{\unu_{\lambda+ \mu}}.
\]
Also, we have $\wt \unw< \wt(\unw\ot\varpi)$ and $\wt\unu_{\lambda+ \mu} \le \wt(\unw\otimes\varpi)$, for all $\mu$ such that $V(\lambda+ \mu)$ is a summand of $V(\lambda)\otimes V(\varpi)$, so $T_{\unw}= C_{\unw}$ and $T_{\unu_{\lambda+ \mu}} = C_{\unu_{\lambda+ \mu}}$. Therefore, 
\begin{align*}
{^{ij}\mathbb{TT}}_{\unw, \varpi}^{\unu_{\lambda+\mu}}&:= (T_{\unw}\ot \id_{\varpi} )\circ (\mathbb{D}(^iELL_{\unw, \varpi}^{\unu_{\lambda+\mu}}))\circ T_{\unu_{\lambda+\mu}} \circ (^jELL_{\unw, \varpi}^{\unu_{\lambda+\mu}})\circ (T_{\unw}\ot \id_{\varpi}) \\
&=( C_{\unw}\ot \id_{\varpi} )\circ (\mathbb{D}(^iELL_{\unw, \varpi}^{\unu_{\lambda+\mu}}))\circ C_{\unu_{\lambda+\mu}} \circ (^jELL_{\unw, \varpi}^{\unu_{\lambda+\mu}})\circ (C_{\unw}\ot \id_{\varpi}) \\
&= {^{ij}\mathbb{LL}}_{\unw, \varpi}^{\unu_{\lambda+\mu}}.
\end{align*}
The claim follows from observing that $C_{\unw\otimes \varpi} = \mathbb{LL}_{\unw,\varpi}^{\unu_{\lambda+ \varpi}}$, which is a consequence of Lemmas \ref{claspabsorption} and \ref{neturalabsorption} (i.e. clasp absorption and neutral absorption). \\
\end{proof}

\begin{lemma}\label{lem:1'imp2'}
If $S_1'(\unw, \varpi)$, then $S_2'(\unw, \varpi)$.  \end{lemma}

\begin{proof}
Write $\lambda= \wt \unw$. For each $\mu$ such that $V(\lambda+ \mu)$ is a summand of $V(\lambda)\otimes V(\varpi)$, consider the linear relation
\[
\sum_{p} \xi_p\cdot {^{p}LL}_{\unw, \varpi}^{\unu_{\lambda+ \mu}} = 0.
\]
We can precompose the relation with $\mathbb{D}({^{\ell}LL}_{\unw, \varpi}^{\unu_{\lambda+ \mu}})$ for all $v_\ell\in \text{Ker}_{\lambda}(V(\varpi)_{\mu})$ to obtain a family of relations
\[
\sum_{p} \xi_p \cdot {^{p\ell}\kappa}_{\unw, \varpi}^{\unu_{\lambda+ \mu}} = 0.
\]
By our hypothesis, we obtain
\[
\sum_{p} \xi_p \cdot {^{p\ell}t}_{\unw, \varpi}^{\unu_{\lambda+ \mu}} = 0,
\]
and it follows from Lemma \ref{tsinvertible} that each $\xi_p= 0$.
\end{proof}

\begin{lemma}\label{lem:2'imp3}
If $S_2'(\unw, \varpi)$, then $S_3(\unw, \varpi)$. \end{lemma}

\begin{proof}
Let $\mu = \wt\unu- \wt\unw$. By combining Corollary \ref{karoubi} with Equation \eqref{PRV}, we may deduce the following
\[
\dim \Hom_{\Kar\DD}(C_{\unw}\otimes \id_{\varpi}, C_{\unu})= \dim \Hom_{\Uq}(V(\wt \unw)\otimes V(\varpi), V(\wt\unu)) = \dim \text{Ker}_{\wt \unw}(V(\varpi)_{\mu}).
\]
The claim follows by observing that a linearly independent set with cardinality equal to the dimension of the vector space must be a spanning set.
\end{proof}

\begin{lemma}\label{lem:41'imp1'}
If $S_4(\unw, \varpi)$, and $S_1'(\unx, \psi)$ whenever $\wt(\unx\otimes\psi) \le \wt \unw$, then $S_1'(\unw, \varpi)$. \end{lemma}

\begin{proof}
The right hand side of the equation ${^{p\ell}\kappa}_{\unw, \varpi}^{\unu} = {^{p\ell}\rho}_{\unw, \varpi}^{\unu}(\kappa)$ only involves terms ${^{ij}\kappa}_{\unx, \psi}^{\uny}$ such that $\wt(\unx\ot \psi) \le \wt \unw$. If we write $ {^{p\ell}\rho}_{\unw, \varpi}^{\unu}(t)$ to denote the same formula with each ${^{ij}\kappa}_{\unx, \psi}^{\uny}$ replaced by ${^{ij}t}_{\unx, \psi}^{\uny}$, then our hypotheses imply that
\[
{^{p\ell}\kappa}_{\unw, \varpi}^{\unu} = {^{p\ell}\rho}_{\unw, \varpi}^{\unu}(t).
\]
Thus, to show that $S_1'(\unw, \varpi)$ holds we must verify the following equality of rational functions in $\mathbb{C}(q)$:
\[
{^{p\ell}t}_{\unw, \varpi}^{\unu} = {^{p\ell}\rho}_{\unw, \varpi}^{\unu}(t),
\]
which we verified using the SAGE code included with the source file of the arXiv submission of this paper.
\end{proof}

\begin{example}
We take verification of Equation \eqref{Recursion(10)} as an example. In order to verify that 

$$t_{(a,b),\varpi_2}^{(1,0)} ={\rho}_{(a,b),\varpi_2}^{(1,0)}(t)= \dfrac{[7]}{[3]}
-\dfrac{1}{t_{(a-1,b),\varpi_1}^{(-1,1)}} t_{(a-2,b+1),\varpi_2}^{(3,-1)}
          -\dfrac{1}{t_{(a-1,b),\varpi_1}^{(2,-1)}} {\ }_{,} $$
we first write the $t_{(x,y), \varpi}^{(s,t)}$'s explicitly using Equations \eqref{texamplea}, \eqref{texampleb}, \eqref{texamplec}, and \eqref{texampled} to obtain 
\begin{align*}
\dfrac{[a+3][a+3b+6]}{[3][a][a+3b+3]}&=\dfrac{[7]}{[3]}-\dfrac{1}{\Bigg(-\dfrac{[a]}{[a-1]}\Bigg)}\Bigg(-\dfrac{[3b+6]}{[3b+3]}\Bigg) - \dfrac{1}{\dfrac{[3b+3][a+3b+3]}{[3b][a+3b+2]}}.
\end{align*}
We can rewrite this as:
\begin{align*}
\dfrac{(q^{a+3}-q^{-a-3})(q^{a+3b+6}-q^{-a-3b-6})(q-q^{-1})}
            {(q^3-q^{-3})(q^a-q^{-a})(q^{a+3b+3}-q^{-a-3b-3})}
&= \dfrac{q^7-q^{-7}}{q^3-q^{-3}} 
- \dfrac{(q^{a-1}-q^{-a+1})(q^{3b+6}-q^{-3b-6})}
          {(q^a-q^{-a})(q^{(3b+3)-q^{-3b-3}})} \\
&- \dfrac{(q^{3b}-q^{-3b})(q^{a+3b+2}-q^{-a-3b-2})}
           {(q^{3b+3}-q^{-3b-3})(q^{a+3b+3}-q^{-a-3b-3})}.
\end{align*}
Making the substitutions $A=q^a$ and $B=q^b$, we obtain:
\begin{align*}
\dfrac{(Aq^3-A^{-1}q^{-3})(AB^3q^6-A^{-1}B^{-3}q^{-6})(q-q^{-1})}
        {(q^3-q^{-3})(A-A^{-1})(AB^3q^3-A^{-1}B^{-3}q^{-3})}
&= \dfrac{q^7-q^{-7}}{q^3-q^{-3}} - 
\dfrac{(Aq^{-1}-A^{-1}q)(B^3q^6-B^{-3}q^{-6})}{(A-A^{-1})(B^3q^3-B^{-3}q^{-3})} \\
&-\dfrac{(B^3-B^{-3})(AB^3q^2-A^{-1}B^{-3}q^{-2})}
            {(B^3q^3-B^{-3}q^{-3})(AB^3q^3-A^{-1}B^{-3}q^{-3})}.
\end{align*}
Then we can use .simplify\_full() in SAGE to simplify the rational function of $A$, $B$, and $q$, which is given by the difference of the left hand side and right hand side of the above equation. The result computed by SAGE is equal to 0, which tells us that Equation \eqref{Recursion(10)} holds.     
\end{example}

\begin{lemma}\label{lem:11'imp5}
If $S_1(\unx)$ for all $\unx$ such that $\wt \unx < \wt (\unw\otimes\varpi) $, and $S_1'(\unw, \varpi)$, then $S_5(\unw, \varpi)$. \end{lemma}

\begin{proof}
By the definition of $T_{\unw\otimes \varpi}$ we find
\[
T_{\unw\otimes\varpi} := T_{\unw}\otimes \id_{\varpi}  - \sum_{\substack{\mu \in \wt V(\varpi)\backslash\lbrace{\varpi}\rbrace \\ v_i, v_j\in \text{Ker}_{\wt \unw}(V(\varpi)_{\mu})}} (^{p\ell}t_{\unw, \varpi}^{\unu_{\wt\unw + \mu}})^{-1}_{ij}\cdot{^{ij}\mathbb{TT}}_{\unw, \varpi}^{\unu_{\wt\unw+ \mu}}. 
\]
Then by our hypotheses, we deduce that
\[
T_{\unw\otimes\varpi}= C_{\unw}\otimes \id_{\varpi} - \sum_{\substack{\mu \in \wt V(\varpi)\backslash\lbrace{\varpi}\rbrace \\ v_i, v_j\in \text{Ker}_{\wt \unw}(V(\varpi)_{\mu})}} ({^{p\ell}\kappa}_{\unw, \varpi}^{\unu_{\wt\unw+ \mu}})^{-1}_{ij}\cdot{^{ij}\mathbb{LL}}_{\unw, \varpi}^{\unu_{\wt\unw+ \mu}}.
\]
\end{proof}

\begin{lemma}\label{lem:5imp6}
If $S_5(\unw, \varpi)$, then $S_6(\unw\otimes\varpi)$.
\end{lemma}
\begin{proof}

By Lemma \ref{claspschur} we deduce the following multiplication formula for double ladders:
\[
{^{p\ell}\mathbb{LL}}_{\unw, \varpi}^{\unu}\circ {^{rs}\mathbb{LL}}_{\unw, \varpi}^{\unv} = \delta_{\wt \unu, \wt \unv}{^{\ell r}\kappa}_{\unw, \varpi}^{\unu} \cdot {^{ps}\mathbb{LL}}_{\unw, \varpi}^{\unu}.
\]
Using the expression for $T_{\unw\otimes \varpi}$ from $S_5(\unw, \varpi)$ and the above formula, one can explicitly compute to verify that $T_{\unw, \varpi}$ is idempotent.

\end{proof}

\begin{lemma}\label{lem:53imp6'}               
If $S_5(\unw, \varpi)$ and $S_3(\unw, \varpi)$, then $S_6'(\unw, \varpi)$. 
\end{lemma} 
\begin{proof}
Write $\lambda = \wt\unw$. Let $\mu\in \wt V(\varpi)\backslash\lbrace\varpi\rbrace$ and let $\unu\in \DD$ such that $\wt \unu = \lambda + \mu$. Let $D\in \Hom_{\DD}(\unw\otimes \varpi, \unu)$.  Consider the neutral diagram $\mathsf{H}_{\unu}^{\unu_{\lambda+\mu}}: \unu\rightarrow \unu_{\lambda + \mu}$ and write $D'= \mathsf{H}_{\unu}^{\unu_{\lambda+\mu}} \circ C_{\unu}\circ D$. 

Combining that clasps are idempotent with Lemma \ref{claspschur} we find
\[
C_{\unu_{\lambda+ \mu}}\circ D' \circ {^{ij}\mathbb{LL}}_{\unw, \varpi}^{\unu_{\lambda+ \nu}} = \delta_{\mu, \nu} \cdot C_{\unu_{\lambda+ \mu}} \circ D' \circ (C_{\unw}\otimes \id_{\varpi}) \circ {^{ij}\mathbb{LL}}_{\unw, \varpi}^{\unu_{\lambda+ \mu}}.
\]
By $S_3(\unw, \varpi)$ there are scalars $\xi_k$ such that 
\[
C_{\unu_{\lambda+ \mu}}\circ D' \circ (C_{\unw} \otimes \id_{\varpi})= \sum_{v_k\in \text{Ker}_{\wt \unw}(V(\varpi)_{\mu})} \xi_k\cdot {^kLL}_{\unw, \varpi}^{\unu_{\lambda+ \mu}}.
\]
Thus, using $S_5(\unw, \varpi)$ we can rewrite  
$C_{\unu_{\lambda+ \mu}}\circ D' \circ T_{\unw\otimes \varpi}$ 
as
\begin{align*}
&\ \ \ \ C_{\unu_{\lambda+ \mu}}\circ D' \circ (C_{\unw} \otimes \id_{\varpi}) - \sum_{v_i,v_j\in \text{Ker}_{\lambda}(V(\varpi)_{\mu})} ({^{p\ell}\kappa}_{\unw, \varpi}^{\unu_{\lambda+ \mu}})^{-1}_{ij}\cdot C_{\unu_{\lambda+ \mu}} \circ D' \circ (C_{\unw}\otimes \id_{\varpi}) \circ {^{ij}\mathbb{LL}}_{\unw, \varpi}^{\unu_{\lambda+ \mu}}
\\
&= C_{\unu_{\lambda+ \mu}}\circ D' \circ (C_{\unw} \otimes \id_{\varpi}) -  \sum_{v_i,v_j\in \text{Ker}_{\wt\unw}(V(\varpi)_{\mu})}\sum_k ({^{p\ell}\kappa}_{\unw, \varpi}^{\unu_{\lambda+ \mu}})^{-1}_{ij}\xi_k\cdot {^{k}LL}_{\unw, \varpi}^{\unu_{\lambda+ \mu}} \circ {^{ij}\mathbb{LL}}_{\unw, \varpi}^{\unu_{\lambda+ \mu}} \\
&= C_{\unu_{\lambda+ \mu}}\circ D' \circ (C_{\unw} \otimes \id_{\varpi}) -  \sum_{v_i,v_j\in \text{Ker}_{\wt\unw}(V(\varpi)_{\mu})}\sum_k ({^{p\ell}\kappa}_{\unw, \varpi}^{\unu_{\lambda+ \mu}})^{-1}_{ij}\xi_k{^{ki}\kappa_{\unw, \varpi}^{\unu_{\lambda+ \mu}}}\cdot {^{j}LL}_{\unw, \varpi}^{\unu_{\lambda+ \mu}} \\
&= C_{\unu_{\lambda+ \mu}}\circ D' \circ (C_{\unw} \otimes \id_{\varpi}) - \sum_{v_j\in \text{Ker}_{\wt\unw}(V(\varpi)_{\mu})}\sum_k \xi_k\delta_{k, j}\cdot {^{j}LL}_{\unw, \varpi}^{\unu_{\lambda+ \mu}} \\
&= C_{\unu_{\lambda+ \mu}}\circ D' \circ (C_{\unw} \otimes \id_{\varpi})-\sum_{v_j\in \text{Ker}_{\wt \unw}(V(\varpi)_{\mu})} \xi_j \cdot {^jLL}_{\unw, \varpi}^{\unu_{\lambda+ \mu}} \\
&= 0.
\end{align*}

Using Lemmas \ref{claspabsorption} and \ref{neturalabsorption} it is not hard to see that  $C_{\unu}\circ \mathsf{H}_{\unu_{\lambda+\mu}}^{\unu} \circ C_{\unu_{\lambda+ \mu}}\circ D'= C_{\unu}\circ D$, and it follows that $C_{\unu}\circ D\circ T_{\unw\otimes \varpi} = 0$.
\end{proof}

\begin{lemma}\label{lem:22'imp4}
Let $\unw\in \DD$ and let $\varpi$ be a fundamental weight. If $S_2(\unx, \psi)$ and $S_2'(\unx, \psi)$ whenever $\wt(\unx\otimes\psi)\le \wt\unw$, then $S_4(\unw, \varpi)$. \end{lemma}

\begin{proof}
Consider $\unx ,\psi$ such that $\wt(\unx \otimes \psi)\le \wt \unw$. By $S_2(\unx, \psi)$ we obtain the following. 
\[
C_{\unx}\otimes \id_{\psi} = \sum_{\substack{\mu \in \wt V(\psi) \\ v_i, v_j\in \text{Ker}_{\wt \unx} (V(\psi)_{\mu})}} {^{ij}\xi_{\unx, \psi}^{\unu_{\wt\unx+ \mu}}} \cdot {^{ij}\mathbb{LL}}_{\unx, \psi}^{\unu_{\wt\unx+ \mu}} 
\]
Postcomposing with ${^pLL}_{\unx, \psi}^{\unu_{\wt\unx+ \mu}}$ and using Lemma \ref{claspschur} results in the next sequence of equalities.
\begin{align*}
{^pLL}_{\unx, \psi}^{\unu_{\wt\unx+ \mu}} &= \sum_{v_i, v_j\in \text{Ker}_{\wt \unx} (V(\psi)_{\mu})}  {^{ij}\xi_{\unx, \psi}^{\unu_{\wt\unx+ \mu}}} \cdot {^pLL}_{\unx, \psi}^{\unu_{\wt\unx+ \mu}}\circ {^{ij}\mathbb{LL}}_{\unx, \psi}^{\unu_{\wt\unx+ \mu}} \\
&= \sum_{v_i, v_j\in \text{Ker}_{\wt \unx} (V(\psi)_{\mu})} {^{ij}\xi_{\unx, \psi}^{\unu_{\wt\unx+ \mu}}} \cdot {^{pi}\kappa}_{\unx, \psi}^{\unu_{\wt\unx+ \mu}} \cdot {^jLL}_{\unx, \psi}^{\unu_{\wt\unx+ \mu}}
\end{align*}
By $S_2'(\unx, \psi)$ it follows that
\[
 \sum_{v_i\in \text{Ker}_{\wt \unx} (V(\psi)_{\mu})} {^{ij}\xi}_{\unx, \psi}^{\unu_{\wt\unx+ \mu}}\cdot {^{pi}\kappa}_{\unx, \psi}^{\unu_{\wt\unx+ \mu}} = \delta_{jp}.
\]
Moreover, using that clasps are idempotent along with Lemma \ref{claspabsorption} we find
\[
LL_{\unx, \psi}^{\unu_{\wt\unx+\psi}}\circ \mathbb{D}(LL_{\unx, \psi}^{\unu_{\wt\unx+ \psi}})= C_{\unu_{\wt\unx+ \psi}}\circ (C_{\unx}\otimes\id_{\psi}) \circ C_{\unu_{\wt\unx+ \psi}} = C_{\unu_{\wt\unx+ \psi}},
\]
and it follows that $\kappa_{\unx, \psi}^{\unu_{\wt\unx + \psi}}= 1$. 
Thus, $\xi_{\unx, \psi}^{\unu_{\wt\unx + \psi}} = 1$ and 
\begin{equation}\label{lemma3clasp}
C_{\unx\otimes \psi} = C_{\unx}\otimes \id_{\psi}- \sum_{\substack{\mu \in \wt V(\psi) \backslash\lbrace{\psi}\rbrace \\ v_i, v_j\in \text{Ker}_{\wt \unx} (V(\psi)_{\mu})}} ({^{p\ell}\kappa}_{\unx, \psi}^{\unu_{\wt\unx+ \mu}})^{-1}_{ij} \cdot {^{ij}\mathbb{LL}}_{\unx, \psi}^{\unu_{\wt\unx+ \mu}}.
\end{equation}

Observe that
\begin{align*}
{^{p\ell}\kappa_{\unw, \varpi}^{\unu}}C_{\unu} &= {^pLL}_{\unw, \varpi}^{\unu} \circ \mathbb{D}({^\ell LL}_{\unw, \varpi}^{\unu}) \\
&= C_{\unu}\circ {^pELL}_{\unw, \varpi}^{\unu} \circ (C_{\unw}\otimes \id_{\varpi}) \circ \mathbb{D}({^\ell ELL}_{\unw, \varpi}^{\unu}) \circ C_{\unu}.
\end{align*}
Then use Equation \eqref{lemma3clasp} for $\unw = \unx\otimes \psi$ to rewrite the $C_{\unw}$ term on the right hand side. This new sum will reduce to a scalar multiple of $C_{\unu}$ by repeatedly applying graphical reductions or by replacing another clasp, necessarily of the form $C_{\uny\otimes\varpi}$ for some $\uny , \varpi$ such that $\wt(\uny\ot \varpi)\le \wt \unw$, using Equation \eqref{lemma3clasp}. The exact form of the coefficient is determined via the calculations in Appendix \ref{GraphicalCalc}, where it is shown to be equal to ${^{p\ell}\rho}_{\unw, \varpi}^{\unu}(\kappa)$. Therefore, ${^{p\ell}\kappa_{\unw, \varpi}^{\unu} C_{\unu} = {^{p\ell}\rho}}_{\unw, \varpi}^{\unu}(\kappa) C_{\unu}$ and the desired result follows from looking at the neutral coefficient of each map.
\end{proof}

\begin{example}\label{ex:diagrammatic-calc-example}
The above argument is best illustrated by example. Consider $\unw$ with $\wt \unw= (a, b)$ and assume $S_2(\unx, \psi)$ and $S_2'(\unx, \psi)$ whenever $\wt(\unx\otimes\psi)\le \wt\unw$. Note that 
\[
\rho_{(a,b), \varpi_1}^{(-1, 1)}:= -[2] - \frac{1}{x_{(a-1, b), \varpi_1}^{(-1, 1)}}.
\]
We will show that $\kappa_{\unw, \varpi_1}^{\unu_{(a-1, b+1)}}= \rho_{(a, b), \varpi_1}^{(-1,1)}(\kappa)$.

Let $\wt(\unv\ot\varpi_1)= \wt \unw$. By definition we have
$\kappa_{\unw, \varpi_1}^{\unu_{(a-1, b+1)}}C_{\unu_{(a-1, b+1)}}$ is equal to
\[
C_{\unu_{(a-1, b+1)}}\circ {ELL}_{\unw, \varpi_1}^{\unu_{(a-1, b+1)}} \circ \left(C_{\unv\otimes \varpi_1}\otimes \id_{\varpi_1}\right) \circ \mathbb{D}({ELL}_{\unw, \varpi_1}^{\unu_{(a-1, b+1)}}) \circ C_{\unu_{(a-1, b+1)}}.
\]
As in the first half of the proof of Lemma \ref{lem:22'imp4}, our hypotheses allow us to write
\[
C_{\unv\otimes {\varpi_1}} = C_{\unv}\otimes \id_{\varpi_1} - \sum_{\substack{\mu \in \wt V(\varpi_1)\setminus \{\varpi_1\} \\ v_i, v_j\in \text{Ker}_{\wt \unv} (V(\varpi_1)_{\mu})}} ({^{p\ell}\kappa}_{\unv, \varpi_1}^{\unu_{\wt\unv + \mu}})^{-1}_{ij} \cdot {^{ij}\mathbb{LL}}_{\unv, \varpi_1}^{\unu_{\wt\unv + \mu}}.
\]
Using Lemma \ref{claspedwebspace}, we observe that if $\mu\ne (-1,1)$, then
\[
C_{\unu_{(a-1, b+1)}}\circ {ELL}_{\unw, \varpi_1}^{\unu_{(a-1, b+1)}} \circ \left ( ({^{p\ell}\kappa}_{\unv, \varpi_1}^{\unu_{\wt\unv + \mu}})^{-1}_{ij} \cdot {^{ij}\mathbb{LL}}_{\unv, \varpi_1}^{\unu_{\wt\unv + \mu}}\otimes \id_{\varpi_1}\right) \circ \mathbb{D}({ELL}_{\unw, \varpi_1}^{\unu_{(a-1, b+1)}}) \circ C_{\unu_{(a-1, b+1)}}
\]
is zero. Finally, applying web relations (and properties of clasps) we find
\begin{equation*}
\kappa_{(a,b), \varpi_1}^{(-1,1)}C_{\unu_{(a-1, b+1)}}=\left(-[2]-\frac{1}{\kappa_{(a-1,b), \varpi_1}^{(-1,1)}}\right) C_{\unu_{(a-1, b+1)}}.
\end{equation*}
We conclude with a schematic of the graphical calculations involved.

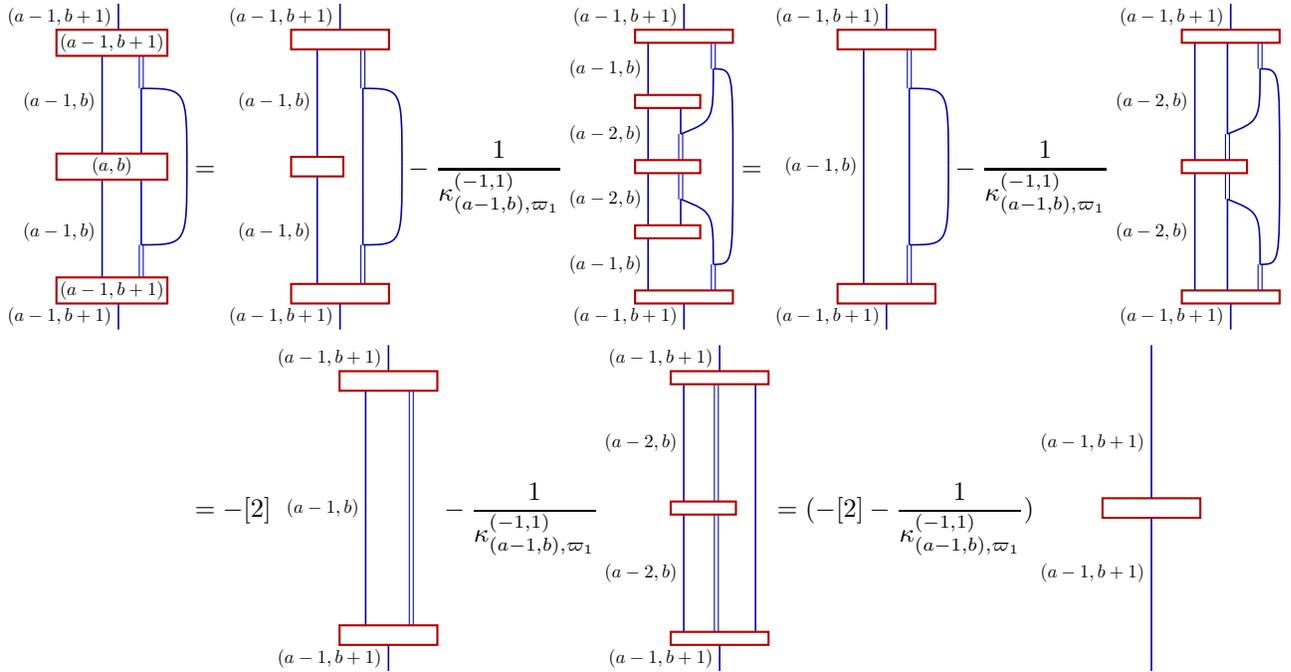
\def\Kba{
\begin{tric}
\draw[double,thin] (1.1,1.7)--(1.1,1.2);
\draw (1.1,1.2)--(1.1,0.2);
\draw[double,thin] (1.1,-1.7)--(1.1,-1.2);
\draw (1.1,-1.2)--(1.1,-0.2);
\draw (1.1,1.2)..controls(1.8,1.2)..(1.8,0);
\draw (1.1,-1.2)..controls(1.8,-1.2)..(1.8,0);

\draw (0.75,2.1)--(0.75,2.5);
\draw (0.75,-2.1)--(0.75,-2.5);
\draw (0.5,1.7)--(0.5,0.2);
\draw (0.5,-1.7)--(0.5,-0.2);
\draw (0.75,2.3) node[anchor=east,black,scale=0.7] {$(a-1,b+1)$};
\draw (0.75,-2.3) node[anchor=east,black,scale=0.7] {$(a-1,b+1)$};
\draw (0.5,1) node[anchor=east,black,scale=0.7] {$(a-1,b)$};
\draw (0.5,-1) node[anchor=east,black,scale=0.7] {$(a-1,b)$};
\draw[darkred,thick] (-0.2,-0.2) rectangle (1.5,0.2) node[pos=0.5,scale=0.7,black] {$(a,b)$};
\draw[darkred,thick] (-0.2,1.7) rectangle (1.5,2.1) node[pos=0.5,scale=0.7,black] {$(a-1,b+1)$};
\draw[darkred, thick] (-0.2,-2.1) rectangle (1.5,-1.7) node[pos=0.5,scale=0.7,black] {$(a-1,b+1)$};
\end{tric}
}

\def\Kbba{
\begin{tric}
\draw (0.75,2.1)--(0.75,2.5);
\draw (0.75,-2.1)--(0.75,-2.5);
\draw (0.4,1.8)--(0.4,0.15); \draw(0.4,-0.15)--(0.4,-1.8);
\draw (0.75,2.3) node[anchor=east,black,scale=0.7] {$(a-1,b+1)$};
\draw (0.75,-2.3) node[anchor=east,black,scale=0.7] {$(a-1,b+1)$};
\draw (0.4,1) node[anchor=east,black,scale=0.7] {$(a-1,b)$};
\draw (0.4,-1) node[anchor=east,black,scale=0.7] {$(a-1,b)$};
\draw[double,thin] (1.1,1.8)--(1.1,1.2);
\draw (1.1,1.2)--(1.1,-1.2);
\draw[double,thin] (1.1,-1.8)--(1.1,-1.2);
\draw (1.1,1.2)..controls(1.7,1.2)..(1.7,0);
\draw (1.1,-1.2)..controls(1.7,-1.2)..(1.7,0);
\draw[darkred,thick] (0,0.15)rectangle(0.8,-0.15);
\draw[darkred,thick] (0,1.8) rectangle (1.5,2.1) ;
\draw[darkred, thick] (0,-2.1) rectangle (1.5,-1.8) ;
\end{tric}
}

\def\Kbbb{
\begin{tric}
\draw (0.75,2.1)--(0.75,2.5);
\draw (0.75,-2.1)--(0.75,-2.5);
\draw (0.2,1.9)--(0.2,1.1);
\draw (0.2,-1.9)--(0.2,-1.1);
\draw (0.2,0.1)--(0.2,0.9);
\draw (0.2,-0.1)--(0.2,-0.9);
\draw (0.75,2.3) node[anchor=east,black,scale=0.7] {$(a-1,b+1)$};
\draw (0.75,-2.3) node[anchor=east,black,scale=0.7] {$(a-1,b+1)$};
\draw (0.2,1.5) node[anchor=east,black,scale=0.7] {$(a-1,b)$};
\draw (0.2,-1.5) node[anchor=east,black,scale=0.7] {$(a-1,b)$};
\draw (0.2,0.5) node[anchor=east,black,scale=0.7] {$(a-2,b)$};
\draw (0.2,-0.5) node[anchor=east,black,scale=0.7] {$(a-2,b)$};

\draw [double,thin] (0.7,0.1)--(0.7,0.5);
\draw  (0.7,0.5)--(0.7,0.9);
\draw (0.7,0.5)..controls(1.2,0.7)..(1.2,1.5);
\draw [double,thin](1.2,1.5)--(1.2,1.9);
\draw (1.2,1.5).. controls (1.5,1.5)..(1.5,0);
\draw [double,thin] (0.7,-0.1)--(0.7,-0.5);
\draw  (0.7,-0.5)--(0.7,-0.9);
\draw (0.7,-0.5)..controls(1.2,-0.7)..(1.2,-1.5);
\draw [double,thin](1.2,-1.5)--(1.2,-1.9);
\draw (1.2,-1.5).. controls (1.5,-1.5)..(1.5,0);

\draw[darkred,thick] (0,-0.1) rectangle (1,0.1) ;
\draw[darkred,thick] (0,0.9) rectangle (1,1.1) ;
\draw[darkred, thick] (0,-0.9) rectangle (1,-1.1) ;
\draw[darkred, thick] (0,1.9) rectangle (1.5,2.1) ;
\draw[darkred, thick] (0,-1.9) rectangle (1.5,-2.1) ;
\end{tric}
}

\def\Kbca{
\begin{tric}
\draw (0.75,2.1)--(0.75,2.5);
\draw (0.75,-2.1)--(0.75,-2.5);
\draw (0.4,1.8)--(0.4,-1.8);
\draw (0.75,2.3) node[anchor=east,black,scale=0.7] {$(a-1,b+1)$};
\draw (0.75,-2.3) node[anchor=east,black,scale=0.7] {$(a-1,b+1)$};
\draw (0.4,0) node[anchor=east,black,scale=0.7] {$(a-1,b)$};
\draw[double,thin] (1.1,1.8)--(1.1,1.2);
\draw (1.1,1.2)--(1.1,-1.2);
\draw[double,thin] (1.1,-1.8)--(1.1,-1.2);
\draw (1.1,1.2)..controls(1.7,1.2)..(1.7,0);
\draw (1.1,-1.2)..controls(1.7,-1.2)..(1.7,0);
\draw[darkred,thick] (0,1.8) rectangle (1.5,2.1) ;
\draw[darkred, thick] (0,-2.1) rectangle (1.5,-1.8) ;
\end{tric}
}

\def\Kbcb{
\begin{tric}
\draw (0.75,2.1)--(0.75,2.5);
\draw (0.75,-2.1)--(0.75,-2.5);
\draw (0.2,0.1)--(0.2,1.9);
\draw (0.2,-0.1)--(0.2,-1.9);
\draw (0.75,2.3) node[anchor=east,black,scale=0.7] {$(a-1,b+1)$};
\draw (0.75,-2.3) node[anchor=east,black,scale=0.7] {$(a-1,b+1)$};
\draw (0.2,1) node[anchor=east,black,scale=0.7] {$(a-2,b)$};
\draw (0.2,-1) node[anchor=east,black,scale=0.7] {$(a-2,b)$};

\draw [double,thin] (0.7,0.1)--(0.7,0.5);
\draw  (0.7,0.5)--(0.7,1.9);
\draw (0.7,0.5)..controls(1.2,0.7)..(1.2,1.5);
\draw [double,thin](1.2,1.5)--(1.2,1.9);
\draw (1.2,1.5).. controls (1.5,1.5)..(1.5,0);
\draw [double,thin] (0.7,-0.1)--(0.7,-0.5);
\draw  (0.7,-0.5)--(0.7,-1.9);
\draw (0.7,-0.5)..controls(1.2,-0.7)..(1.2,-1.5);
\draw [double,thin](1.2,-1.5)--(1.2,-1.9);
\draw (1.2,-1.5).. controls (1.5,-1.5)..(1.5,0);

\draw[darkred,thick] (0,-0.1) rectangle (1,0.1) ;
\draw[darkred, thick] (0,1.9) rectangle (1.5,2.1) ;
\draw[darkred, thick] (0,-1.9) rectangle (1.5,-2.1) ;
\end{tric}
}

\def\Kbda{
\begin{tric}
\draw (0.75,2.1)--(0.75,2.5);
\draw (0.75,-2.1)--(0.75,-2.5);
\draw (0.4,1.8)--(0.4,-1.8);
\draw (0.75,2.3) node[anchor=east,black,scale=0.7] {$(a-1,b+1)$};
\draw (0.75,-2.3) node[anchor=east,black,scale=0.7] {$(a-1,b+1)$};
\draw (0.4,0) node[anchor=east,black,scale=0.7] {$(a-1,b)$};
\draw (1.1,1.2)--(1.1,-1.2);
\draw[double,thin] (1.1,-1.8)--(1.1,1.8);
\draw[darkred,thick] (0,1.8) rectangle (1.5,2.1) ;
\draw[darkred, thick] (0,-2.1) rectangle (1.5,-1.8) ;
\end{tric}
}

\def\Kbdb{
\begin{tric}
\draw (0.75,2.1)--(0.75,2.5);
\draw (0.75,-2.1)--(0.75,-2.5);
\draw (0.2,0.1)--(0.2,1.9);
\draw (0.2,-0.1)--(0.2,-1.9);
\draw (0.75,2.3) node[anchor=east,black,scale=0.7] {$(a-1,b+1)$};
\draw (0.75,-2.3) node[anchor=east,black,scale=0.7] {$(a-1,b+1)$};
\draw (0.2,1) node[anchor=east,black,scale=0.7] {$(a-2,b)$};
\draw (0.2,-1) node[anchor=east,black,scale=0.7] {$(a-2,b)$};

\draw [double,thin] (0.7,0.1)--(0.7,1.9);
\draw [double,thin] (0.7,-0.1)--(0.7,-1.9);
\draw (1.3,1.9)--(1.3,-1.9);

\draw[darkred,thick] (0,-0.1) rectangle (1,0.1) ;
\draw[darkred, thick] (0,1.9) rectangle (1.5,2.1) ;
\draw[darkred, thick] (0,-1.9) rectangle (1.5,-2.1) ;
\end{tric}
}

\def\Kbe{
\begin{tric}
\draw (0.75,0.15)--(0.75,2.5);
\draw (0.75,-0.15)--(0.75,-2.5);
\draw (0.75,1) node[anchor=east,black,scale=0.7] {$(a-1,b+1)$};
\draw (0.75,-1) node[anchor=east,black,scale=0.7] {$(a-1,b+1)$};
\draw[darkred,thick] (0,-0.15) rectangle (1.5,0.15);
\end{tric}
}

\begin{align*}
\Kba&=\Kbba-\frac{1}{\kappa_{(a-1,b), \varpi_1}^{(-1,1)}} \Kbbb
= \Kbca -\frac{1}{\kappa_{(a-1,b), \varpi_1}^{(-1,1)}}\Kbcb \\
&=-[2] \Kbda -\frac{1}{\kappa_{(a-1,b), \varpi_1}^{(-1,1)}}\Kbdb
= (-[2]-\frac{1}{\kappa_{(a-1,b), \varpi_1}^{(-1,1)}})\Kbe
\end{align*}
\end{example}

Finally, we combine the previous lemmas to deduce the result of our main theorem. 

\begin{proof}[Proof of Theorem \ref{mainthmtoprove}]
We will prove the result by induction on $\wt\unw$ with respect to $\le$. The base case follows from observing that $T_{\emptyset}$ is $1$ times the empty diagram, which agrees with $C_{\emptyset}$. Assume that $S_1(\unx)$ holds for all $\unx$ such that $\wt\unx<\wt\unw$ and assume that $S_1'(\uny, \psi)$ holds for all $\uny,\psi$ such that $\wt(\uny\ot \psi)< \wt\unw$. We will show $S_1(\unw)$ and $S_1'(\unw', \varpi)$, where $\unw = \unw'\otimes \varpi$.

Consider $\uny,\psi$ such that $\wt(\uny\otimes \psi)< \wt\unw$. Then $S_1(\unx)$ holds whenever $\wt\unx\le \wt(\uny\otimes \psi)$, and by Lemma \ref{lem:1imp2} we deduce $S_2(\uny, \psi)$. Thus, $S_2(\uny, \psi)$ holds for all $\uny,\psi$ such that $\wt(\uny\otimes\psi)<\wt \unw$. 

If $\uny,\psi$ is such that $\wt(\uny\otimes \psi)< \wt\unw$, then our inductive hypothesis also says that $S_1'(\uny, \psi)$ holds. Along with Lemmas \ref{lem:1'imp2'} and \ref{lem:2'imp3}, this implies $S_2'(\uny, \psi)$ and $S_3(\uny, \psi)$. Hence, $S_2'(\uny, \psi)$ and $S_3(\uny, \psi)$ holds for all $\uny,\psi$ such that $\wt(\uny\otimes \psi)< \wt \unw$. 

For $\unx$ such that $\wt\unx< \wt \unw$, we have $S_2(\uny, \psi)$ and $S_2'(\uny, \psi)$ whenever $\wt(\uny\ot \psi)\le \wt \unx$. So from Lemma \ref{lem:22'imp4} we deduce $S_4(\unx, \varpi)$ for all $\unx$ such that $\wt\unx< \wt\unw$ and for arbitrary $\varpi$. 

If $\unw= \unw'\otimes\varpi$, then $\wt \unw'<\wt \unw$, so $S_4(\unw', \varpi)$. Also, if $\wt(\uny\otimes\psi)\le \wt\unw'$, then $\wt(\uny\otimes\psi)< \wt\unw$ so $S_1'(\uny, \psi)$ holds whenever $\wt(\uny\otimes \psi)\le\wt \unw'$. Thus, Lemma \ref{lem:41'imp1'} implies $S_1'(\unw', \varpi)$.

At this point, we know that $S_1(\unx)$ whenever $\wt\unx< \wt(\unw'\otimes \varpi)$ and that $S_1'(\unw', \varpi)$ holds. Therefore, Lemma \ref{lem:11'imp5} implies that $S_5(\unw', \varpi)$ holds too. Then from Lemma \ref{lem:5imp6} we deduce $S_6(\unw'\otimes\varpi)$ is true.

Moreover, since $S_1'(\unw', \varpi)$ is true, Lemmas \ref{lem:1'imp2'} and \ref{lem:2'imp3} together imply $S_3(\unw', \varpi)$. Therefore, we can use Lemma \ref{lem:53imp6'} to deduce $S_6'(\unw', \varpi)$. 

If we show that $T_{\unw'\otimes\varpi}\ne 0$, then Definition \ref{gclaspdef} and Lemma \ref{uniqueessofclasps} will tell us that $S_6(\unw'\otimes\varpi)$ and $S_6'(\unw', \varpi)$ imply $S_1(\unw'\otimes \varpi)$, so we are then done by induction. To see that $T_{\unw'\otimes\varpi}$ is not $0$, we apply $\Phi$ from Theorem \ref{T:Kupes-functor-Phi} and evaluate on a weight vector in $V(\unw'\otimes\varpi)_{\wt\unw'+ \varpi}$. Using $S_5(\unw', \varpi)$, along with the observations that for all $\unu$ such that $\wt\unu- \wt \unw'= \mu\in \wt V(\varpi) \backslash\lbrace{\varpi}\rbrace$, the map $\Phi\left({^{ij}\mathbb{LL}}_{\unw', \varpi}^{\unu}\right)$ acts as zero on $V(\unw'\otimes\varpi)_{\wt\unw'+ \varpi}$ (since these maps factor through representations which do not have $\wt\unw'+ \varpi$ as a weight) and $\Phi\left(C_{\unw'}\otimes \id_{\varpi}\right)$ acts on $V(\unw'\otimes\varpi)_{\wt\unw'+ \varpi}$ as multiplication by $1$, we deduce that $T_{\unw'\otimes \varpi}$ is non-zero.
\end{proof}


\bibliographystyle{plain}
\bibliography{mastercopy}

\begin{thebibliography}{10}

\bibitem{Sp4tilt}
Elijah Bodish.
\newblock Web calculus and tilting modules in type {$C_2$}.
\newblock Preprint, 2020.
\newblock arXiv:2009.13786.

\bibitem{Sp4Clasp}
Elijah Bodish.
\newblock Triple clasp formulas for {$C_2$} webs.
\newblock {\em J. Algebra}, 604:324--361, 2022.

\bibitem{bodish2021type}
Elijah Bodish, Ben Elias, David E.~V. Rose, and Logan Tatham.
\newblock Type {$C$} webs.
\newblock Preprint, 2021.
\newblock arXiv:2103.14997.

\bibitem{CKM}
Sabin Cautis, Joel Kamnitzer, and Scott Morrison.
\newblock Webs and quantum skew {H}owe duality.
\newblock {\em Math. Ann.}, 360(1-2):351--390, 2014.

\bibitem{MR2901969}
Benjamin Cooper and Vyacheslav Krushkal.
\newblock Categorification of the {J}ones-{W}enzl projectors.
\newblock {\em Quantum Topol.}, 3(2):139--180, 2012.

\bibitem{MR2873427}
Ben Elias.
\newblock {\em Soergel {D}iagrammatics for {D}ihedral {G}roups}.
\newblock ProQuest LLC, Ann Arbor, MI, 2011.
\newblock Thesis (Ph.D.)--Columbia University.

\bibitem{elias2015light}
Ben Elias.
\newblock Light ladders and clasp conjectures.
\newblock Preprint, 2015.
\newblock arXiv:1510.06840.

\bibitem{JantzenQgps}
Jens~Carsten Jantzen.
\newblock {\em Lectures on Quantum Groups}, volume~6 of {\em Graduate Studies
  in Mathematics}.
\newblock American Mathematical Society, Providence, RI, first edition, 1996.

\bibitem{Jon2}
Vaughan F.~R. Jones.
\newblock Index for subfactors.
\newblock {\em Invent. Math.}, 72(1):1--25, 1983.

\bibitem{JonesPolynomial}
Vaughan F.~R. Jones.
\newblock A polynomial invariant for knots via von {N}eumann algebras.
\newblock {\em Bull. Amer. Math. Soc. (N.S.)}, 12(1):103--111, 1985.

\bibitem{Kauffman-Lins}
Louis~H. Kauffman and S\'{o}stenes~L. Lins.
\newblock {\em Temperley-{L}ieb recoupling theory and invariants of
  {$3$}-manifolds}, volume 134 of {\em Annals of Mathematics Studies}.
\newblock Princeton University Press, Princeton, NJ, 1994.

\bibitem{MR2695927}
Mikhail Khovanov.
\newblock {\em Graphical calculus, canonical bases and {K}azhdan-{L}usztig
  theory}.
\newblock ProQuest LLC, Ann Arbor, MI, 1997.
\newblock Thesis (Ph.D.)--Yale University.

\bibitem{MR1680395}
Mikhail Khovanov and Greg Kuperberg.
\newblock Web bases for {${\rm sl}(3)$} are not dual canonical.
\newblock {\em Pacific J. Math.}, 188(1):129--153, 1999.

\bibitem{MR2221529}
Dongseok Kim.
\newblock Trihedron coefficients for {$ U_q({{sl}}_3(\mathbb{ C}))$}.
\newblock {\em J. Knot Theory Ramifications}, 15(4):453--469, 2006.

\bibitem{Kim07}
Dongseok Kim.
\newblock Jones-{W}enzl idempotents for rank 2 simple {L}ie algebras.
\newblock {\em Osaka J. Math.}, 44(3):691--722, 2007.

\bibitem{KirbyMove}
Robion Kirby.
\newblock A calculus for framed links in {$S^{3}$}.
\newblock {\em Invent. Math.}, 45(1):35--56, 1978.

\bibitem{Kupe-first-G2}
Greg Kuperberg.
\newblock The quantum {$G_2$} link invariant.
\newblock {\em Internat. J. Math.}, 5(1):61--85, 1994.

\bibitem{Kupe}
Greg Kuperberg.
\newblock Spiders for rank {$2$} {L}ie algebras.
\newblock {\em Comm. Math. Phys.}, 180(1):109--151, 1996.

\bibitem{KupConjecture}
Greg Kuperberg.
\newblock Personal communication, 2021.

\bibitem{Lickorish}
William B.~R. Lickorish.
\newblock The skein method for three-manifold invariants.
\newblock {\em J. Knot Theory Ramifications}, 2(2):171--194, 1993.

\bibitem{MR1035415}
George Lusztig.
\newblock Canonical bases arising from quantized enveloping algebras.
\newblock {\em J. Amer. Math. Soc.}, 3(2):447--498, 1990.

\bibitem{MR1180036}
George Lusztig.
\newblock Canonical bases in tensor products.
\newblock {\em Proc. Nat. Acad. Sci. U.S.A.}, 89(17):8177--8179, 1992.

\bibitem{LuszQuantum}
George Lusztig.
\newblock {\em Introduction to quantum groups}.
\newblock Modern Birkh\"auser Classics. Birkh\"auser/Springer, New York, 2010.
\newblock Reprint of the 1994 edition.

\bibitem{MR1272656}
Gregor Masbaum and Pierre Vogel.
\newblock {$3$}-valent graphs and the {K}auffman bracket.
\newblock {\em Pacific J. Math.}, 164(2):361--381, 1994.

\bibitem{OhtsukiYamada}
Tomotada Ohtsuki and Shuji Yamada.
\newblock Quantum {${\rm SU}(3)$} invariant of {$3$}-manifolds via linear skein
  theory.
\newblock {\em J. Knot Theory Ramifications}, 6(3):373--404, 1997.

\bibitem{PRVformula}
K.~R. Parthasarathy, R.~Ranga~Rao, and V.~S. Varadarajan.
\newblock Representations of complex semi-simple {L}ie groups and {L}ie
  algebras.
\newblock {\em Ann. of Math. (2)}, 85:383--429, 1967.

\bibitem{MR2679382}
Emily Peters.
\newblock A planar algebra construction of the {H}aagerup subfactor.
\newblock {\em Internat. J. Math.}, 21(8):987--1045, 2010.

\bibitem{RTlink}
N.~Yu. Reshetikhin and V.~G. Turaev.
\newblock Ribbon graphs and their invariants derived from quantum groups.
\newblock {\em Comm. Math. Phys.}, 127(1):1--26, 1990.

\bibitem{RTinv}
Nicolai Reshetikhin and Vladimir~G. Turaev.
\newblock Invariants of {$3$}-manifolds via link polynomials and quantum
  groups.
\newblock {\em Invent. Math.}, 103(3):547--597, 1991.

\bibitem{rowellqmod}
Eric~C. Rowell.
\newblock From quantum groups to unitary modular tensor categories.
\newblock In {\em Representations of algebraic groups, quantum groups, and Lie
  algebras}, volume 413 of {\em Contemp. Math.}, pages 215--230, Providence,
  RI, 2006. Amer. Math. Soc.

\bibitem{G2basecasecalc}
Takuro Sakamoto and Yasuyoshi Yonezawa.
\newblock Link invariant and {$G_2$} web space.
\newblock {\em Hiroshima Math. J.}, 47(1):19--41, 2017.

\bibitem{Wenzl}
Hans Wenzl.
\newblock On sequences of projections.
\newblock {\em C. R. Math. Rep. Acad. Sci. Canada}, 9(1):5--9, 1987.

\bibitem{Weyl1932}
H.~Weyl, G.~Rumer, and E.~Teller.
\newblock Eine für die valenztheorie geeignete basis der binären
  vektorinvarianten.
\newblock {\em Nachrichten von der Gesellschaft der Wissenschaften zu
  Göttingen, Mathematisch-Physikalische Klasse}, 1932:499--504, 1932.

\end{thebibliography}

\newpage
\input appendix.tex

\end{document}